\documentclass[a4paper,10pt,reqno]{amsart}
\usepackage{amsmath,amsthm,amssymb,amscd}
\usepackage{amsfonts}
\usepackage[all]{xy}
\usepackage{latexsym,amssymb}
\usepackage{enumerate}
\usepackage[active]{srcltx}
\renewcommand{\baselinestretch}{1.2}
\newenvironment{beweis}{\vspace{4pt}\noindent \bf Proof\,:\ \rm}{\nolinebreak $\Box$
\par \medskip}
\newcommand{\BB}{\begin{beweis}}
\newcommand{\EB}{\end{beweis}}

\newenvironment{main}{\vspace{4pt}\noindent \bf Proof of Theorem~1.1\,:\ \rm}{\nolinebreak $\Box$
\par \medskip}
\newcommand{\Bmain}{\begin{main}}
\newcommand{\Emain}{\end{main}}

\newtheorem{satz}{Theorem}[section]
\newtheorem{prop}[satz]{Proposition}
\newtheorem{cor}[satz]{Corollary}
\newtheorem{lem}[satz]{Lemma}
\newtheorem{rem}[satz]{Remark}
\newtheorem{defis}[satz]{Definition}
\newtheorem{ex}[satz]{Example}
\newtheorem{exas}[satz]{Examples}

\newcommand{\Bth}[1]{\smallskip \begin{satz} \label{#1} \begin{rm}}
\newcommand{\Eth}{\end{rm} \end{satz}}

 \newcommand{\Bprop}[1]{\smallskip \begin{prop} \label{#1} \begin{rm}}
\newcommand{\Eprop}{\end{rm} \end{prop}}

\newcommand{\Bex}[1]{\smallskip \begin{ex} \label{#1} \begin{rm}}
\newcommand{\Eex}{\nolinebreak $\Box$ \end{rm} \end{ex} \medskip}
\newcommand{\Bexas}[1]{\smallskip \begin{exas} \label{#1} \begin{rm}}
\newcommand{\Eexas}{\nolinebreak $\Box$ \end{rm} \end{exas} \medskip}

\newcommand{\Bcor}[1]{\smallskip \begin{cor} \label{#1} \begin{rm}}
\newcommand{\Ecor}{\end{rm} \end{cor} \medskip}

\newcommand{\Brem}[1]{\smallskip \begin{rem} \label{#1} \begin{rm}}
\newcommand{\Erem}{\end{rm} \end{rem} \medskip}

\newcommand{\Bdef}[1]{\smallskip \begin{defis} \label{#1} \begin{rm}}
\newcommand{\Edef}{\end{rm} \end{defis} \medskip}

\newcommand{\Blem}[1]{\smallskip \begin{lem} \label{#1} \begin{rm}}
\newcommand{\Elem}{\end{rm} \end{lem} \medskip}

\DeclareMathOperator{\Hom}{Hom} \DeclareMathOperator{\Ext}{Ext}
\DeclareMathOperator{\End}{End} \DeclareMathOperator{\rk}{rk}
\DeclareMathOperator{\Add}{Add} \DeclareMathOperator{\Tor}{Tor}
\DeclareMathOperator{\Cogen}{Cogen}

\DeclareMathOperator{\Prod}{Prod} \DeclareMathOperator{\im}{im}
\DeclareMathOperator{\add}{add}

\newcommand{\La}{R}
\newcommand{\m}{\mbox{mod\hspace{0.1em}$\La$}}

\newcommand{\ml}{\mbox{$\La$\hspace{0.1em}mod}}
\newcommand{\M}{\mbox{Mod\hspace{0.1em}$\La$}}
\newcommand{\Modu}{\mbox{Mod\hspace{0.1em}$\La_{\mathcal{U}}$}}
\newcommand{\Ml}{\mbox{$\La$\hspace{0.1em}Mod}}

\newcommand{\eme}{{\mathcal M}}

\newcommand{\cala}{\mbox{$\mathcal A$}}
\newcommand{\calb}{\mbox{$\mathcal B$}}
\newcommand{\calc}{{\mathcal C}}
\newcommand{\calf}{{\mathcal F}}

\newcommand{\cals}{{\mathcal S}}
\newcommand{\calt}{{\mathcal T}}
\newcommand{\calu}{{\mathcal U}}

\newcommand{\cald}{{\mathcal D}}

\newcommand{\calv}{\mathcal{V}}

\newcommand{\p}{\ensuremath{\mathbf{p}}}
\newcommand{\q}{\ensuremath{\mathbf{q}}}
\newcommand{\tube}{\ensuremath{\mathbf{t}}}

\newcommand{\Gen}[1]{\mbox{Gen\,$#1$}}
\newcommand{\ltube}{{}_R\ensuremath{\mathbf{t}}}

\newcommand{\ifa}{if and only if}
\newcommand{\feq}{following statements are equivalent}

\newcommand{\ra}{\rightarrow}

\newcommand{\mapr}[1]{\stackrel{#1}{\longrightarrow}}

\newcommand{\exs}[5]{0\longrightarrow #1 \mapr{#2} #3 \mapr{#4} #5
\longrightarrow 0}
\newcommand{\exst}[3]{0\rightarrow #1\rightarrow #2 \rightarrow #3}
\newcommand{\exsl}[3]{#1\rightarrow #2 \rightarrow #3\rightarrow 0}

\newcommand{\N}{\mathbb N}

\usepackage[pdftex]{hyperref}

\begin{document}

\title{
Tilting modules over tame hereditary algebras}
\author{Lidia Angeleri H\"ugel}
\address{Dipartimento di Informatica - Settore di Matematica\\
Universit\`a di Verona\\
Strada le Grazie 15 - Ca' Vignal 2, 37134 Verona, Italy}
\email{lidia.angeleri@univr.it}
\author{Javier S\'anchez }
\address{
Department of Mathematics - IME \\ University of S\~ao Paulo\\ Caixa
Postal 66281\\ S\~ao Paulo, SP\\ 05314-970, Brazil}
\email{jsanchez@ime.usp.br}

\date{\today}
\maketitle

\begin{abstract}
We give a complete classification of the infinite dimensional tilting modules over a tame hereditary algebra $R$.
We start our investigations by considering
 tilting modules of the form $T=R_\calu\oplus R_\calu /R$ where $\calu$ is a union of tubes, and $R_\calu$ denotes the universal localization of $R$ at $\calu$ in the sense of Schofield and Crawley-Boevey. Here  $R_\calu/R$ is a direct sum of the Pr\"ufer modules corresponding to the tubes in $\calu$.
Over the Kronecker algebra, large tilting  modules are of this form in all but one case, the exception being the Lukas tilting module $L$ whose tilting class $\Gen L$ consists of all modules without indecomposable preprojective summands.
Over an arbitrary tame hereditary algebra, $T$ can have finite dimensional summands, but  the infinite dimensional part  of $T$ is still built up from universal localizations, Pr\"ufer modules and (localizations of) the Lukas tilting module. We also recover the classification of the infinite dimensional cotilting $R$-modules due to Buan and Krause.
\end{abstract}


In this paper, we continue our study of tilting modules arising from universal localization started in \cite{AS}.
More precisely, we consider tilting modules over a ring $R$ that have the form $R_\calu\oplus R_\calu /R$ where $\calu$ is a set of finitely presented $R$-modules of projective dimension one, and $R_\calu$ denotes the universal localization of $R$ at $\calu$ in the sense of Schofield. We have seen in \cite{AS} that over certain rings this construction leads to a classification of all tilting modules. For example, over a Dedekind domain, every tilting module is equivalent to a tilting module of the form $R_\calu\oplus R_\calu /R$ for some set of simple $R$-modules $\calu$. Aim of this paper is to prove a similar result for finite dimensional tame hereditary algebras.

\smallskip

Universal localizations of a tame hereditary algebra $R$ were already investigated by Crawley-Boevey in \cite{CB1}. He showed that the normalized defect provides a rank function $\rho$ as studied by Schofield in \cite{Schofieldbook}, and that the $\rho$-torsion modules are precisely the finite dimensional regular modules. He also described the shape of the universal localization $R_\calu$ at a set $\calu$ of quasi-simple modules, proving that there are substantially different situations depending on whether  $\calu$ does contain a complete clique (that is, all quasi-simples belonging to a certain tube) or not. In particular, $R_\calu$ will be an infinite dimensional $R$-module whenever $\calu$ contains a complete clique.

\smallskip

 We now want to employ these results to give a classification of the large tilting modules over a tame hereditary algebra $R$. By \emph{large} we mean tilting modules $T$ that are not equivalent to finite dimensional ones, that is, there is no finite dimensional tilting module $T'$ such that $\Gen T=\Gen T'$.

\smallskip

Recall that by a result of Bazzoni and Herbera \cite{BH} a
large tilting module $T$ is determined up to equivalence by
a set of finite dimensional modules $\cals$, in the sense
that its tilting class  $\Gen T$ coincides with the class
of modules $X\in\M$ such that $\Ext^1_R(\cals,X)=0$.

\smallskip

The set $\cals$ can be chosen to consist of the finite dimensional
modules in ${}^\perp(T^\perp)$, and then it turns out that
$\cals=\add(\p\cup\tube')$ where $\p$ denotes the class of
indecomposable preprojective modules, and $\tube'\subset\tube$ is a
subset of the class of all finite dimensional indecomposable regular
modules (Theorem \ref{class}).

 \smallskip

Notice that, as a consequence, the lattice of large tilting modules
has a largest and a smallest element. Indeed,  the largest tilting
class in $\M$ not generated by a finite dimensional tilting module
is the class $\p^\perp$ of modules without indecomposable
preprojective summands, which  is generated by the Lukas tilting
module $\mathbf L$ (see \cite{KT} and \ref{Lukas}), while the
smallest one  is the class of all divisible modules $\tube^\perp$,
and the corresponding tilting module is the direct sum $\mathbf
W=\bigoplus_{S\in {\mathbb U}} S[\infty]\oplus  G$  of all Pr\"ufer
modules and the generic module (see \cite{RR} and \ref{Ringel}), or
in other words, it is the tilting module $R_\tube\oplus R_\tube/R$
arising from universal localization at the set of all quasi-simple
modules.

\smallskip

Moreover, from the description of $\cals$ we also deduce that a large tilting module over the Kronecker algebra must have  the form $R_\calu\oplus R_\calu /R$ for some set of quasi-simple $R$-modules $\calu$ in all but one case, the exception being the Lukas tilting module $\mathbf L$ (Corollary \ref{resume}).

\smallskip

In the general case, the situation is more involved due to the
possible presence of finite dimensional summands in $T$ coming from
non-homogeneous tubes. On the other hand,  there are at most
finitely many  such indecomposable summands up to isomorphism
(Lemma~\ref{fginAddT}). This allows to reduce the classification
problem to a situation similar to the Kronecker case. More
precisely, we show that $T$ is equivalent to a tilting module of the
form $Y\oplus M$ where $Y$ is finite dimensional, while  $M$ has no
finite dimensional indecomposable direct summands and is a tilting
module over a suitable universal localization $R'$ of $R$. Since
$R'$ will again be a tame hereditary algebra, this will enable us to
conclude that $M$ is either the Lukas tilting module over $R'$, or
it arises from universal localization at a union of tubes over $R'$.
Notice that the finite dimensional part $Y$ can be described
explicitly. It is a regular multiplicity-free exceptional module
whose indecomposable summands are arranged in disjoint wings, and
the number of summands from each wing equals the number of
quasi-simple modules in that wing. A module satisfying these
properties will be called a {\it branch module}.

\medskip

 Summarizing, we obtain two disjoint families of large tilting modules  as described below.

\smallskip


{\bf Theorem A} (cf.~Theorem \ref{main}) {\it Let $R$ be  a finite dimensional tame hereditary algebra,
and let $\tube=\bigcup_{\lambda\in\mathfrak{T}}\tube_\lambda$ where the
$\tube_\lambda$ are the tubes in the Auslander-Reiten quiver of
$\La$.
\begin{enumerate}[(1)]
\item For every branch module $Y$ there is a tilting module $$T_{(Y,\emptyset)}=Y\oplus (\mathbf L\otimes_R R_{\mathcal U})$$
where $\calu$ is a suitable set of quasi-simple modules determined by $Y$.
\item  For every branch module $Y$ and every non-empty subset $\Lambda\subseteq\mathfrak{T}$ there is a tilting module $$T_{(Y,\Lambda)}=Y\oplus R_{\mathcal V}\oplus
R_{\mathcal V}/R_{\mathcal U}$$ where   $\mathcal U$, $\mathcal V$ are suitable sets of quasi-simple modules determined by $Y$ and $\Lambda$.
\end{enumerate}
 Every large
tilting module is  equivalent  to
 precisely one module
 in this list. }

\medskip

Observe that there are only finitely many branch modules up to
isomorphism (Lemma \ref{fginAddT}). So, if $\mathcal
Y=\{Y_1,\ldots,Y_t\}$ is a complete irredundant set of branch
modules, and $\mathcal P(\mathfrak{T})$ denotes the power set of
$\mathfrak T$, then the  large tilting modules  are  parametrized,
up to equivalence, by the elements of $\mathcal Y\times\mathcal
P(\mathfrak{T})$.

\bigskip

Combining this with  decomposition  results from \cite{R}, we  obtain the following structure result.

\smallskip

{\bf Theorem B.} (cf.~Theorem~\ref{structure} and Corollary~\ref{structure2})
 {\it Let $R$ be  a finite dimensional tame hereditary algebra.
Let $T$ be
a tilting $R$-module  which is not equivalent to a finite dimensional tilting module.
Then $T$  has a unique decomposition
$$T=\bigoplus_{\lambda\in{\mathfrak{T}}} \tube_\lambda(T)\oplus
\overline{T}$$ where  $T'=\bigoplus_{\lambda\in{\mathfrak{T}}} \tube_\lambda(T)$ is a torsion module, hence a direct sum of Pr\"ufer modules and finite dimensional regular modules, and $\overline{T}$ is a torsion-free
module.

More precisely,
for each tube $\tube_\lambda$ of rank $r$, the summand $\tube_\lambda(T)$  is given as follows:
\begin{enumerate}
\item[(i)] if
${}^\perp(T^\perp)$  contains some modules from  $\tube_\lambda$,
but no complete ray, then $\tube_\lambda(T)$ is a branch module
which is a direct sum of at most $r-1$  modules from
$\tube_\lambda$;
\item[(ii)] if ${}^\perp(T^\perp)$ contains some   rays
from $\tube_\lambda$, then $\tube_\lambda(T)$ has precisely $r$ pairwise non-isomorphic indecomposable summands: these are  the $s$
Pr\"ufer modules corresponding to the $s\le r$ rays from $\tube_\lambda$ contained in ${}^\perp(T^\perp)$, and  $r-s$ modules from $\tube_\lambda$;
\item[(iii)] $\tube_\lambda(T)=0$ whenever $\tube_\lambda\cap{}^\perp(T^\perp)=\emptyset$.
\end{enumerate}
Moreover, the torsion-free summand $\overline{T}$ is given as follows:
\begin{enumerate}
\item[(i)]  if
${}^\perp(T^\perp)$  contains  no complete ray, then there is  a set
${\mathcal U}$ of quasi-simple $R$-modules containing no complete
cliques such that $\overline{T}$ is a tilting module over  the
universal  localization $R_{\mathcal U}$  which is equivalent to the
Lukas tilting $R_\calu$-module $\mathbf L\otimes R_{\mathcal U}$;
\item[(ii)]  if
${}^\perp(T^\perp)$  contains  some rays, then there is a set ${\mathcal V}$ of quasi-simple $R$-modules containing  complete cliques such that  $\overline{T}$ is    a projective generator over the universal localization $R_{\mathcal V}$.
\end{enumerate}
}

\bigskip

In particular, we see that a large tilting module $T$ is equivalent to some $T_{(Y,\emptyset)}$ if ${}^\perp(T^\perp)$  contains  no complete ray, and it is  equivalent to some $T_{(Y,\Lambda)}$ with $\Lambda\not=\emptyset$ if ${}^\perp(T^\perp)$ contains  some rays. Indeed, $\Lambda$  consists  of those $\lambda\in\mathfrak{T}$ for which $\tube_\lambda$ has some ray in ${}^\perp(T^\perp)$. Moreover, in the first case the torsion part $T'$ of $T$ coincides with $Y$ up to multiplicities, while in the second case $T'$ also has Pr\"ufer modules as infinite dimensional summands.
In fact, any combination of Pr\"ufer modules $S[\infty]$  can occur in the torsion part as long as the corresponding  quasi-simples $S$  are not regular composition factors of the Auslander-Reiten translate $\tau^-Y$.  Notice furthermore that in both cases the torsion-free part $\overline{T}$ of $T$ is determined by a suitable localization of the Lukas tilting module. For details we refer to Remark~\ref{Lukalization}.

\bigskip

Recall that the large cotilting modules over $R$ have been classified by Buan and Krause in \cite{BK,BK2}, given Bazzoni's result \cite{Bazzonicotilting} that establishes the pure-injectivity of  cotilting modules.
By using the fact that every cotilting module over a finite dimensional algebra is equivalent to the dual of a tilting module \cite{T}, we can
now recover this classification. Let us remark that the other direction does not work: one cannot use the classification of  cotilting modules for studying the tilting modules, as duals of (large) cotilting modules  need not even be tilting, cf.~\ref{Lukas}.

\medskip

The paper is organized as follows. In Section 1, we collect some preliminaries on infinite dimensional modules, tilting theory, and   universal localization. In Section 2, we prove that a large tilting module $T$ is determined  by a set   $\cals=\add(\p\cup\tube')$ as described above, and we settle the cases where $\tube'=\emptyset$ (then $T$ is equivalent to the Lukas tilting module) or $\tube'$ is a union of tubes (then $T$ arises from universal localization). Section 3 is devoted to the finite dimensional summands of $T$. In Section 4, we show that $T$ has a canonical decomposition as above. The description of the torsion-free part $\overline{T}$ is achieved in Section 5, where we also prove our classification and discuss the cases when $T$ is noetherian over its endomorphism ring or ($\Sigma$-)pure-injective.  In the Appendix we deal with the classification of  cotilting modules.

\bigskip

{\bf Acknowledgements:} We thank the referee for careful reading and
for many valuable suggestions.
The first author would like to thank Jan Trlifaj for valuable
discussions that led to the results in  Section 2.

This research was partially carried out during a visit of the second
named author at Universit\`a dell'Insubria, Varese, in 2006,
supported by the Departament d'Universitats, Recerca i Societat de
la Informaci\'o de la Generalitat de Catalunya. The second named
author would like to thank Dipartimento di Informatica e
Comunicazione dell'Universit\`a degli Studi dell'Insubria di Varese
for its hospitality.

The second named author was partially supported by  Funda\c{c}\~ao de Amparo \`a Pesquisa do Estado de S\~ao Paulo (FAPESP) processo n\'umero 09/50886-0.

Furthermore, we acknowledge partial support from   by the DGI and the European
Regional Development Fund, jointly, through Project
 MTM2008--06201--C02--01, and by the Comissionat per Universitats i Recerca
of the Generalitat de Ca\-ta\-lunya, Project 2009 SGR 1389.
\section{Preliminaries}\label{preli}
{\it Throughout this note,} let $\La$ be a finite dimensional
tame hereditary (w.l.o.g.~indecomposable) algebra over a
field $k$. We denote by $\M$
(respectively, \Ml) the category of all right
(respectively, left) $R$-modules and by $\m$ (respectively,
$\ml$) the category of finitely generated right
(respectively, left) $R$-modules.  Let $D:\m\rightarrow
\La\,\mbox{mod}$ be the usual duality. Given a subcategory
$\cals\subset\m$, the subcategory of $\La\,\mbox{mod}$
consisting of the dual modules $D(S)$ with $S\in\cals$ will
be denoted by $\cals^\ast$.

\medskip

We adopt terminology and notation from \cite{R,RR}. In particular,
we denote by $\p,\tube,\q$ the classes
of  indecomposable preprojective, regular, and preinjective
right $R$-modules of finite length, respectively. The corresponding classes of left $R$-modules
will be denoted by
$_{\La}\p,\,_{\La}\tube,\,_{\La}\q$. An arbitrary
  $R$-module  will be called  \emph{regular} if
  it has neither preinjective nor preprojective
direct summands.

\smallskip

We fix a complete
irredundant set of quasi-simple (i.e.~simple regular) modules ${\mathbb U}$, and for
each  $S\in\mathbb U$, we denote by $S[m]$  the
module of regular length $m$ on the ray $$S=S[1]\subset
S[2]\subset \dotsb\subset S[m]\subset S[m+1]\subset\dotsb$$ and let $S[\infty]=\varinjlim S[m]$ be the corresponding   {\em Pr\"ufer
module}. The   {\em adic module}   $S[-\infty]$ corresponding to $S\in\mathbb U$ is defined dually as the inverse limit along the coray ending at $S$.

\smallskip

We
write $\tube=\bigcup_{\lambda\in\mathfrak{T}}\tube_\lambda$, where
 $\tube_\lambda$ denotes the class of indecomposable modules in a tube of the Auslander-Reiten quiver of
$\m$. The tubes in $\ml$
 will be denoted by ${}_R\tube_\lambda$.
 It is well known that almost all tubes are homogeneous, that is, they contain a unique quasi-simple module up to isomorphism. In order to deal with the (at most three) non-homogeneous tubes, we
consider the equivalence relation  $\sim$ on ${\mathbb U}$ generated
by
\begin{quote} $S\sim S'$ if   $\Ext_R^1(S,S')\not=0.$
\end{quote}
According to \cite{CB1}, we call the equivalence  classes
of this relation \emph{cliques}. In other words, two
quasi-simple modules belong to the same clique iff they are
in the same tube. The   order of the clique is the
\emph{rank} of the tube.

\medskip

We will  need a combinatorial description of the {\em extension closure}
 of a set of quasi-simples $\mathcal U\subset {\mathbb U}$, that is, of the smallest
subcategory $\mathcal W\subset\m$ that contains $\mathcal U$ and is closed under extensions. Given a tube
$\tube_\lambda$  of rank $r>1$ and a module $X\in\tube_\lambda$ of
regular length $m<r$, we consider the full subquiver $\mathcal
W_X$ of $\tube_\lambda$ which is isomorphic to the
Auslander-Reiten-quiver  $\Theta(m)$ of the linearly oriented
quiver of type $\mathbb A_m$ with $X$ corresponding to the
projective-injective vertex of $\Theta(m)$. Following
\cite[3.3]{Rbook}, we call $\mathcal W_X$ the \emph{wing} of
$\tube_\lambda$  with \emph{vertex} $X$. The following result is straightforward.

\Blem{wing}  Let $\tube_\lambda$ be a tube of rank $r>1$, and let
$\mathcal U=\{U_1,\ldots,U_m\}\subset \mathbb U$ be a set of $m<r$
quasi-simples in $\tube_\lambda$ where $U_{i+1}=\tau^-U_i$ for all
$1\le i<m$. Then  the extension closure $\mathcal W$ of
$\mathcal U$ consists of all finite direct sums of modules in the wing
 $\mathcal W_{U_1[m]}=\{U_i[k]\,\mid\,1\le i\le m,\,1\le k\le m-i+1\}.$ $\Box$ \Elem

\medskip

Let us introduce some further notation. Let $\eme \subset $
Mod$\La$ be a class of modules.
Denote by {Add}\,{$\eme$ } (respectively,  add\,$\eme$) the class
 consisting of all modules isomorphic to direct summands
of (finite) direct sums of elements of $\eme$. The class consisting
of   all modules isomorphic to direct  summands of products of
modules of ${\mathcal M}$ is denoted by $\Prod{\mathcal M}$. The
class  consisting of the right $R$-modules which are epimorphic
images of arbitrary direct sums of elements in $\mathcal M$ is
denoted by $\Gen{\mathcal{M}}$. Dually, we define $\Cogen{\mathcal
M}$ as the class of all submodules of arbitrary direct products of elements in
$\mathcal{M}$. We further write
$$\eme^o = \{ X_{\La} \mid  {\rm Hom}_\La (M,X) = 0 \;{\rm for\; each}\; M \in
\eme \}$$
$$\eme^\perp = \{ X_{\La} \mid  \Ext^1_\La (M,X) = 0 \;{\rm for\; each}\; M \in
\eme \}$$
$$\eme^{\wedge}=\{ X_{\La} \mid  \Ext^i_\La (M,X) = 0 \;{\rm for\; each}\; i\ge 0,\, M \in
\eme \}$$
$$\eme^\intercal = \{ {}_{\La}X \mid  \Tor_1^\La (M,X) = 0 \;{\rm for\; each}\; M \in
\eme \}$$ and define dually ${}^o\eme,\,{}^\perp
\eme,\,{}^{\wedge}\eme,\, {}^\intercal\eme$. If $\eme$ contains a
unique module $M$, then we shall denote these subcategories by {
Add}{ $M$}, $M^o$, { $M^\perp$}, etc.

\medskip

Finally, we denote by $G$ the \emph{generic module}. It is the
unique indecomposable infinite dimensional module which has finite length over its endomorphism ring. In the notation of
\ref{Ringel} and \ref{Lukas}, it is the unique indecomposable in $\mathcal{F}\cap\mathcal{D}$, that is, the unique
indecomposable torsion-free divisible module, cf.~\cite[5.3 and p.408]{R}.

\smallskip

We now  collect some tools we will freely use when working with infinite dimensional modules.

\Blem{inf}\begin{enumerate}[(1)] \item If $M\in\M$ and $X$
is a finitely generated  indecomposable module in $\Add M$,
then $X$ is isomorphic to a  direct summand  of $M$.

\item Every finite-dimensional $R$-module is {\em
endofinite}, that is, it has finite length as a module over
its endomorphism ring.
 Every direct sum of copies of finitely many endofinite modules  is endofinite. Every dual $D(M)$ of an endofinite module $M$ is endofinite.

\item Suppose $M$ is endofinite. Then $\Add M=\Prod M$. In
particular, if $M\in X^\perp$ for some $X\in\M$, then $
\Add M\subset X^\perp$.

\item
 Every endofinite $R$-module $M$ is {\em pure-injective}, that is, pure-exact sequences starting at $M$ split.

 \item Every indecomposable pure-injective
$R$-module is isomorphic to a module in the
 following list:
\begin{enumerate}[-]
 \item the finitely generated indecomposable modules,
 \item the Pr\"ufer modules $S[\infty]$, $S\in {\mathbb U}$,
 \item the adic modules $S[-\infty]$, $S\in {\mathbb U}$,
 \item the generic module $G$.
\end{enumerate}
\item Let $M$ and $N$ be infinite dimensional indecomposable
pure-injective modules. Then \begin{enumerate}[-]
\item $\Ext_R^1(M,P)\not=0$ for every $P\in \p$.
\item $\Ext_R^1(Q,M)\not=0$ for every $Q\in \q$.
\item $\Ext_R^1(M,N)\not=0$ \ifa\ there are  $S\sim S'$ such that
$M= S[\infty]$ and $N=  S'[-\infty]$.
\end{enumerate}
\end{enumerate}
\Elem \BB (1) Since $X$ is  a finitely generated  module,
 being (isomorphic to) a direct summand of a direct sum $M^{(I)}$  of copies of $M$ means  that $X$ is (isomorphic to) a  direct summand in a finite subsum $M^{(I_0)}$. Now the claim follows from the fact that $X$ has a local endomorphism ring.

 The first statement in (2) is clear because every finite-dimensional $R$-module is  finitely generated over its endomorphism ring, which is again a finite-dimensional $k$-algebra. For the other statements on endofinite modules, we refer to \cite{CB2}.  Details on pure-injective modules can be found in  \cite[Chapter 7]{JL}.
The classification of  the indecomposable
pure-injective $\La$-modules is contained in \cite{CB}. Statement (6) is shown in
 \cite[2.5 and 2.7]{BK}
 \EB

\medskip

Recall from \cite{CT} that a module $T$ is {\em tilting} provided that $\Gen T=T^\perp$, or equivalently, $T$ satisfies\\
(T1) $\mbox{proj.dim}(T) \le 1$;\\
(T2) $\Ext^1_R(T,{T^{(\kappa)}}) = 0$ for any cardinal $\kappa$;\\
(T3) There is an exact sequence $0 \to \La \to T_0 \to T_1 \to 0$
with $T_0,T_1 \in \mbox{Add}(T)$.

Note that every tilting module $T$  satisfies $\Add
T=T^\perp\cap{}^\perp(T^\perp)$. Moreover, $T$ gives  rise to a
torsion pair with torsion class $T^\perp$ and torsion-free class
$T^o$. The class $T^\perp$ is called a {\em tilting class}. Tilting
modules  having the same tilting classes are said to be {\em
equivalent}. {\em Cotilting modules} and {\em cotilting classes} are
defined dually, and  equivalence of cotilting modules is defined
correspondingly.

By \cite[5.1.12]{GT}, two tilting modules $T,T'$ are equivalent if and only if
$\Add T=\Add T'$,  while two cotilting modules $C,C'$ are
equivalent if and only if $\Prod C=\Prod C'$.

\medskip

Here are some examples of infinite-dimensional  tilting or
cotilting modules.

\Bex{Ringel} {\em The Reiten-Ringel tilting module.} It is
shown in \cite{RR} that the module  $$\mathbf
W=\bigoplus_{S\in {\mathbb U}} S[\infty]\oplus  G$$  is an
infinite dimensional tilting module  whose tilting class
$\Gen{\mathbf W}={\mathbf W}^\perp$ coincides with the
class
$$\mathcal D={}^o\tube=\tube^\perp$$  of all {\em
divisible} modules, and moreover, $\mathbf W$ is a
cotilting module whose cotilting class $\Cogen{\mathbf
W}={}^\perp\mathbf W$ coincides with the class $$\mathcal
C={}^\perp \q=\q^o$$ of all modules without indecomposable
preinjective direct summands. \Eex

\Bex{Lukas} {\em The Lukas tilting module.} Based on a
construction due to F.~Lukas \cite[2.1]{L2}, Kerner and
Trlifaj showed  in \cite{KT} that there is a countably
infinitely generated $\p$-filtered tilting module ${\mathbf
L}\in \p^\perp$ whose tilting class   $$\Gen{\mathbf
L}=\p^\perp= {}^o\p$$ coincides with the class of all
modules without indecomposable preprojective direct
summands. The corresponding torsion-free class $\mathbf L^o$
coincides  with the class of preprojective modules in the
sense of \cite[Section~2]{R}. In particular, $\mathbf L^o$
is contained in the class $$\mathcal
F=\tube^o={}^\perp\tube$$ of all {\em torsion-free} modules,
which is a cotilting class with cotilting module $D(_\La
\mathbf W)$, cf. \cite[Prop.7]{aht3}. Here $_\La\mathbf W$
denotes the Reiten-Ringel tilting module in the category
$\La\,\mbox{Mod}$. The torsion class corresponding to the
torsion-free class $\mathcal F$ is the class $\Gen\tube$ of
all {\em torsion} modules.

Note that the dual $D(_\La\mathbf W)$ of the cotilting module $_\La\mathbf W$ is not tilting as it does not satisfy condition (T2). Indeed, $G$ and the adic modules are summands of $D(_\La\mathbf W)$, but no countable  direct sum of  copies of an adic module belongs to $G^\perp$, see \cite[Prop.1 and Remark on p.265]{O}.
\Eex

\medskip

Next, let us recall Schofield's notion of universal localization
{\cite[Theorem~4.1]{Schofieldbook}}.

\Bth{scho} \label{def:universallocalization}
 Let  $\Sigma$ be a set of morphisms between
finitely generated projective right $R$-modules. Then there are a
ring $R_\Sigma$ and a morphism of rings $\lambda\colon
R\rightarrow R_\Sigma$ such that
\begin{enumerate}
\item $\lambda$ is \emph{$\Sigma$-inverting,} i.e. if
$\alpha\colon P\rightarrow Q$ belongs to  $\Sigma$, then
$\alpha\otimes_R 1_{R_\Sigma}\colon P\otimes_R R_\Sigma\rightarrow
Q\otimes_R R_\Sigma$ is an isomorphism of right
$R_\Sigma$-modules, and \item $\lambda$ is \emph{universal
$\Sigma$-inverting}, i.e. if $S$ is a ring such that there exists
a $\Sigma$-inverting morphism $\psi\colon R\rightarrow S$, then
there exists a unique morphism of rings $\bar{\psi}\colon
R_\Sigma\rightarrow S$ such that $\bar{\psi}\lambda=\psi$.
\end{enumerate}
\Eth

The morphism $\lambda\colon R\rightarrow R_\Sigma$ is an epimorphism
in the category of rings  with  ${\rm
Tor}^R_1({R_\Sigma},{R_\Sigma})=0.$ It
 is called the \emph{universal localization of $R$ at
$\Sigma$}.

\smallskip

Let now $\mathcal{U}$ be a set of  finitely presented right
$R$-modules. For each
$U\in\mathcal{U},$ consider a morphism $\alpha_U$ between finitely
generated projective right $R$-modules such that
$$0\to P\stackrel{{\alpha_U}}{\to} Q\to U\to 0$$
 We will denote by
$\lambda_{\mathcal{U}}:R\to R_{\mathcal{U}}$ the universal
localization of $R$ at the set $\Sigma=\{\alpha_U\mid
U\in\mathcal{U}\},$  and we will  call it the \emph{universal
localization of $R$ at ${\mathcal{U}}$}. Note that   $R_{\mathcal
U}$ does not depend on the choice of  $\Sigma$.

\medskip

\Bex{Var} {\em Tilting modules arising from universal
localization.} Let now $\mathcal U\subset\mathbb U$ be a set of
quasi-simple  modules.
 Then, as shown in \cite[4.7]{AS},  the module $$T_{\mathcal U}=\La_{\mathcal U}\oplus \La_{\mathcal U}/\La$$ is a tilting module
with tilting class $\mathcal U^\perp$. In particular, if $\mathcal
U=\mathbb U$, then $T_{\mathcal U}$ is equivalent to the
Reiten-Ringel tilting module $\mathbf W=\bigoplus_{S\in
\mathbb{U}}S[\infty]\oplus G$. \Eex

 More generally, if $\mathcal U$ is a
union of cliques, then $\La_{\mathcal U}$ is a torsion-free module,
and $\La_{\mathcal U}/\La$ is a direct sum of the Pr\"ufer modules
corresponding to the quasi-simples in $\mathcal{U}$, as we are going
to see below in Propositions~\ref{RU} and \ref{prop:sumofprufer}
(compare also \cite[2.4]{T}).

\smallskip

We first collect some facts on
universal localization which we will also need later.
Recall that, given a set of $R$-modules $\calu$, the {\it torsion pair
generated by} $\calu$ is the pair $(\mathcal{T}_\calu,\calu^\circ)$
where $\mathcal{T}_\calu={}^\circ(\calu^\circ)$.

\Bprop{univloc} Let $\mathcal U $ be a set of quasi-simple modules
and let $\mathcal W$ be the extension closure of $\mathcal U$. Let
further $t$ be the torsion
radical associated to the torsion pair $(\calt_{\mathcal U}, \mathcal U^o)$
generated by $\mathcal U$. The following statements hold true.
\begin{enumerate}[(1)]
\item $\mathcal W$  is a full exact abelian subcategory of \m.
\item $R_{\mathcal U}$ coincides with $R_\mathcal{W}$, the universal localization of $R$ at $\mathcal W$.
\item The torsion pair $(\calt_{\mathcal U}, \mathcal U^o)$ generated by $\mathcal U$ coincides with the torsion pair $(\mathcal{T}_\mathcal{W},\mathcal{W}^\circ)$ generated
by $\mathcal W$.
\item $\mathcal U^\wedge=\mathcal W^\wedge$ is the essential image of the restriction functor $\M_\calu\to\M$. In other words, an $R$-module $X$ is an $R_{\mathcal U}$-module if and only if $X\in\mathcal U^\wedge$.
\item $\calt_{\mathcal U}=\Gen {\mathcal W}=\{X\in\M\,\mid\, X\otimes_R R_\calu=0\}$.
\item  $R_\calu/R$ is a directed union of
finite extensions of modules in $\mathcal{U}$.
\item For every $A\in\M$ there is a short exact sequence
$$0\to A/tA\to A\otimes_R R_\calu\to A\otimes_R R_\calu/R\to 0$$
where $A\otimes_R R_\calu\in\mathcal U^\wedge$ and $A\otimes_R
R_\calu/R\in{}^\wedge(\mathcal U^\wedge)=\calt_{\mathcal
U}\cap{}^\perp(\mathcal U^\perp)$.
\end{enumerate}
\Eprop
\BB (1), (2) We adopt Schofield's terminology from
\cite{S2}. Since $\mathcal U$ is a Hom-perpendicular set,
$\mathcal W$ is well-placed, cf.~\cite[p.4]{S2}. Then
$\mathcal W={}^\wedge(\mathcal U^\wedge)\cap\m$ is the
well-placed closure of $\mathcal U$, and $R_{\mathcal
U}=R_{\mathcal W}$, cf.~\cite[2.3]{S2}.

(3), (4) We claim  $\mathcal U^o=\mathcal W^o$. The inclusion
`$\supset$' follows from $\mathcal U\subset\mathcal W$.
Conversely,  ${}^o(\mathcal U^o)$ contains $\mathcal U$, and also
its extension closure $\mathcal W$, hence $\mathcal
U^o\subset\mathcal W^o$.
 Similarly, we prove $\mathcal U^\perp=\mathcal W^\perp$. We then deduce $\mathcal U^\wedge=\mathcal W^\wedge$. For the second statement see \cite[1.7]{AA}.

(5) $\{X\in\M\,\mid\, X\otimes_R R_\calu=0\}$ is closed under
extensions, direct sums and epimorphic images, hence it is a torsion
class containing $\mathcal U$ and thus also $\calt_{\mathcal U}$,
which in turn contains $\Gen{ \mathcal{W}}$. The converse inclusions
follow from \cite[5.1 and 5.5]{S1}.

(6) is a consequence of
\cite[Theorem~12.6]{Schofieldbook},
\cite[Theorem~3]{Schofieldquivers} and \cite[Lemma~4.4]{CB1}.
Another proof can be found in \cite[Theorem~2.6]{S2}.

(7) is contained in  \cite[page 2349]{KS}. We give a direct proof for
the reader's convenience. Applying $A\otimes_R -$ on the short exact
sequence $0\to R\to R_\calu\to R_\calu/R\to 0$, we obtain an exact
sequence $A\to A\otimes_R R_\calu\to A\otimes_R R_\calu/R\to 0$,
which gives rise to the short exact sequence  $0\to A/tA\to
A\otimes_R R_\calu\to A\otimes_R R_\calu/R\to 0$ because $tA$ is the
kernel of the canonical map $A\to A\otimes_R R_\calu$, cf.~(5) and
\cite[5.5]{S1}. Since $A\otimes_RR_\calu$ is an $R_\calu$-module, it
follows from (4)  that $A\otimes_R R_\calu\in\mathcal U^\wedge$.

Let us show that $A\otimes_R R_\calu/R\in{}^\wedge(\mathcal
U^\wedge)$. Take $M\in\calu^\wedge$. First of all, note that
$A\otimes_R R_\calu/R$ is generated by $R_\calu/R$, which belongs to
$\calt_{\mathcal U}$ by (5) and (6). Thus $A\otimes_R
R_\calu/R\in\calt_{\mathcal U}$, and since $M\in\calu^o$, we have
$\Hom_R(A\otimes R_\calu/R,M)= 0$.

Next, note that $A/tA\otimes_RR_\calu\cong A\otimes_R R_\calu$ by
(5). Then $\Hom_R(A\otimes_RR_\calu,M)\cong\Hom_R(A/tA,M)$ because
$M$ is an $R_\calu$-module by (4), and we have the exact sequence
$${0\rightarrow
\Ext_R^1(A\otimes_RR_\calu/R,M)\rightarrow\Ext_R^1(A\otimes_RR_\calu,M)\stackrel{\psi}{\longrightarrow}\Ext_R^1(A/tA,M)\rightarrow
0}$$ Now, if we prove that $\psi$ is injective, we obtain that
$\Ext_R^1(A\otimes_RR_\calu/R,M)=0$, as desired.

Given an extension
$$\epsilon:\exs {M}{}{X}{}{A\otimes_RR_\calu},$$ its image under $\psi$ is given by pullback
$$\xymatrix{\epsilon:&0\ar[r]&M \ar[r]\ar@{=}[d]& X \ar[r] & A\otimes_R
R_\calu\ar[r]& 0 \\
\psi(\epsilon):&0\ar[r]& M\ar[r] & Z\ar@{^{(}->}[u]\ar[r] &
A/tA\ar[r]\ar@{^{(}->}[u] & 0}$$ Observe that $X$ is an
$R_\calu$-module because $\calu^\wedge$ is closed under extensions.
Therefore we obtain the commutative diagram in $\M_\calu$
$$\xymatrix{ 0\ar[r] & M \ar[r]\ar@{=}[d] &
X\ar[r] & A\otimes_RR_\calu\ar[r] & 0\\
0\ar[r] & M\ar[r] & Z\otimes_RR_\calu\ar[u]\ar[r] &
A/tA\otimes_RR_\calu\ar[r]\ar[u]^\cong & 0}$$ showing that
$\epsilon$ is uniquely determined by $\psi(\epsilon)$. Thus $\psi$
is injective.

Finally, we show that ${}^\wedge(\mathcal U^\wedge)=\calt_{\mathcal
U}\cap{}^\perp(\mathcal U^\perp)$. The inclusion `$\supset$' follows
from the definition. For `$\subset$', consider
$X\in{}^\wedge(\mathcal U^\wedge)$. If $A\in\mathcal U^o$, then
there is an embedding $0\to A\to A\otimes_R R_\calu\in\mathcal
U^\wedge$, hence $0\to{\rm Hom}_R(X,A)\to{\rm
Hom}_R(X,A\otimes_RR_\calu)=0$, which proves $X\in{}^o(\mathcal
U^o)=\calt_{\mathcal U}$. Moreover, if $A\in\mathcal U^\perp$, then
by \cite[4.7]{AS} there is an epimorphism $R_\calu^{(\alpha)}\to
A\to 0$ where $R_\calu^{(\alpha)}\in\mathcal U^\wedge$, thus
$0=\Ext^1_R(X,R_\calu^{(\alpha)})\to\Ext^1_R(X,A)\to 0$, showing
$X\in{}^\perp(\mathcal U^\perp)$. \EB

\Bprop{RU} Let $\mathcal{U}\subseteq\mathbb{U}$ be a set of
quasi-simple modules.
The following statements hold true.
\begin{enumerate}[(1)]
\item The $R$-module $R_\calu$ is torsion-free, and the $R$-module $R_\calu/R$ is torsion regular.
\item The  $R$-module $R_\mathbb{U}$ is torsion-free and divisible, and it is a direct sum of $\alpha={\rm dim}_{\End_RG}\,G
$ copies of  the generic module $G$.
Moreover, $R_\mathbb{U}$ is a simple artinian ring isomorphic to the
ring of \mbox{$\alpha\times\alpha$}--matrices over the division ring
$\End_RG$,  and $G$ is  the unique indecomposable
$R_\mathbb{U}$-module.
\item For any module $V$ in the extension closure of $\calu$ there is an isomorphism of $k$-$\End_R V$-bimodules $$\Hom_R(V,R_\calu/R) \cong
\Ext_R^1(V,R).$$
\end{enumerate}
\Eprop
\BB
(1) and (2): Let $\calu\subseteq\mathbb{U}$. First of all, we show that $R_\calu/R$
is a torsion regular module, that is, it belongs to $\Gen \tube$ and has no summands in $\p\cup\q$. By Proposition~\ref{univloc}(6),   we can write $R_\calu/R$ as a directed union $\varinjlim N_i$ with the $N_i$'s
finite extensions of elements in $\calu$.  Then
$R_\calu/R\in \Gen \tube$. Moreover, if $P\in\p$
(respectively, $Q\in\q)$ were a direct summand of $R_\calu/R$, then $P$ (respectively, $Q)$
would be a direct summand of some regular module $N_i$, a contradiction.

 If  $\calu=\mathbb{U}$, then the fact that $R_\mathbb U$ is an
$R_\mathbb U$-module yields by Proposition~\ref{univloc}(4) that
$\Hom_R(\mathbb U,R_\mathbb U)=\Ext_R^1(\mathbb U,R_\mathbb U)=0$,
that is, $R_\mathbb U$ is a torsion-free divisible module. So
\cite[5.4 and 5.6]{R} imply that $R_\mathbb U$ is a direct sum of
$-\delta(R)$ copies of $G$, where  $\delta$  denote the defect (cf.\cite[p.333]{R}).  The ring $R_\mathbb{U}$ is simple
artinian by \cite[Lemma~4.4]{CB1}. The $R$-module $G$ is an
$R_\mathbb{U}$-module because it is a torsion-free divisible
$R$-module, and it is the only simple $R_\mathbb{U}$-module because
it is indecomposable over $R$ and hence over $R_\mathbb{U}$. Now
$R_\mathbb{U}\cong \End_{R_\mathbb{U}} R_\mathbb{U}\cong \End_R
R_\mathbb{U}\cong \End_R(G^{(\alpha)})\cong M_\alpha (\End_RG)$, and
$\End_RG$ is a division ring by \cite[5.3]{R}. Finally, $\alpha={\rm
dim}_{\End_RG}\,G$ by the Theorem of Wedderburn-Artin.

In the general case, we have that $R_\calu$ is torsion-free because
$R_\calu\subset R_\mathbb{U}$ by Propostion~\ref{useful}(4).

(3) follows from Proposition~\ref{univloc}(4) by applying
$\Hom_R(V,-)$ to the exact sequence $\exs R {} {R_\calu} {}
{R_\calu/R}$. \EB

\Blem{samedimension} Let $\tube_\lambda$ be a tube of rank $r$. If $X$ and $Y$
are indecomposable regular modules in $\tube_\lambda$ of regular
length at most $r$, then $\End_R X$ and $\End_R Y$ are isomorphic
division rings. Moreover, \begin{enumerate}[(1)] \item if
$X\subseteq Y$, then $f(X)\subseteq X$ for all $f\in \End_RY$ and
the map $\End_R Y\rightarrow\End_R X$, given by $f\mapsto f_{|X}$ is
an isomorphism.
\item if $X\cong Y/K$ for some (unique regular) $R$-submodule $K$ of $Y$, then the map
$\End_R Y\rightarrow \End_R X$ given by $f\mapsto \bar{f}$ is an
isomorphism where $\bar{f}$ is the induced map on $Y/K$ by $f$.
\end{enumerate}
\Elem \BB Let $\{U_1,\dotsc,U_r\}\subseteq\mathbb{U}$ be the set of
$r$ quasi-simples in $\tube_\lambda$ where $U_{i+1}=\tau^-U_i$ for
all $1\leq i<r$. That $\End_R(U_i[j])$ is a division ring for $1\leq
i\leq r,\ 1\leq j\leq r$ follows from the fact that $\tube$ is an
abelian category and that every indecomposable regular module has
unique regular composition series. By the same reason, the maps in (1) and (2) are well-defined injective morphisms of $k$-algebras.

Fix $i\in \{1,\dotsc,r\}$ and $1\le s<r$. Then the exact sequence
$$\exs{U_i}{}{U_i[s+1]}{}{U_{i+1}[s]}$$
induces the following ones
$$\exst{\Hom_R(U_i,U_i)}{\Hom_R(U_i,U_i[s+1])}{\Hom_R(U_i,U_{i+1}[s])=0}$$
$$\exsl{0=\Hom_R(U_{i+1}[s],U_i[s+1])}{\Hom_R(U_i[s+1],U_i[s+1])}{\Hom_R(U_i, U_i[s+1])}$$
$$\exst{\Hom_R(U_{i+1}[s],U_{i+1}[s])}{\Hom_R(U_i[s+1],U_{i+1}[s])}{\Hom_R(U_i,U_{i+1}[s])=0}$$
$$\exsl{0=\Hom_R(U_i[s+1],U_i)}{\Hom_R(U_i[s+1],U_i[s+1])}{\Hom_R(U_i[s+1],U_{i+1}[s])}$$

Hence $\dim_k \End_R(U_i)=\dim_k \End_R(U_i[s+1])=\dim_k\End_R(U_{i+1}[s]).$ Since we have the injective morphisms of rings $\End_R(U_i[s+1])\rightarrow \End(U_i)$, $f\mapsto f_{|U_i}$, and $\End_R(U_i[s+1])\rightarrow\End(U_{i+1}[s])$, $f\mapsto \bar{f}$, it turns out that $\End_R(U_i)$, $\End_k(U_i[s+1])$ and $\End_R(U_{i+1}[s])$ are isomorphic $k$-algebras for any $1\leq s< r$. The result now follows because $i\in\{1,\dotsc,r\}$ is arbitrary.
\EB

\Bprop{prop:sumofprufer} Let $\mathcal{U}\subseteq\mathbb{U}$ be a
set of quasi-simple modules. Set $\alpha={\rm dim}_{\End_RG}\,G$ and $\alpha_U=\dim_{\End_RU} \Ext_R^1(U,R)$
for each $U\in\mathcal{U}$. The following statements hold true.
\begin{enumerate}[(1)]
\item If $\mathcal{U}$ is a union of cliques, then
$$R_\mathcal{U}/R\cong\bigoplus_{U\in\mathcal{U}}
U[\infty]^{(\alpha_U)}.$$ In particular, if
$\calu=\mathbb{U}$, then
$$T_\mathcal{U}=G^{(\alpha)}\oplus(\bigoplus_{U\in\mathbb{U}}U[\infty]^{(\alpha_U)}).$$
\item Let $\tube_\lambda$ be a tube of
rank $r>1$, let $\mathcal{U}=\{U_1,\dotsc,U_m\}\subseteq\mathbb{U}$
be a set of $m<r$ quasi-simples in $\tube_\lambda$ where
$U_{i+1}=\tau^-U_i$ for all $1\leq i<m$. Then $R_\mathcal{U}/R$ is a direct sum of modules on the coray ending at $U_m$. More precisely,
$$R_\mathcal{U}/R\cong U_1[m]^{(\alpha_{U_1})}\oplus U_2[m-1]^{(\alpha_{U_2})}\oplus\dotsb\oplus U_m^{(\alpha_{U_m})}.$$
 \end{enumerate}
\Eprop

\BB First of all,  by Proposition~\ref{univloc}(6), we can write $R_\calu/R$ as a directed union $\bigcup N_i=\varinjlim N_i$ with the $N_i$'s
finite extensions of elements in $\calu$.

(1) Suppose  that $\calu$ is a union of cliques. Then
$R_\calu/R$ is divisible. Indeed, if $V$ is a quasi-simple not in $\calu$,
then $\Ext_R^1(V,R_\calu/R)=\varinjlim\Ext_R^1(V,N_i)=0$. On the
other hand, if $U\in\calu$, then $\Ext_R^1(U,R_\calu/R)=0$ because $T_\calu=R_\calu\oplus R_\calu/R$ is a tilting module with tilting
class $\calu^\perp$.
So $R_\calu/R$ is a divisible torsion regular
module, hence a direct sum of Pr\"ufer modules by \cite[4.5, Lemma~3]{R}.

Observe that
 for $U,V\in\mathbb U$ we have $$\Ext^1_R(U[\infty],\tau V)\cong D\Hom_R(V, U[\infty])=0\;\text{iff}\; V\not= U.$$
So, if
$V\in\mathbb{U}\setminus\calu$, then as $\calu$ is a union of cliques,
$\tau V\in\calu^\perp= (R_\calu/R)^\perp$, which implies that the Pr\"ufer module $V[\infty]$ cannot occur
as a direct summand of $R_\calu/R$. Similarly, if $U\in\calu$, then clearly $\tau U\notin\calu^\perp=(R_\calu/R)^\perp$,
thus the Pr\"ufer module $U[\infty]$
must be a direct summand of $R_\calu/R$.  Therefore
$$R_\calu/R\cong\bigoplus_{U\in\calu}U[\infty]^{(\beta_U)}$$
for some cardinals $\beta_U$. Recall now that  $\End_R(U)$ is a
division ring for $U\in\mathbb{U}$. Furthermore,
$\dim_{\End_R(U)}\Hom_R(U,U[n])=1$ for all $n\ge 1$, and thus
$\dim_{\End_R(U)}\Hom_R(U,U[\infty])=1$. Then, for a fixed
$U\in\calu$, the number of direct summands of $R_\calu/R$ isomorphic
to $U[\infty]$ coincides with $\dim_{\End_R(U)}\Hom_R(U,R_\calu/R)$,
which by Proposition~\ref{RU}(3) equals $\alpha_U$. Therefore
$\beta_U=\alpha_U$ for all $U\in\calu$, as desired.

The statement for $\calu=\mathbb U$ follows from Proposition \ref{RU}(2).

\smallskip

(2) Suppose now that $\calu$ is defined as in (2).
Note that  the modules $N_i$ above are
finite direct sums of elements in the wing $\mathcal
W_{U_1[m]}$ of $\calu$, see Lemma~\ref{wing}. Set $Y=U_1\oplus\dotsb\oplus
U_m$. By Proposition \ref{RU}(3) we get that
$\Hom_R(Y,R_\calu/R)\cong\Ext_R^1(Y,R)$, which implies that
$\dim_k\Hom_R(Y,R_\calu/R)<\infty$. Therefore the directed union
$R_\calu/R=\bigcup N_i$ is finite, which means that
$R_\calu/R$ is a finite direct sum of elements in $\mathcal W_{U_1[m]}$. In
particular $R_\calu/R$ and $R_\calu$ are finite dimensional over
$k$ (this  is well known, see \cite[Theorem~4.2]{CB1} and
\cite[Theorem~13]{Schofieldquivers} or \cite[10.1]{GL}).

Since the number of direct summands of $R_\mathcal{U}/R$ isomorphic to some module in the ray determined by $U_i$ equals  $\dim_{\End_R(U_i)}\Hom_R(U_i,R_{\mathcal{U}}/R)$,
the total number of indecomposable direct summands of $R_\calu/R$
equals $\alpha_{U_1}+\dotsb+\alpha_{U_m}$ by
Proposition~\ref{RU}(3).

Let us consider the modules in the wing $\mathcal
W_{U_1[m]}$ that lie on the coray ending at $U_m\in\calu$. These are $\{ U_m, U_{m-1}[2],\ldots, U_2[m-1], U_1[m]\}=\{U_{m-i+1}[i]\mid i=1,\dotsc,m\}$.

For a fixed $1\le i\le m$, the number of direct summands of $R_\mathcal{U}/R$ admitting a non-zero morphism from $U_{m-i+1}[i]$ equals $\dim_{\End_R(U_{m-i+1}[i])} \Hom_R(U_{m-i+1}[i],R_{\mathcal{U}}/R)$. If $i=1$, this number agrees with  $\alpha_{U_m}$ by Proposition \ref{RU}(3).  This shows that $R_\mathcal{U}/R$ has $\alpha_{U_m}$ summands isomorphic to $U_m$.

For $i\geq2$, we observe that any morphism from $U_{m-i+1}[i]$ to $R_\mathcal{U}/R$ which is not injective factors through $U_{m-i+2}[i-1]$. Thus the number of direct summands of $R_\mathcal{U}/R$ which are isomorphic to $U_{m-i+1}[i]$ equals
$$\dim_{\End_R(U_{m-i+1}[i])} \Hom_R(U_{m-i+1}[i], R_\mathcal{U}/R) -\dim_{\End_R(U_{m-i+2}[i-1])}\Hom_R(U_{m-i+2}[i-1],R_\mathcal{U}/R).$$
We want to show that this number agrees with  $\alpha_{U_{m-i+1}}$. To this end, we claim that $$\dim_{\End_R(U_{m-i+1}[i])}\Hom_R(U_{m-i+1}[i], R_\mathcal{U}/R)=\alpha_{U_m}+\dotsb+\alpha_{U_{m-i+1}}$$ for $i=1,\dotsc,m$.
This is clear for $i=1$.
We proceed by recurrence and suppose our claim holds true for $i$.
From the exact sequence $\exs{U_{m-i}}{}{U_{m-i}[i+1]}{}{U_{m-i+1}[i]}$ we obtain
the exact sequence $$\exs {\Hom_R(U_{m-i+1}[i],R_\mathcal{U}/R)}{}{\Hom_R(U_{m-i}[i+1],R_\mathcal{U}/R)}{}{\Hom_R(U_{m-i},R_\mathcal{U}/R)},$$
hence
$$
\dim_k\Hom_R(U_{m-i}[i+1],R_\calu/R)=\dim_k\Hom_R(U_{m-i+1}[i],R_\calu/R)+\dim_k\Hom_R(U_{m-i},R_\calu/R).
$$
But for every indecomposable module $X\in \mathcal{W}_{U_1[m]}$, we can compute
$\dim_k\Hom_R(X,R_\mathcal{U}/R)= \gamma\cdot \dim_{\End_R(X)}\Hom_R(X,R_\mathcal{U}/R)$ where $\gamma=\dim_k\End_R(X)$ does not depend on $X$ by Lemma~\ref{samedimension}. Therefore,
dividing by $\gamma$, and using that $\dim_{\End_R(U_{m-i})}\Hom_R(U_{m-i},R_\calu/R)=\alpha_{U_{m-i}}$,
 we obtain the claim for $i+1$.

 So $R_\mathcal{U}/R$ has   $\alpha_{U_{m-i+1}}$ summands isomorphic to $U_{m-i+1}[i]$ for any $1\le i\le m$, and (2) is proven.
\EB

Here are some further results on universal localization that will be needed in Sections 5 and 6.

\Bprop{useful}{\bf \cite{CB1,Schofieldquivers,GL}}
\begin{enumerate}[(1)]
\item Let $\mathcal{Y}$ be a set of quasi-simple modules, and set $\tau\mathcal{Y}=\{\tau V\mid
V\in\mathcal{Y}\}$.
\begin{enumerate}[(a)]
\item  If $S\in\mathbb{U}\setminus(\mathcal{Y}\cup\tau\mathcal{Y})$, then
$S\otimes_RR_\mathcal{Y}\cong S$.

\item If $S\in\tau\mathcal{Y}\setminus\mathcal{Y}$, then $S$ belongs
to a tube $\tube_\lambda$ of rank $r>1$, and there exists  $2\le m\le r$ such
that
 $S\otimes_RR_\mathcal{Y}\cong S[m]$.
\end{enumerate}

\item Assume that $\mathcal{U}\subset\mathbb U$ is a set of
quasi-simple $R$-modules that contains no complete cliques. Then:
\begin{enumerate}[(a)]
\item The universal localization $R_\mathcal{U}$ is a tame hereditary
$k$-algebra with $\rk K_0(R_\mathcal{U})=\rk K_0(R)-|\mathcal{U}|$.
\item The set $\{S\otimes_RR_\mathcal{U}\mid S\in\mathbb{U}\setminus
\mathcal{U}\}$ is a complete irredundant set of quasi-simple
$R_\mathcal{U}$-modules.

\item The set
$\mathbf{\tube}_\calu=\{V\otimes_R R_\calu\mid V \in\mathbf{\tube}
\textrm{ with }  \Hom_R(V,U)=\Hom_R(U,V)=0\ \textrm{ for all }
U\in\mathcal{U}\}$ is a complete irredundant set of  finite
dimensional indecomposable regular $R_\mathcal{U}$-modules.

\item In particular, if $\tube_\lambda$ is a tube of rank $r>1$ with
quasi-simples $U_1,\, U_2=\tau^-U_1, \dotsc, U_{r}=\tau^-U_{r-1}$,
and  $\calu=\{U_2,\dotsc,U_{m+1}\}$ for some $m<r$, then the tube
$\tube_\lambda\otimes R_\calu$ in the Auslander-Reiten quiver of
$R_\calu$ is given by the quasi-simple $R_\calu$-modules
$$U_1\otimes_RR_\calu,\; \tau^-(U_1\otimes
R_\calu)=U_{m+2}\otimes_RR_\calu,\;\dotsc\;
,\tau^-(U_{r-1}\otimes_RR_\calu)=U_r\otimes_RR_\calu.$$

\item The set $\{S[\infty]\mid S\in\mathbb{U}\setminus\calu\},$ is a complete
irredundant set of Pr\"ufer $R_\calu$-modules. We have
$(S\otimes_RR_\calu)[\infty]=S[\infty]$ for each
$S\in\mathbb{U}\setminus\calu$.
\end{enumerate}

\item Assume that $\mathcal{V}\subset\mathbb U$ is a set of
quasi-simple $R$-modules that contains  a complete clique. Then
$R_\calv$ is a hereditary order. Moreover, $\{S\otimes_RR_\calv\mid
S\in \mathbb{U}\setminus\calv\}$ is a  complete irredundant set of
simple $R_\calv$-modules, and  $\{S[\infty]\mid
S\in\mathbb{U}\setminus\calv\}$,
 is a complete irredundant set of injective envelopes of
simple $R_\calv$-modules. We have injective envelopes
$E(S\otimes_RR_\calv)=S[\infty]$ for each
$S\in\mathbb{U}\setminus\calv$.

\item
Assume that $\mathcal{U}\subset\mathbb U$ and $\mathcal{V}\subseteq \mathbb{U}\setminus\calu$. Then
$R_{\calu\cup\mathcal{V}}=(R_\calu)_{\mathcal{V}'}$ where
$\mathcal{V}'=\{V\otimes_RR_\calu\mid V\in\mathcal{V}\}$.
In particular, there is an injective ring epimorphism $R_{\mathcal U}\to R_{\mathbb U}$.
\end{enumerate}
\Eprop

\BB (1) (a) If
$S\in\mathbb{U}\setminus(\mathcal{Y}\cup\tau\mathcal{Y})$, then
$\Hom_R(V,S)=\Ext_R^1(V,S)=\Hom_R(S,\tau V)=0$ for all
$V\in\mathcal{Y}$. That is, $S$ is an $R_\mathcal{Y}$-module by
Proposition~\ref{univloc}(4), and therefore
$S\otimes_RR_\mathcal{Y}=S$.

(b) Let  $S\in\tau\mathcal{Y}\setminus\mathcal{Y}$. Then $S$ belongs
to   a tube $\tube_\lambda$ of rank $r>1$. Choose the numbering
$S=U_1,\, U_2=\tau^-U_1, \ldots, U_r=\tau^-U_{r-1}$ for the
quasi-simples in $\tube_\lambda$. Since, by assumption,
$S\notin\mathcal{Y}$, there is $m$ with $2\le m\le r$ such that
$U_2,\ldots, U_m\in\mathcal Y$ and $\tau^-U_{m}\notin\mathcal Y$.

For each $p=1,\dotsc,m-1$, the exact sequence $0\to
S[p]\to S[p+1]\to U_{p+1}\to 0$ induces
$$\dotsb\to\Tor_1^R(U_{p+1},R_\mathcal{Y})\to S[p]\otimes_R R_\mathcal{Y}\to S[p+1]\otimes_RR_\mathcal{Y}\to U_{p+1}\otimes_RR_\mathcal{Y}\to 0.$$
Clearly
$\Tor_1^R(U_{p+1},R_\mathcal{Y})=U_{p+1}\otimes_RR_\mathcal{Y}=0$
as $U_{p+1}\in\mathcal{Y}$.
Hence we obtain that $S\otimes_RR_\mathcal{Y}\cong S[2]\otimes_RR_\mathcal{Y}\cong\ldots\cong
S[m]\otimes_RR_\mathcal{Y}$ .
Note that $S[m]$ is an $R_\mathcal{Y}$-module because
$\Hom_R(\mathcal{Y},S[m])=0$ and $\Ext_R^1(\mathcal{Y},S[m])\cong D\Hom_R(S[m],\tau\mathcal{Y})=0$. Thus
$S[m]\otimes_R R_\mathcal{Y}\cong S[m]$ as
desired.

(2) Statement (a) is  \cite[Theorem~4.2(1)]{CB1}. The shape of the
quasi-simple and the regular $R_\calu$-modules  follows from
\cite[Theorem~10]{Schofieldquivers} (cf. \cite[Theorem~3.5]{S2}), as
noted in \cite[2.3, 2.4, Section~4]{CB1}. See also \cite[10.1]{GL}.
The statement on $\tube_\lambda\otimes R_\calu$ is shown in
\cite[Section~4.2]{CB1}.

It remains to prove (e). Let $S$ be a quasi-simple $R$-module not in
$\calu$. Since Pr\"ufer modules are divisible and $\Hom_R(\calu,
S[\infty])=0$, it follows that $S[\infty]\in\calu^{\wedge}$ is a right
$R_\calu$-module by
Proposition~\ref{univloc}(4). Further, if $S$ belongs to the $R$-tube
$\tube_\lambda$, $S[\infty]$ is filtered by the quasi-simple
$R_\calu$-modules $\{S\otimes_RR_\mathcal{U}\mid
S\in\tube_\lambda\setminus \mathcal{U}\}$ by (1). Then
$\{(S\otimes_RR_\calu)[n]\,\mid\, n\in\N\}$ is a ray on the
$R_\calu$-tube $\tube_\lambda\otimes R_\calu$, and
$S[\infty]=\varinjlim\limits_n (S\otimes_RR_\calu)[n]$.

(3) By \cite[4.2]{CB1}, $R_\calv$ is a hereditary order, and by
\cite[Section~3]{CB1} (or \cite[6.5]{AS}) and
\cite[Theorem~10]{Schofieldquivers}, the set
$\{S\otimes_RR_\calv\mid S\in\mathbb{U}\setminus\calv \}$ is an
irredundant set of simple $R_\calv$-modules.

Let $S\in\mathbb{U}\setminus\calv$, and suppose that $S$ belongs to
the $R$-tube $\tube_\lambda$. By (1), $S[\infty]$ is filtered by the
simple $R_\calv$-modules $U\otimes_RR_\calv$, where  $U$ runs
through the quasi-simple modules in $\tube_\lambda\setminus\calv$.
By \cite[Theorem~10]{Schofieldquivers}, there exists an equivalence
of categories from the category of bound $R$-modules $M$ such that
\begin{eqnarray}\label{eq:equivalence}
\Hom_R(M,V)=\Hom_R(V,M)=0 \textrm{ for all } V\in\calv
\end{eqnarray}
to the category of bound $R_\calv$-modules that restricts to an equivalence from the category of regular $R$-modules satysfying \eqref{eq:equivalence} to the category of torsion $R_\calv$-modules.
Thus $S[\infty]$ is a uniserial
$R_\calv$-module that contains $S\otimes_RR_\calv$, and the injective envelope $E(S\otimes_RR_\calv)$ of
$S\otimes_RR_\calv$ has to contain $S[\infty]$. But by
\cite[Theorem~19(c)]{GoodearlWarfield},
$E(S\otimes_RR_\calv)$ is also uniserial and has the same
filtration as $S[\infty]$, so they must coincide.

(4) is shown in
\cite[2.4]{CB1} as a consequence of \cite{Schofieldbook} and \cite{Schofieldquivers}.
\EB

\bigskip

\section{Parametrizing tilting modules.}

Tilting classes  are in one-to-one-correspondence with certain
subcategories of \m. Recall that a subcategory $\cals\subset\m$ is
said to be {\it resolving} provided $\cals$ is closed under direct
summands,  extensions, and kernels of epimorphisms, and $\La$
belongs to $\cals$. Observe that, since $\La$ is hereditary,   a
subcategory $\cals\subset\m$ is resolving whenever it is closed
under direct summands and extensions and
 contains $\La$, see \cite[1.1]{aht}.

\smallskip

Bazzoni and Herbera proved in  \cite{BH} that every tilting class
$\mathcal B=T^\perp$ is determined by a class of finitely
presented modules. More precisely, $\mathcal B = \mathcal S^\perp$
where $\mathcal S= {}^\perp \mathcal B\cap \m$. Combining this
with \cite[Theorem 2.2]{aht} and \cite[Theorem 4.14]{T} we obtain

\Bth{oneone} (1) The tilting classes in $\M$ correspond
bijectively to the resolving subcategories of $\m$. The
correspondence is given by the mutually inverse assignments
$$\alpha: \calb\mapsto
{}^\perp\!\calb\cap\m\quad\mbox{and}\quad \beta: \cals\mapsto
\cals^\perp$$ (2) The cotilting classes   in $\La\,\mbox{Mod}$
correspond bijectively to the resolving subcategories of $\m$. The
correspondence is given by the mutually inverse assignments
$$\gamma: \calc\mapsto
{}^\intercal\!\calc\cap\m\quad\mbox{and}\quad \delta: \cals\mapsto
\cals^\intercal={} ^\perp\!(\cals^\ast)$$ (3) The above
correspondences yield a one-to-one-correspondence between tilting
classes in $\M$ and cotilting classes   in $\La\,\mbox{Mod}$. \Eth

\Brem{cals} (1) $\alpha, \beta, \gamma, \delta$ are
order-reversing: If $\calb_1,\calb_2$ are two tilting classes
with $\calb_1\subset\calb_2$, then
$\alpha(\calb_2)\subset\alpha(\calb_1)$, and the analogous
property holds for the remaining assignments.

(2) Any resolving subcategory of $\m$ is closed under submodules,
since it occurs as $^\perp\!\calb\cap\m$ for some class
$\calb\subset\M$ and all modules in $\M$ have injective dimension at
most one.

(3) Let $\cals$ be a subcategory of $\m$ containing $\La$, and
assume that $\cals$  is closed under predecessors, that is, if
$X\in\m$ is an indecomposable module with a nonzero map $X\ra S$
to a module $S\in\cals$, then $X\in\cals$. Then it is easy to see
that $\cals$ is resolving. \Erem

In particular we have the following examples:

\Bex{p} The category $\add{\p}$ is a resolving subcategory of \m\
with $\beta (\add{\p})=\p^\perp=\Gen{\mathbf L},$ and $\delta
(\add{\p})={}^\perp\!_{\La}\q = \Cogen {}_\La\mathbf W$ where
$\mathbf L$ and $\mathbf W$ are the Lukas and the Reiten-Ringel
tilting modules respectively. \Eex

\Bex{put} Let $\tube'$ be a nonempty union of  tubes, and let
$\mathcal U$ be the set of quasi-simple  modules in $\tube'$. Then
the category $\add{(\p\cup \tube')}$ is a resolving subcategory of
\m\ with $\beta (\add{(\p\cup \tube')})={\tube'}^\perp=\Gen
T_{\mathcal U}$ where $T_\calu=R_\calu\oplus R_\calu/R$.

In fact, if $Z\in \p$ and $S$ is quasi-simple, then there is a
nonzero map from $Z$ to the ray $\{S[n]\mid n\in\N\}$ defined by
$S$, cf. \cite[XII, 3.6]{SS}. So, by the Auslander-Reiten formula we
deduce that the modules in $\tube'^\perp$ cannot have direct
summands in $\p$, and therefore $\tube'^\perp\subset \p^\perp$ and $
(\add{(\p\cup \tube')})^\perp=\tube'^\perp$. This implies $\beta
(\add{(\p\cup \tube')})=\mathcal U^\perp=\Gen T_{\mathcal U}$ by
Example \ref{Var}.

In particular,  $\beta (\add{(\p\cup
\tube)})=\tube^\perp=\Gen {\mathbf W}$. Moreover, with dual
arguments one proves that
 $\delta (\add{(\p\cup \tube)})={}^\perp\,_{\La}\tube= _{\La}\tube\,^o$ is the category of all torsion-free left ${\La}$-modules.
\Eex

The examples above give a complete
list of large tilting modules over the Kronecker algebra, as we are  going to  see  in
Corollary~\ref{resume}, as a consequence of the general
Theorem~\ref{class}.

\Blem{firstprop} Let $T$ be a tilting $\La$-module with tilting
class $\calb=T^\perp$, and
 $\cals={}^\perp\calb\cap \m$. Then
$D(T)$ is a cotilting module with cotilting class $ ^\perp
(\cals^\ast)=\{_{\La}X \mid  D(X)\in \calb\}$. \Elem \BB
The well-known Ext-Tor relations yield $^\perp
(\cals^\ast)=\cals^\intercal$ and $^\perp D(T)= T^\intercal
= \{_{\La}X \mid  D(X)\in \calb\}$. Now,  if
$D(X)\in\calb$, and $\cala={}^\perp\calb$, then
$\Ext_R^1(\cala, D(X))=0$, hence $\Tor_1^R(\cala, X)=0$,
and in particular $X\in\cals^\intercal.$
Conversely, since $T$ is a direct limit of modules from $\cals$ by \cite[4.4]{T}, we have  $\cals^\intercal\subset T^\intercal$. So, we have shown $^\perp D(T)={}^\perp (\cals^\ast)$. \\
We now deduce that $D(T)$ is a cotilting module. In fact, the
conditions (T1) and (T3) for $T$ yield the dual conditions (C1)
and (C3) for $D(T)$. Moreover, applying the Ext-Tor relations we
obtain that $D(T)\in\cals^\intercal$ since $T\in\cals^\perp$. So,
$D(T)\in{}^\perp D(T)$, and since $^\perp D(T)=\cals^\intercal$ is
closed under products, we infer
  $\Ext_R^1(D(T)^{\kappa}, D(T)) = 0$ for any cardinal $\kappa$, that is, the dual condition (C2)
is also satisfied. \EB

\Blem{fg} The \feq\ for a tilting $\La$-module $T$.
\begin{enumerate}[(1)]
\item $T$ is equivalent to a finitely generated tilting module.
\item $D(T)$ is equivalent to a finitely generated cotilting module.
\item All indecomposable direct  summands of $D(T)$ are
finitely generated.
\end{enumerate}\Elem \BB We will freely
use the results on {endofinite} modules collected in Lemma
\ref{inf}.

(1) $\Rightarrow$ (3): Let $T'$  be a
finite-dimensional tilting module equivalent to $T$. Clearly,
$\Add{T} =\Add{T'}$ implies $\Prod{D(T)} =\Prod{D(T')}$.  Then the
indecomposable direct summands of $D(T)$ belong to $\Prod
D(T')=\Add{D(T')}$, and are therefore isomorphic to
indecomposable direct summands of $D(T')$.\\
(3) $\Rightarrow$ (2): By a well-known result of Bongartz
\cite{B},
the number of isoclasses of indecomposable direct summands of $D(T)$ is bounded by the number of isoclasses of simple $\La$-modules, and $D(T)$ is equivalent to a finitely generated cotilting module.\\
(2) $\Rightarrow$ (1): Let $_{\La}C$  be a finite-dimensional
cotilting module equivalent to $D(T)$. Then $D(T)$ belongs to
$\Prod C=\Add{C}$, and is thus isomorphic to a direct sum of copies
of  a finite number of indecomposable finitely generated modules.
In particular, this implies that $D(T)$ is endofinite. But then
$T$ is a pure submodule of the endofinite module $D^2(T)$ and is
therefore a direct summand of $D^2(T)$ by \cite[4.3]{CB2}. In particular, also $T$ is
isomorphic to a direct sum of copies of  a finite number of
indecomposable finitely generated modules, which proves (1).
\EB

\Bth{class} Let $T$ be a tilting $\La$-module with tilting
class $\calb=T^\perp$, and $\cals={}^\perp\calb\cap \m$.
Assume that $T$ is not  equivalent to a finitely generated
tilting module. Then the following hold true.
\begin{enumerate}[(1)]
\item $T$ is a regular module
and $\Gen{\mathbf W}\subset \Gen T\subset \Gen{\mathbf L}$.
\item There is a subset $\tube'\subset \tube$ such that  $\cals=\add{(\p\cup \tube')}$.
\item If $\tube'=\emptyset$, then $T$ is equivalent to the Lukas tilting module $\mathbf L$.
\item If $\tube'$ is a non-empty union of tubes, and
$\mathcal U$ is the set of quasi-simple  modules in
$\tube'$, then  $T$ is equivalent to $T_{\mathcal U}$.
\end{enumerate}\Eth \BB By assumption and Lemma \ref{fg}, the module
$_{\La}D(T)$ has an indecomposable direct summand $M$ which
is infinite dimensional. Observe that $M$ is pure-injective
as it is a summand of a dual module. From Lemma
\ref{inf}(6) and Lemma \ref{firstprop}, we infer that
$\cals^\ast$ cannot contain modules from $_{\La}\p$, hence
$\cals$ cannot contain modules from $\q$. Similarly,
${}^\perp D(T)$ cannot contain modules from $_{\La}\q$,
hence $\calb=\cals^\perp$ cannot contain modules from $\p$.
But then $\calb\subset \p^\perp= \Gen {\mathbf L}$, thus
$\p\subset {}^\perp\calb\cap\m=\cals$. So, $T$ is a regular
module, and we have verified (1) and (2). Now (3) and (4)
follow immediately from Examples \ref{p} and \ref{put}. \EB

\Bcor{resume} Over the Kronecker-algebra, every tilting
module is either equivalent to a finitely generated tilting module, or to precisely one of the modules in the following list:
\begin{enumerate}[\ -]
\item the Lukas tilting module $\mathbf L$,
\item the tilting modules of the form $T_{\mathcal U}$ for a non-empty
set of quasi-simples $\mathcal U$. \end{enumerate}
In other words, there is a one-one-correspondence between the subsets of $\mathfrak{T}$ and the equivalence classes of large tilting modules.\Ecor \BB
Assume that $T$ is not equivalent to a finitely generated
tilting module. With the notation of Theorem~\ref{class},
we note that $\tube'$ can only contain modules from
homogeneous tubes. Then, with any regular module $M\in
\tube'$, the resolving subcategory $\cals$ contains also
its regular socle $S$ by Remark \ref{cals}(2), and  so it
contains the whole (homogeneous) tube  $S$ belongs to. This
shows that $\tube'$ is a union of tubes, so the claim
follows from Theorem~\ref{class}.
In particular, the large tilting modules are  parametrized, up to equivalence, by the subsets of $\mathfrak{T}$; hereby, the empty set corresponds to the equivalence class of $\mathbf L$. \EB


\section{Finite dimensional direct summands}\label{sec:finitedimensional}

In this section we describe the finite dimensional direct summands
of a large tilting module $T$.  They are regular modules
whose indecomposable summands belong to non-homogeneous tubes. We show that these summands are arranged in disjoint wings, and that
the number of summands from each wing equals the number of
quasi-simple modules in that wing. Moreover, the summands  contributed by each tube $\tube_\lambda$ are determined by the  intersection $\tube_\lambda\cap \cals$ of the tube with the resolving subcategory $\cals$ corresponding to $T$. Special attention will  be devoted to the case when $\cals$ contains a complete ray from $\tube_\lambda$.

\Blem{fginAddT}
If $T$ is a large tilting $R$-module, then
every finitely generated indecomposable module
$X\in \Add T$ is a regular module  from  a non-homogeneous tube, and its
regular length  $m<r$ is bounded by
the rank $r$ of the tube. Thus there are at most finitely many
non-isomorphic finitely generated indecomposable modules that can occur as direct summands of large tilting modules.
\Elem

\BB Suppose that $T$ has tilting class $\calb=T^\perp$ and set
$\cals={}^\perp\calb\cap\m$.

Notice that $X$ is   isomorphic to a direct summand of $T$
(cf.~Lemma~\ref{inf}), so it follows from Theorem~\ref{class}(1)
that  $X$ is a regular module, and there exist a tube
$\tube_\lambda$ and  a quasi-simple module $S\in\tube_\lambda$ such
that $X=S[m]$. Now  $0=\Ext_R^1(X,X)\cong D\Hom(S[m],\tau S[m])$
implies that  the tube $\tube_\lambda$ has rank $r>1$.

Choose the numbering $S=U_1,\, U_2=\tau^-U_1, \ldots,
U_r=\tau^-U_{r-1}$ for the quasi-simples in
$\tube_\lambda$. Recall that  $\Hom_R(S[m],
U_{i}[m-i+1])\neq 0$ for all $1\leq i\leq m$, where we
suppose that $U_i=U_j$ whenever $i\equiv j \mod r$. Now, if
$m\ge r$, we consider the module $S[m-r+1]$. Since it is a
submodule of $X\in\cals$ and $\cals$ is closed under
submodules, we have $S[m-r+1]\in\cals$. On the other hand,
$\Ext_R^1(S[m-r+1], S[m])=D\Hom_R(S[m],U_r[m-r+1])\neq 0$,
contradicting the fact that $S[m]\in \Add T\subseteq
\cals^\perp$. So, we conclude that $m<r$.

Since there are at most finitely many  (at most three)
non-homogeneous tubes, the foregoing shows that there are at most
finitely many non-isomorphic finitely generated indecomposable
modules that can occur as direct summands of large tilting modules.
\EB

{\it From now on in this section, we fix  a tilting $R$-module $T$
with tilting class $\mathcal{B}=T^\perp$}. We  work in a more
general setting which is needed for the proof of our main result
Theorem~\ref{main}: {\it we assume that
$\mathcal{S}={}^\perp\mathcal{B}\cap\m$ does not contain any
non-zero preinjective module, thus $\mathcal{S}=\add(\p'\cup\tube')$
where $\p'\subseteq \p$ and $\tube'\subseteq\tube$}. Of course,
every large tilting module satisfies this assumption by
Theorem~\ref{class}(2).

\Brem{trivial}
If $X\in \Add T\cap\m$, then $X\in\cals$. Indeed, $\Add T\cap\m= \calb\cap{}^\perp \calb\cap\m=\calb\cap\cals$.
\Erem

\Blem{old}
 Let $\tube_\lambda$ be a tube of rank $r>1$, and let $S$ be a
quasi-simple module in $\tube_\lambda$. Choose the numbering $S=U_1,\,
U_2=\tau^-U_1, \ldots, U_r=\tau^-U_{r-1}$ for the quasi-simples in
$\tube_\lambda$.
\begin{enumerate}[(1)]
\item If $\cals$ contains some, but not all modules
from the ray $\{S[n]\mid n\in\N\}$, then there is $m<r$ such that
$S[m]\in\Add T$. More precisely, if $S[m]$ is the module of
maximal regular length  in
$\cals\cap\{S[n]\mid n\in\N\}$, then $S[m]\in\Add T$.

\item If $S[m]\in
\tube_\lambda\cap\Add T$, then  the rays starting at $U_2, \ldots, U_{m+1}$
are not completely contained in $\cals$. In fact, $U_2[m],\
U_3[m-1],\dotsc,\ U_{m+1}\notin\mathcal{S}$.

\item  If $S[m]\in
\tube_\lambda\cap\Add T$, then
$\mathcal{W}_{S[m]}\cap \Add T$ contains
precisely $m$ modules which are uniquely determined by $\cals\cap \mathcal{W}_{S[m]}$.
\end{enumerate}
\Elem
 \BB
 (1)  Assume that $S[r]\in\cals$. We claim that $S[n]\in\cals$ for
all $n\geq 1$. For $\cals$ is closed under submodules, thus
$S[l]\in\cals$ for all $1\leq l\leq r$. If $n>r$, write $n=kr+l$
with $r\leq kr<n$ and $1\leq l\leq r$ and consider the
exact sequence $\exs {S[kr]}{}{S[n]}{}{S[l]}$. Now the claim follows by induction on $n$ since $\cals$ is closed
under extensions.

Thus there exists $m<r$ such that $S[m]\in\cals$ and
$S[m+1]\notin\cals$. We prove that $S[m]\in\Add T$. By Remark \ref{trivial}, we have
 to show that $S[m]\in\cals^\perp$. Take a module
$Z\in\cals=\add(\p'\cup\mathbf{t}')$, w.l.o.g.~assume that $Z$
is indecomposable. If $Z\in\p'$, then
$\Ext_R^1(Z,S[m])=D\Hom_R(\tau^- S[m],Z)=0$. If
$Z\in\mathbf{t}'$, we can assume w.l.o.g.~that $Z$
belongs to $\mathbf{t}_\lambda$. If
$\Ext_R^1(Z,S[m])=D\Hom_R(S[m],\tau Z)\neq 0$, we  would have
$Z=U_{i+1}[m-i+1+l]$ for $1\leq i\leq m$ and $0\leq l$. But then the
exact sequence $\exs{U_1[i]}{}{U_1[m+1+l]}{}{U_{i+1}[m-i+1+l]}$
together with the fact that $\cals$ is closed under extensions would
imply that $S[m+1+l]\in\cals$,  contradicting the choice of $m$. We
conclude that $S[m]$ belongs to $\cals^\perp$ and thus to $\Add T$.

(2) All modules of regular length at most $m$ on the coray
ending at $U_m$ are quotients of $S[m]$ and therefore
belong to the tilting class $\calb$. Hence the modules of
regular length at most $m$ on the coray ending at $U_{m+1}$
cannot be in $\cals$ by the AR-formula. This yields the
claim, because these modules lie on the rays starting at
$U_2, \ldots, U_{m+1}$.

(3) We show by induction on $m$ that $\mathcal{W}_{S[m]}\cap \Add T$ contains $m$ modules. Our proof will show how the $m$ modules are determined by $\cals\cap\mathcal{W}_{S[m]}$. The result
clearly holds for $m=1.$
Let $m>1$. First of all, note that the modules
$U_2[m-1],U_3[m-2],\dotsc,U_m$ on the coray ending at $U_m$ are in $\cals^\perp$ because they are
epimorphic images of $S[m]$.

Suppose that none of the modules
$U_2[m-1],U_3[m-2],\dotsc,U_m$ belongs to $ \cals$. Then no regular module
containing any of these modules can belong to $\cals$. On the other hand,  for
 $X\in \p\cup\mathbf{t}$ we have
$\Ext_R^1(X,S[m-1])=D\Hom_R(U_2[m-1],X)\neq 0$ if and only if $X$
is a regular module that contains one of the modules
$U_2[m-1],U_3[m-2],\dotsc,U_m$. Hence $S[m-1]$ belongs to
$\cals^\perp$, and as a submodule of $S[m]\in\cals$ it also belongs to $\cals$, therefore $S[m-1]\in\Add T$. So $\Add T$
contains precisely $m$ modules in $\mathcal{W}_{S[m]}$: these are
$S[m]$ and the $m-1$ modules in $\mathcal{W}_{S[m-1]}$ given by
the induction hypothesis.

Suppose now that one of  $U_2[m-1],U_3[m-2],\dotsc,U_m$ belongs to
$\cals$. Choose
$U_{i+1}[m-i]\in\cals$ of maximal regular length. Then $U_{i+1}[m-i]\in \Add T$, and the
induction hypothesis implies that $\Add T$ contains precisely
$m-i$ modules in $\mathcal{W}_{U_{i+1}[m-i]}$.

Since $U_{i+1}[m-i]$ (and
its submodules on the ray starting at $U_{i+1}$) are in $\cals$, no module of regular length at most $m-i$ on the ray starting at $U_i$ can belong to $\cals^\perp$. This shows that $S[i],S[i+1],\dotsc,S[m-1]\notin\cals^\perp$. We claim that $S[i-1]\in\Add T$.
To this end, we note  that  for
$X\in\p\cup\mathbf{t}$ we have
$\Ext_R^1(X,S[i-1])=D\Hom_R(U_2[i-1],X)\neq 0$ if and only if $X$ is regular and contains
one of the modules $U_2[i-1], U_3[i-2],\dots,
U_i$ as a submodule. But none of $U_2[i-1], U_3[i-2],\dotsc,
U_i$ can belong to $\cals$. Indeed, this follows from  the choice of $U_{i+1}[m-i]$, by using that
each of the modules $U_2[m-1],U_3[m-2],\ldots,U_i[m-i+1]\notin\cals$
can be written as an extension of one of the modules $U_2[i-1], U_3[i-2],\dots,
U_i$
by the module $U_{i+1}[m-i]\in \cals$.

So, we infer that $S[i-1]\in\Add T$, and the induction hypothesis
implies that $\Add T$ contains precisely $i-1$ modules in
$\mathcal{W}_{S[i-1]}$.
We conclude that $\Add T$ contains precisely $m$ modules in
$\mathcal{W}_{S[m]}$: these are the $m-i$ modules in
$\mathcal{W}_{U_{i+1}[m-i]}$, the $i-1$ modules in
$\mathcal{W}_{S[i-1]}$, and $S[m ]$.\EB

\medskip

The following result shows that the indecomposable summands of $T$
from a tube $\tube_\lambda$ are arranged in disjoint wings, and that
the union of such wings does not contain all quasi-simples from
$\tube_\lambda$.

\Bcor{disjointwings}
Let $X,X'$ be two finitely generated indecomposable modules in $\Add T$, and let $\mathcal W_X,  \mathcal W_{X'}$ be the corresponding wings. Then either $\mathcal{W}_X\subset \mathcal{W}_{X'}$ or $\mathcal{W}_{X'}\subset \mathcal{W}_X$ or $\mathcal{W}_X\cap\mathcal{W}_{X'}=\emptyset$. Moreover, given a tube $\tube_\lambda$ of rank $r>1$,
 the  quasi-simple modules
in the union of all wings $\mathcal W_X$ with $X\in\tube_\lambda\cap\Add T$
do not form a complete clique, and there are at most $r-1$ isomorphism classes of modules in $\tube_\lambda\cap\Add T$.
\Ecor
\BB
We can assume w.l.o.g.~that $X,X'$
belong to the same tube $\mathbf{t}_\lambda$. Let $S,S'$ be
quasi-simples in $\mathbf{t}_\lambda$ such that $X=S[m]$ and
$X'=S'[m']$. Assume that $m\leq m'$, and suppose that
$\mathcal{W}_X\nsubseteq\mathcal{W}_{X'}$ and
$\mathcal{W}_X\cap\mathcal{W}_{X'}\neq\emptyset$.

We have to consider two cases. In the first case,  the coray
$\mathfrak{c}'$ that contains $S'[m']$ meets the ray
$\mathfrak{r}$ determined by $S$ in a module $S[l]\in\mathfrak{c}'\cap\mathfrak{r}$ with $1\leq
l \le m$. We even have $l<m$ since otherwise $\mathcal{W}_X\subseteq\mathcal{W}_{X'}$. Then
$S[l+1]\in\mathcal{S}$ and $\tau S[l+1]\in\mathcal{S}^\perp$
because $\mathcal{S}$ and $\mathcal{S}^\perp$ are closed under
submodules and images respectively. But $\Ext_R^1(S[l+1],\tau
S[l+1])=\Hom_R(\tau S[l+1],\tau S[l+1])\neq 0$, a
contradiction.

In the second case, the coray $\mathfrak{c}$ that
contains $S[m]$ meets the ray $\mathfrak{r}'$ determined by
$S'$ in a module $S'[l]\in
\mathfrak{c}\cap\mathfrak{r}'$, where again  $l<m$ (otherwise
$\mathcal{W}_X\subseteq\mathcal{W}_{X'}$). Then $S'[l+1]\in\mathcal{S}$
and $\tau S'[l+1]\in\mathcal{S}^\perp$. But
$\Ext_R^1(S'[l+1],\tau S'[l+1])=\Hom_R(S'[l+1],S'[l+1])\neq
0$, again a contradiction.

For the proof of the second statement, let $U_1=S,
U_{2}=\tau^{-}U_1,\ldots, U_{m}$ be the quasi-simple modules in
$\mathcal W_X$. Then it follows from Lemma~\ref{fginAddT} that $m<r$
and from Lemma~\ref{old}(2) that $\tau^-U_{m}\notin \cals$ cannot be
a submodule of a module $X'\in\Add T$. Thus it cannot belong to any
wing $\mathcal W_{X'}$ with $X'\in\tube_\lambda\cap\Add T$.

Finally,  by Lemma \ref{old}(3) the number of isomorphism classes of modules in $\tube_\lambda\cap\Add T$ equals the number of quasi-simple modules in the union of all wings involved, hence it is at most $r-1$.
\EB

\medskip

Let us  now deal with the case that $\cals$ contains a complete ray from $\tube_\lambda$.

\Blem{lem:completeray} Let $\tube_\lambda$ be a tube of rank $r>1$, and let $S$ be a
quasi-simple module in $\tube_\lambda$. Suppose that the ray $\{S[n]\mid n\geq 1\}$
starting at  $S$ is completely contained in
$\cals$. Choose the numbering $S=U_1,
U_2=\tau^-U_1,\dotsc,U_r=\tau^-U_{r-1}$ for the quasi-simples in
$\mathbf{t}_\lambda$. The following assertions hold true.
\begin{enumerate}[(1)]
\item If the ray $\{U_2[n]\mid n\geq 1\}$
starting at $U_2$ is completely contained in $\cals$, then
$S[n]\notin\Add T$ for all $n\geq 1$.
 \item If $2<i\le r$ is the least number such
that the ray $\{U_i[n]\mid n\geq 1\}$ starting at $U_i$ is
completely contained in $\cals$, then $S[i-2]$ is the module
of maximal regular length in $\{S[n]\mid n\geq 1\}\cap\Add T$.
\item If
$\{S[n]\mid n\geq 1\}$ is the only ray of $\tube_\lambda$ which is completely
contained in $\cals$, then $S[r-1]\in\Add T$.
\end{enumerate}
\Elem

\BB (1) Clearly $\Ext_R^1(U_2[n],S[n])=D\Hom_R(S[n],S[n])\neq 0$.

(2) We have to verify $S[i-2]\in\cals^\perp$. Observe that, since $\cals$ is closed under submodules,
$\Ext_R^1(Z,S[i-2])=D\Hom_R(U_2[i-2],Z)=0$ for all $Z\in\cals$ if
and only if $U_2[i-2], U_3[i-3],\dotsc,U_{i-1}\notin\cals$. So, assume that one of
 the modules $U_2[i-2],
U_3[i-3],\dotsc,U_{i-1}$ belongs to $\cals$, say $U_j[i-j]\in\cals$ with $2\leq j\leq i-1$. Since  the rays
starting at $U_2,\dotsc,U_{i-1}$,  are not completely contained in
$\cals$, it follows from Lemma~\ref{old}(1) that $U_j[l]\in\Add T$ for some  $l\geq
i-j$. As $\cals^\perp$ is closed under epimorphic images,   there exists a module in $U_{i-1}[t]\in\cals^\perp$ on the ray starting at $U_{i-1}$. But this is a contradiction
because
$\Ext_R^1(U_i[t],U_{i-1}[t])=D\Hom_R(U_{i-1}[t],U_{i-1}[t])\neq
0$.

Hence $S[i-2]\in\Add T$. Moreover, no module of the form $S[i-2+t]$ with $t>0$ is in $\Add T$, because otherwise
its epimorphic image $U_{i-1}[t]$ would be a module in $\cals^\perp$ on the ray starting at $U_{i-1}$.

(3) Proceed as  in (2) and show that $U_2[r-2],\dotsc,
U_{r-1}\notin\cals$.\EB

If $\cals$ contains some, but not all rays from a tube $\tube_\lambda$, then it certainly contains the rays with modules of maximal
regular length in $\mathbf{t}_\lambda\cap\Add T$, as we are going to see next.

\Blem{completeraymaxlength} Let $\mathbf{t}_\lambda$ be a
nonhomogeneous tube. Suppose that $\cals$ contains a complete ray
from $\mathbf{t}_\lambda$.
For every module $X\in \mathbf{t}_\lambda\cap\Add T$ there is a module
$S[m]\in \tube_\lambda\cap\Add T$ lying on a ray $\{S[n]\mid n\geq 1\}$ which is
completely contained in $\cals$ such that $X$ belongs to the wing $\mathcal{W}_{S[m]}$.
More
precisely, $S[m]$ can be chosen to be either $S[i-2]$ as in
Lemma~\ref{lem:completeray}(2) or $S[r-1]$ as in
Lemma~\ref{lem:completeray}(3). \Elem

\BB Let $S'\in\mathbf{t}_\lambda$ be a
quasi-simple such that $S'[m']\in\mathbf{t}_\lambda\cap\Add T$ for
some $m'\geq 1$. Choose the numbering $U_1,
U_2=\tau^-U_1,\dotsc,U_r=\tau^-U_{r-1}$ for the quasi-simples in
$\tube_\lambda$ where the ray starting at $U_1$ is completely
contained in $\mathcal{S}$, $S'=U_j$ for some
$j\in\{1,\dotsc,r-1\}$, but no ray starting at $U_l$ is completely
contained in $\mathcal{S}$ for $2\leq l\leq j$. Note that  also the ray starting at $\tau^-S'=U_{j+1}$ is not
completely contained in $\mathcal{S}$ by Lemma~\ref{old}(2).

Set $S=U_1$. If there is no other $i\in\{1,\dotsc,r\}$ such that the
ray starting at $U_i$ is completely contained in $\mathcal{S}$, then
$S[r-1]\in\Add T$ by Lemma~\ref{lem:completeray}(3). The result then
holds by Lemma~\ref{fginAddT} and Lemma~\ref{old}(2).

If $i\in\{j+2,\dotsc,r\}$ is the first number such that
the ray  $\{U_i[n]\mid n\geq1\}$ is completely
contained in $\mathcal{S}$, then $S[i-2]\in\Add T$ by
Lemma~\ref{lem:completeray}(2). Since $U_j[m']=S'[m']\in \Add T$, we
know that the rays starting at $U_{j+1},\dotsc,U_{j+m'}$ are not
completely contained in $\mathcal{S}$ by Lemma~\ref{old}(2). Hence
$i\geq j+m'+1$. Thus $i-2\geq j$ and $i-2\geq m'$. The first inequality implies that
$U_j=S'\in\mathcal{W}_{S'[m']}\cap \mathcal{W}_{S[i-2]}$. By Corollary~\ref{disjointwings}, the second inequality implies that
$\mathcal{W}_{S'[m']}\subseteq \mathcal{W}_{S[i-2]}$. Therefore
$S'[m']\in\mathcal{W}_{S[i-2]}$. \EB

Let us summarize our discussion on $\tube_\lambda\cap\Add T$.

\Bprop{branchmodules} Let $\tube_\lambda$ be a  tube of rank $r$.
Then $\tube_\lambda\cap \cals$ determines $\tube_\lambda\cap\Add T$.
More precisely:
\begin{enumerate}[(1)]
\item If $\tube_\lambda\cap\cals=\emptyset$, then $\tube_\lambda\cap\Add T=\emptyset$.
\item If $\tube_\lambda\subseteq\cals$, then $\tube_\lambda\cap\Add T=\emptyset$.
\item If $\emptyset\neq \tube_\lambda\cap\cals \subsetneq \tube_\lambda$, then $\tube_\lambda\cap\cals$ determines unique quasi-simples $S_1,\dotsc,S_l\in\tube_\lambda$ and unique $m_1,\dotsc,m_l\in\mathbb{N}$
 such that \begin{enumerate}
\item[(a)] $S_j[m_j]\in \tube_\lambda\cap\Add T$ for $j=1,\dotsc,l$.
\item[(b)] $\mathcal{W}_{S_{j_1}[m_{j_1}]}\cap
\mathcal{W}_{S_{j_2}[m_{j_2}]}=\emptyset$ if $j_1\neq j_2$.
\item[(c)] $\tube_\lambda \cap\Add T \subseteq \bigcup_{j=1}^l\mathcal{W}_{S_j[m_j]}.$
\end{enumerate}
For each $j\in\{1,\dotsc,l\}$, there are exactly $m_j$ modules from $\mathcal{W}_{S_j[m_j]}$ in $\tube_\lambda\cap\Add T$  and they are uniquely determined by $\cals\cap\mathcal{W}_{S_j[m_j]}$. Therefore there are exactly
$m_1+\dotsb+m_l<r$ modules in  $\tube_\lambda\cap\Add T$.
\end{enumerate}
\Eprop

\BB By Remark~\ref{trivial}, every finite dimensional indecomposable module in $\Add T$ belongs to $\cals$. Thus (1) follows.

(2) holds by Lemma~\ref{lem:completeray}(1).

(3)
If $\tube_\lambda\cap\cals\not=\emptyset$ contains no complete ray, then there exist unique
quasi-simples $S_1,\dotsc,S_l$ and $m_1,\dotsc,m_l\in\mathbb{N}$
verifying (a), (b) and (c)  by
Lemma~\ref{old}(1) and Corollary~\ref{disjointwings}.

If $\tube_\lambda\cap \cals$ contains a complete ray, then there exist unique quasi-simples
$S_1,\dotsc,S_l$ and $m_1,\dotsc,m_l\in\mathbb{N}$ verifying (a),
(b) and (c) by
Lemma~\ref{completeraymaxlength} and Corollary~\ref{disjointwings}.

In both cases 
Lemma~\ref{old}(3) implies that $\tube_\lambda\cap\Add T$ contains exactly
$m_j$ modules from each $\mathcal{W}_{S_j[m_j]}$ and that these $m_j$ modules are uniquely determined by $\cals\cap\mathcal{W}_{S_j[m_j]}$.

Altogether, $\tube_\lambda\cap\Add T$ consists of  $m_1+\dotsb+m_l$ modules, and  $m_1+\dotsb+m_l<r$ by Corollary~\ref{disjointwings}.
\EB

\Bdef{def:vertices} Let $\tube_\lambda$ be a tube. The modules
$S_1[m_1],\dotsc,S_l[m_l]$ satisfying (a), (b) and (c) in
Proposition~\ref{branchmodules} will be called the \emph{vertices}
of $T$ in $\tube_\lambda$. \Edef

We now want to describe the  regular modules that can occur as the finite dimensional part of $T$.

\Bdef{def:Y} Recall that a module $Y$ is said to be {\it exceptional} if $\Ext_R^1(Y,Y)=0$. Inspired by \cite[4.4]{Rbook}, we will say that a finite dimensional regular multiplicity free exceptional
$R$-module  $Y$ is  a \emph{branch module} if it satisfies the
following condition: \begin{enumerate}[(B)]
\item For each quasi-simple module $S$ and $m\in\N$ such that $S[m]$
is a direct summand of $Y$, there exist precisely $m$ direct summands of $Y$
that belong to $\mathcal{W}_{S[m]}$.
\end{enumerate} \Edef

Let $T$ be a tilting module
with tilting class $\calb=T^\perp$ such that $\cals={}^\perp\calb\cap\m$ 
 does not contain any non-zero preinjective module. By Lemma~\ref{old}(3),  the
direct sum  $Y$ of a complete irredundant set of  finitely generated
indecomposable direct summands of $T$ is a branch module.
 The following result shows
that there do not exist any other branch modules. We will even see
in Theorem~\ref{main} that every branch module does occur as a direct
summmand of a large tilting module.

\Blem{finitedimensionaltilting} Every finite dimensional regular
multiplicity free exceptional  module $Z$  is a direct summand of
 a finite dimensional tilting $R$-module $H=H_0\oplus Y$
satisfying the following properties:
\begin{enumerate}[(a)]
\item $H_0\not=0$ is a preprojective module.
\item $Y$ is a branch module  with  the same
quasi-simple composition factors as $Z$.
\item $H^\perp=Z^\perp=Y^\perp$.
\item $\cals_H={}^\perp(H^\perp)\cap\m$  does not contain any non-zero preinjective module.
\end{enumerate}
In particular, if $Z$ is a branch module, then $H=H_0\oplus Z$.
\Elem

\BB The module $Z$ is a partial tilting module, so by a well known construction due to Bongartz,
 taking a universal extension
 $\exs R \iota {R_0} \pi {Z^{(c)}}$
where $c={\rm dim}_k\Ext_R^1(Z,R)$, we obtain a finitely generated
tilting $R$-module $H=R_0\oplus Z$ with $H^\perp=Z^\perp$. Hence
$\q\subseteq H^\perp$ and therefore
$\cals_H={}^\perp(H^\perp)\cap\m$  does not contain any non-zero
preinjective module. So $R_0=H_0\oplus Y_0$ with $H_0$ preprojective
and $Y_0$ regular, and $H_0\not=0$ since there are no finite
dimensional regular tilting modules over $R$ (indeed, if $H_0=0$,
then $H$ is a direct sum of modules from non-homogeneous tubes, and
the number of isomorphism classes from each such tube is smaller
than the rank of the tube by Proposition~\ref{branchmodules}, so the
number of pairwise non-isomorphic indecomposable summands of $H$ is
strictly smaller than the  number of isomorphism classes of simples,
see the table in \cite[p.335]{R}). Observe that the regular module
$\ker\pi_{|{Y_0}}$ is contained in the preprojective module $\im
\iota$. Thus $\ker\pi_{|{Y_0}}=0$ and $Y_0\subseteq Z^{(c)}$.

We can suppose that $H=H_0\oplus Y$ where $Y=Y'\oplus Z$ is a direct
sum of a complete irredundant set of the indecomposable direct
summands of $Y_0\oplus Z$. Then $Y$ is a branch module by
Lemma~\ref{old}(3), and it has the same quasi-simple composition
factors as $Z$
 since $Y_0\subseteq
Z^{(c)}$.

Finally, note that any exceptional module which is a direct sum of
modules from a wing $\mathcal W_{S[m]}$  can have at most $m$
non-isomorphic indecomposable summands. So, if $Z$ is a branch
module, then   the fact that
$Y_0\subseteq Z^{(c)}$ implies $Y_0\in\add Z$ and therefore $Y=Z$.
\EB

\Brem{findimversion} Let $Z$ be a branch module. It can be proved
that the tilting module $H$ above is equivalent to $ R_\calu\oplus
Z$ where $R_\calu$ is the universal localization of $R$ at the set
$\calu$ of quasi-simple composition factors of $Z$. \Erem


\section{Decomposing tilting
modules}\label{section:decomposing}

{\it Throughout this section, we fix a tilting $\La$-module $T$ with
tilting class $\calb=T^\perp$ and $\cals={}^\perp\calb\cap \m$. We
assume that $T$ is not equivalent to a finite dimensional tilting
module}. We prove a structure result for the modules in $\calb$,
from which we derive
 a  canonical decomposition for $T$.

 \medskip

We are going to use two torsion pairs first studied by Ringel in
\cite{R}. The first is the split torsion pair $(\cald,\cald^\circ)$
whose torsion class is the class $\cald$ of the divisible modules.
We call a module \emph{reduced} if it belongs to the corresponding
torsion-free class $\cald^o$. The second  is the non-split torsion
pair $(\Gen\tube,\mathcal{F})$ with torsion class $\Gen\tube$. Here
$(\Gen\tube)\,^o=\calf$ is the class of  all torsion-free modules,
cf.~\ref{Ringel} and
 \ref{Lukas}.

We will further need the following canonical decomposition
of the regular modules in $\Gen \tube$. Write
$$\tube=\bigcup_{\lambda\in\mathfrak{T}}\tube_\lambda$$ where the
$\tube_\lambda$ are the tubes in the Auslander-Reiten
quiver of $\La$, and set $\calt_\lambda=\varinjlim \text{
add } \tube_\lambda$.
For $X\in\M$ denote by $\tube_\lambda(X)$  the maximal submodule of $X$ belonging to
$\calt_\lambda$.
As shown in \cite[4.5]{R}, every
regular module $X\in  \Gen\tube$ has a unique decomposition
$$X=\bigoplus_{\lambda\in\mathfrak{T}} \tube_\lambda(X).$$

We will say that a Pr\"ufer module $S[\infty]$
(or an adic module $S[-\infty]$) {\em belongs to a tube} $\tube_\lambda$ if $S$ is a quasi-simple module in (the mouth of) $\tube_\lambda$.

 \medskip

 Let us start by investigating the modules in the tilting class $\calb= \Gen
 T$. Since $\cals$ consists of finitely presented modules,
 the class $\mathcal{B}$ is \emph{definable}, i.~e., it is closed under direct limits, direct products,
and pure submodules.

\Bprop{Sperp} Let $X\in\calb=\cals^\perp$. Then \\
(1) $X=X_\cald\oplus X_{\rm red}$ where $X_\cald\in\cald$ is divisible, and $ X_{\rm red}$ is reduced.\\
(2) There is a pure-exact sequence $0\to X'\to X_{\rm red} \to \overline{X}\to 0$ where $\overline{X}\in\calb$ is torsion-free, and $X'\in \Gen \tube$.\\
(3) $X'=\bigoplus_{\lambda\in\mathfrak{T}} X_\lambda$, and for
each $\lambda$ there is a pure-exact sequence $0\to
A_\lambda\to X_\lambda\to Z_\lambda\to 0$ where $A_\lambda$
is a direct sum of modules in $\tube_\lambda\cap\calb$, and
$Z_\lambda\in\mathcal{B}$ is a direct sum of Pr\"ufer
modules belonging to the tube $\tube_\lambda$. \Eprop \BB
For (1) and (2), we refer to \cite[4.7 and 4.1]{R}.

(3) Note that the torsion-free class of reduced modules is closed
under submodules, and the tilting class $\calb=\cals^\perp$ is
definable, hence closed under pure-submodules. So, we infer from (1)
and (2) that $X'$ is a reduced module in $\calb$. Since preinjective
modules are divisible, it follows that $X'$ has no indecomposable
summands from $\q$. Moreover, $X'$ has no indecomposable summands
from $\p$ because $X'\in \Gen\tube$. Thus $X'$ is a regular module
in $\Gen \tube$ and has a decomposition
$X'=\bigoplus_{\lambda\in\mathfrak{T}} \tube_\lambda(X')$ as above
by \cite[4.5]{R}. We set $X_\lambda= \tube_\lambda(X')$. From
\cite[Theorem G and 4.8]{R} we know that there is a pure-exact
sequence $0\to A_\lambda\to X_\lambda\to Z_\lambda\to 0$ where
$A_\lambda$ is a direct sum of indecomposable modules of finite
length, and $Z_\lambda$ has no indecomposable direct summand of
finite length. Thus $Z_\lambda$ is regular,  and is therefore a
direct sum of Pr\"ufer modules. Again, we see that $A_\lambda$ is a
regular module in $\calb$, and since ${\rm
Hom}_\La(\tube_\nu,\calt_\lambda)=\varinjlim {\rm
Hom}_\La(\tube_\nu,\tube_\lambda)=0$ for $\nu\not=\lambda$, we infer
that $A_\lambda$ is a direct sum of modules in
$\tube_\lambda\cap\calb$. Similarly, we see that
$Z_\lambda\in\mathcal{B}$ and that  the Pr\"ufer modules occurring
as direct summands of $Z_\lambda$ admit non-zero maps from
$\tube_\lambda$ and therefore belong to the tube $\tube_\lambda$.
\EB

We can now refine the structure result of Proposition~\ref{Sperp} to
the modules in $\Add T$.
By Lemma~\ref{fginAddT}, there are at most finitely many
non-isomorphic finitely generated indecomposable modules in $\Add
T$. {\it We denote by $Y$ the direct sum of a complete irredundant
set of such modules}, which is a branch module by
Lemma~\ref{old}(3). Of course, $Y=\bigoplus_{\lambda\in\mathfrak{T}}
\tube_\lambda(Y)$ where $\tube_\lambda(Y)$ is the direct sum of  a
complete irredundant set of modules in $\tube_\lambda\cap\Add T$.

\Bprop{Add}  Every module  $X\in\Add T$  has a unique direct sum
decomposition $$ X=\bigoplus_{\lambda\in\mathfrak{T}}
\tube_\lambda(X)\oplus \overline{X}$$ where $\overline{X}$ is
torsion-free and each $\tube_\lambda(X)$ has a  decomposition in
torsion modules with local endomorphism ring. The indecomposable
summands of $\tube_\lambda(X)$ are isomorphic  to  direct summands
of $\tube_\lambda(T)$ and are either modules from $\tube_\lambda$ or
Pr\"ufer modules belonging to $\tube_\lambda$.

Moreover, every torsion (respectively, torsion-free) direct summand of $X$ is a direct summand of the torsion part $\bigoplus_{\lambda\in\mathfrak{T}}
\tube_\lambda(X)$ (respectively, of the torsion-free part $\overline{X}$).
\Eprop
\BB
 Let $X\in\Add T$. We know from \cite[4.1]{R} that there is a pure-exact
sequence $0\to X'\to X \to \overline{X}\to 0$ where
$\overline{X}\in\calb$ is torsion-free,  and
$X'\in\Gen\tube$. Note that $X\in\Add
T={}^\perp\calb\cap\calb$, and  ${}^\perp\calb$ is closed
under submodules, while $\calb$ is  closed under pure
submodules because it is a definable class. So, we infer that
$X'\in\Add T$ is a regular module in $\Gen \tube$, which by
\cite[4.5]{R} has a decomposition
$X'=\bigoplus_{\lambda\in\mathfrak{T}} \tube_\lambda(X)$. As in
the proof of Proposition \ref{Sperp}, we deduce from
\cite[Theorem G and 4.8]{R} that for each $\lambda$ there
is a pure-exact sequence $0\to A_\lambda\to
\tube_\lambda(X)\to Z_\lambda\to 0$ where $A_\lambda$ is a
direct sum of modules in $\tube_\lambda\cap\calb$, and
$Z_\lambda$ is a direct sum of Pr\"ufer modules belonging
to the tube $\tube_\lambda$. Again, we see that
$A_\lambda\in\Add T$, which implies by Lemma \ref{fginAddT}
that $A_\lambda$ has only finitely many non-isomorphic
indecomposable direct summands. In particular, this shows
that $A_\lambda$ is endofinite, thus pure-injective
(cf.~Lemma \ref{inf}), so the pure-exact sequence $0\to
A_\lambda\to \tube_\lambda(X)\to Z_\lambda\to 0$ splits,
and $\tube_\lambda(X)$ is a direct sum of modules in $\tube_\lambda\cap\Add T$ and Pr\"ufer modules
belonging to $\tube_\lambda$. In particular, $\tube_\lambda(X)$ has a decomposition in modules with local endomorphism ring.

 We infer that   $X'=\bigoplus_{\lambda\in\mathfrak{T}}
\tube_\lambda(X)$, being a direct sum of modules isomorphic to indecomposable direct summands of $Y$ or to
Pr\"ufer
modules, belongs to $\Add (Y\oplus\mathbf W)$. Now  $Y$ is
finite dimensional and therefore $\Sigma$-pure-injective (that is, every direct sum of copies of $Y$ is pure-injective),  and
$\mathbf W$ is $\Sigma$-pure-injective because $\Add {\mathbf
W}=\Prod {\mathbf W}$, see \cite[10.1]{RR}. Thus  $X'$ is
pure-injective, and  the pure-exact sequence $0\to X'\to X \to
\overline{X}\to 0$ splits, that is,
 $X=\bigoplus_{\lambda\in\mathfrak{T}}\tube_\lambda(X)\oplus\overline{X}$ has the stated decomposition. The uniqueness of $\overline{X}$ and the $\tube_\lambda(X)$ follows directly from torsion theory.

Let $A$ be a direct summand of  $X=X'\oplus \overline{X}$. Then
there are morphisms $\iota=(\iota',\overline{\iota}):A\to X$ and
$\pi=(\pi',\overline{\pi}):X\to A$ such that
$1_A=\pi\iota=\pi'\iota'+\overline{\pi}\,\overline{\iota}$. If  $A$
is torsion, then $\overline{\iota}=0$, so   $A$ is a direct summand
of $X'$. Similarly, if $A$ is torsion-free, then $\pi'=0$ and $A$ is
a direct summand of $\overline{X}$. In particular, each summand $A$
of $\tube_\lambda(X)$ belongs to $\Add \tube_\lambda(T)$. As
$\tube_\lambda(T)$ has a decomposition in modules with local
endomorphism ring, we deduce  from  the Theorem of
Krull-Remak-Schmidt-Azumaya that $A$ is isomorphic to an
indecomposable direct summand of $\tube_\lambda(T)$, see
e.g.~\cite[7.3.4]{Kasch}. \EB


\medskip

The following result  will be useful when dealing with the
torsion-free part $\overline{X}$ in the structure results from Propositions
\ref{Sperp} and \ref{Add}.

\Blem{pi} Let $\tube_\lambda$ be  a  tube.
\begin{enumerate}[(1)]
\item  $\cals$ contains a complete ray
$\{S[n]\,\mid\,n\geq1\}$ from $\tube_\lambda$ if and only
if $\calb$ does not contain any adic module belonging to
$\tube_\lambda$.

\item Suppose that $\tube_\lambda$ is a tube of
rank $r>1$ such that $\cals$  contains no complete ray from
$\tube_\lambda$. Let $\calu$ denote the set  of  quasi-simple modules in the union of all wings determined by the vertices of $T$ in $\tube_\lambda$.
Then for a quasi-simple
module $S\in\tube_\lambda$, the adic module $S[-\infty]$ belongs to $\calb$ if and
only if $S\notin\tau\calu=\{\tau U\mid U\in\calu\}$. Thus
$\calb$ contains precisely $r-|\calu|$ pairwise non-isomorphic adic modules belonging to
$\tube_\lambda$.

\item Let $\calu$ be a set of quasi-simple modules in
$\tube_\lambda$. Every torsion-free module in $\calb$ is
contained in $\calu^\perp$ if and only if  all adic modules
in $\calb$ belonging to $\tube_\lambda$  are contained in
$\calu^\perp$ (equivalently, every torsion-free module in
$\calb$ is an $R_\calu$-module if and only if all the adic
modules in $\calb$ belonging to $\tube_\lambda$ are
$R_\calu$-modules.)
\end{enumerate}
\Elem \BB We start by proving  the only-if part of (1). Suppose that
$\cals$ contains the complete ray $\{S[n] \mid n \geq   1 \}$.
Choose the numbering $S=U_1,\, U_2=\tau^-U_1, \ldots,
U_r=\tau^-U_{r-1}$, $r\geq 1$,  for the quasi-simples in
$\tube_\lambda$.  Consider $U_i[-\infty]$ for some $1\le i\le r$.
Then $\Ext^1_R(S[i+1], U_i[-\infty])\cong D\Hom_R(U_i[-\infty],\tau
S[i+1])= D\Hom_R(U_i[-\infty],U_r[i+1])\not=0$. Therefore
$U_i[-\infty]\notin\cals^\perp$ and hence
$U_i[-\infty]\notin\mathcal{B}$.

Next, we prove (2) and the  if-part of (1). First of all, observe
that  for any quasi-simple $S\in\tube_\lambda$ we have
$\Ext_R^1(\p,S[-\infty])=0$ because $S[-\infty]$ has no non-zero
preprojective summands. Also $\Ext_R^1(\tube_\mu,S[-\infty])=0$ for
all $\mu\neq\lambda$.

So, if $\tube_\lambda\cap\cals=\emptyset$,
then all adic modules belonging to $\tube_\lambda$ are  in $\calb$ (and indeed, this is the case $\mathcal U=\tau\mathcal U=\emptyset$ in (2)).

Assume now   $\tube_\lambda\cap\cals\not=\emptyset$, and suppose
that $\cals$ does not contain a complete ray from $\tube_\lambda$.
We know from Lemma \ref{old} that $\tube_\lambda\cap\cals$ is
contained in   the extension closure $\mathcal W$ of $\mathcal{U}$.
By Proposition~\ref{branchmodules}, $\calu$ does not contain a
complete clique. So, there are  quasi-simple modules
$S\in\tube_\lambda\setminus\tau\calu$. For such $S$ we have
$\Ext_R^1(U,S[-\infty])=D\Hom_R(S[-\infty],\tau U)=0$ for all
$U\in\calu$, so $S[-\infty]\in\mathcal U^\perp=   \mathcal W^\perp$
by Lemma \ref{wing}, and combining this with our first observation,
we conclude that $S[-\infty]\in\cals^\perp=\calb$. On the other
hand, if $S\in\tau\calu$, it is easy to see that
$\Ext_R^1(X_{j},S[-\infty])=D\Hom_R(S[-\infty],\tau X_{j})\neq0$ for
$X_j$ a vertex of $T$ in $\tube_\lambda$, which shows
$S[-\infty]\not\in\calb$ and completes the proof of (1) and (2).

(3) First of all, we note that the class of all torsion-free
modules $\mathcal F=\tube^o$, as well as the classes
$\calb$ and $\mathcal U^\perp$, are definable classes.
Indeed, $\mathcal F$ is clearly closed under direct
products and submodules, and it is closed under direct
limits since $\tube$ consists of finitely presented
modules. As for $\calb=\cals^\perp$ and $\mathcal U^\perp$,
closure under direct products is clear, and closure under
direct limits and pure submodules follows from the fact
that $\cals$ and $\mathcal U$ consist of finitely presented
modules.

We are now ready to consider a torsion-free module $X\in\calb$. Take
the pure-injective envelope $I$ of ${X}$, which is again  a
torsion-free module in $\calb$ as definable classes are closed under
pure-injective envelopes, see \cite[3.1.10]{GT}. Moreover, $I$ is
the pure-injective envelope of  $\bigoplus_{l\in L}I_l$, where
$\{I_l\,\mid\,{l\in L}\}$ is a complete irredundant set of
indecomposable summands of $I$, cf.~\cite[Chapter 8]{JL}. Now  the
$I_l$ are   indecomposable pure-injective torsion-free modules in
$\calb$, and  they are in $\Gen {\mathbf{L}}=\p^\perp$ by
Theorem~\ref{class}(1). We infer from the classification of the
indecomposable pure-injective modules reviewed in Lemma \ref{inf}
that $I_l$ is either  the generic module $G$ or an adic module.
Notice that $G$ is divisible and thus belongs to $\calu^\perp$.
Moreover, using the Auslander-Reiten formula, it is easy to see that
$I_l\in\calu^\perp$ if $I_l$ is an adic module belonging to a tube
$\tube_\mu$ with $\mu\neq\lambda$.  So,  if we assume that all adic
modules $I_l$ belonging to $\tube_\lambda$ are contained in
$\mathcal U^\perp$, then also $\bigoplus_{l\in L}I_l$ and its
pure-injective envelope $I$ are in  $\mathcal U^\perp$, and
therefore also the pure submodule ${X}$.

Conversely, recall that any adic module is torsion-free.\EB

\medskip

Let us determine the branch module $Y$ when   the tilting class
$\calb=\cals^\perp$ is the class of modules that are Ext-orthogonal
to a ray in a non-homogeneous tube, or in other words,  $\cals=\add
(\p\cup\tube')$ with $\tube'$ being a ray. This is a  special case
that will play an important role in the sequel.

\Bex{ray} Let  $S$ be a quasi-simple module, and assume
$\cals=\add(\p\cup\{\,S[n]\,\mid\,n\in\N\,\})$. Then
$$\calb=S[\infty]\,^\perp,\; \text{ and }\;S[\infty]\in\Add
T.$$ Moreover,  if $S$ belongs to a tube $\tube_\lambda$ of
rank $r>1$, then $$Y=S\oplus S[2]\oplus\dotsb\oplus
S[r-1]$$ (and $T\sim S\oplus S[2]\oplus\ldots\oplus
S[r-1]\oplus S[\infty]\oplus R_{\tube_\lambda}$,  as we
will see later in \ref{ray2}).

\smallskip

Indeed, we show as in Example \ref{put}  that
$\cals^\perp=\bigcap_{n\ge 1} S[n]^\perp$. Moreover, as $S[n]\subset
S[\infty]$, we have $S[\infty]^\perp\subseteq S[n]^\perp$ for all
$n\geq 1$, hence $S[\infty]^\perp\subseteq\calb.$ For the reverse
inclusion,
 note that $S[\infty]$ is filtered by $S[r]$,
thus every
 $X\in\calb\subset S[r]^\perp$ also belongs to $S[\infty]\,^\perp$ by \cite[3.1.2]{GT}.

Now we deduce $S[\infty]\in {}^\perp(S[\infty]^\perp)={}^\perp\calb$. Moreover,  $S[\infty]$ also belongs to $\calb$ as it is divisible and thus satisfies $\Ext_R^1(S[n],S[\infty])=0$ for all $n\in\N$. So, we conclude $S[\infty]\in\Add T.$

For the last claim, observe
that every finitely generated indecomposable module in $\Add T$ must belong to $\calb$, and also to ${}^\perp \calb\cap\m=\cals$ and thus to the ray  $\{\,S[n]\,\mid\,n\in\N\,\}$. So, we deduce from $\tube_\lambda\cap\calb=\mathcal W_{S[r-1]}$ that $S, S[2],\ldots, S[r-1]$ is a complete irredundant set  of the finitely generated indecomposable modules in $\Add T$.
\Eex


We are now ready for the main result of this section.

\Bth{structure}
There is a unique direct sum decomposition
$$T=\bigoplus_{\lambda\in\mathfrak{T}} \tube_\lambda(T)\oplus
\overline{T}$$ where $\overline{T}$ is torsion-free, and
$\tube_\lambda(T)$ is a direct sum of copies of the indecomposable direct summands of
$\tube_\lambda(Y)$ and of Pr\"ufer modules belonging to
$\tube_\lambda$.
More precisely,
for each tube $\tube_\lambda$ of rank $r$, the summand $\tube_\lambda(T)$  is given as follows:
\begin{enumerate}[(i)]
\item if
$\cals$  contains some modules from  $\tube_\lambda$, but no complete ray, then $\tube_\lambda(T)$ is a direct sum of at most $r-1$ pairwise non-isomorphic modules from $\tube_\lambda$ that are arranged in the disjoint wings determined by the vertices of $T$ in $\tube_\lambda$, and the number of non-isomorphic summands from each wing equals the number of quasi-simple modules in that wing;
\item if $\cals$ contains some   rays
from $\tube_\lambda$, then $\tube_\lambda(T)$ has precisely $r$ pairwise non-isomorphic indecomposable summands: these are  the $s$
Pr\"ufer modules corresponding to the $s\le r$ rays from $\tube_\lambda$ contained in $\cals$, and  $r-s$ modules from $\tube_\lambda$ which are arranged in the disjoint wings determined by the vertices of $T$ in $\tube_\lambda$;
\item $\tube_\lambda(T)=0$ whenever $\tube_\lambda\cap\cals=\emptyset$.
\end{enumerate}
\Eth
\BB
The existence of the decomposition follows from  Proposition \ref{Add}. Observe that
every  indecomposable direct summand of
$\tube_\lambda(Y)$ lies in $\tube_\lambda\cap\Add T$ and therefore
occurs as a direct summand in $\tube_\lambda(T)$ by Lemma~\ref{inf}(1).

We now turn to the additional statements. Note first that  the
finite dimensional direct summands  of $\tube_\lambda(T)$ are
contained in $\tube_\lambda\cap\cals$,  and further, recall that
${}^\perp\calb$ is closed under submodules, so  with every Pr\"ufer
module $S[\infty]$ it contains also the corresponding ray
$\{S[n]\mid n\in\N\}$. This proves (iii) and shows that
$\tube_\lambda(T)$ cannot contain infinite dimensional direct
summands if $\tube_\lambda\cap\cals$ does not contain complete rays.
Thus, in case (i),  the non-isomorphic indecomposable summands of
$\tube_\lambda(T)$ are  precisely the    indecomposable direct
summands of $\tube_\lambda(Y)$, and they have the stated properties
by Proposition~\ref{branchmodules}.

It remains to prove (ii). Assume that $\tube_\lambda\cap\cals$  contains
a complete ray $\{S[n]\mid n\in\N\}$. Then
  $\calb\subset\operatornamewithlimits{\cap}\limits_{n\geq
1}S[n]^\perp=S[\infty]^\perp$ by Example \ref{ray}, so $S[\infty]\in
{}^\perp\mathcal{B}$.
Further, using that
 $\mathcal{S}=\add(\p\cup\tube')$ for some $\emptyset\not=\tube'\subseteq\tube$,  we see
  that
$S[\infty]$ lies in $\calb$, hence in $\Add T$. We then infer from Proposition \ref{Add}  that $S[\infty]$ is a direct summand in
$\tube_\lambda(T)$.

Now we determine the remaining indecomposable summands of $\tube_\lambda(T)$.  If $\cals$ contains the whole tube $\tube_\lambda$, then it follows from Lemma~\ref{lem:completeray}(1) that $\tube_\lambda(T)$ has no finite dimensional summands, and therefore it is a direct sum of all Pr\"ufer modules belonging to $\tube_\lambda$.
  If $\tube_\lambda$ has rank $r>1$, and
$\cals$ contains $1\leq s<r$ complete rays form $\tube_\lambda$, we
get from Lemma~\ref{completeraymaxlength} and
Proposition~\ref{branchmodules} that there are exactly $r-s$ finite
dimensional indecomposable summands in $\tube_\lambda(T)$.
 \EB

\Brem{missing} There seems to be an asymmetry between case (i) and
(ii) in Theorem \ref{structure} above: the number $s$ of pairwise
non-isomorphic indecomposable summands of $\tube_\lambda(T)$ equals
the rank $r$ of $\tube_\lambda$ when $\cals$ contains some rays from
$\tube_\lambda$, but is smaller than $r$ otherwise. Note however
that in the latter case $s$ coincides with the number of
quasi-simple modules in the union of the wings determined by the
vertices of $T$ in $\tube_\lambda$.  So, the ``missing'' summands
are somehow ``replaced'' by the $r-s$ adic modules in $\calb$
established by Lemma \ref{pi}(2). This aspect will become more clear
in Remark~\ref{alternative} of the Appendix.
 \Erem

\section{Classifying tilting modules}

 Let again {\it $T$ be a tilting $\La$-module with tilting class
$\calb=T^\perp$, and $\cals={}^\perp\calb\cap \m$}. We assume that
$T$ is  {\it not equivalent to a finitely generated tilting module}. Then
we know from Theorem \ref{class} that {\it there is a subset
$\tube'\subset \tube$ such that  $\cals=\add{(\p\cup \tube')}$}. We
have seen in Theorem \ref{class} that $T$ is equivalent either to
the Lukas tilting modules if $\tube'$ is empty, or to a tilting
module arising from universal localization in case $\tube'$ is a
non-empty union of tubes. We now discuss the general case.

Recall that  we denote by {\it $Y$ the branch module defined as  the
direct sum of a complete irredundant set of the finitely generated
indecomposable modules in $\Add T$}. Thus $Y$ is a finite
dimensional direct summand of $T$ by Lemma~\ref{fginAddT} and
Lemma~\ref{inf}(1).

Our  aim is to reduce the classification problem to the
situation  considered in Theorem \ref{class}. To this end, we will
show that $T$ is equivalent to a tilting module of the form
$Y\oplus M$, where  $M$ has no finite dimensional indecomposable
direct summands and is a tilting module over a suitable universal
localization of $R$.  We will prove this step by step, by
considering the finitely many non-homogeneous tubes
$\tube_\lambda$ where $\tube_\lambda\cap\Add T\not=\emptyset$.

\smallskip

We first give a general criterion for constructing a tilting
module of the desired form.

\Blem{critilting} Let $Y'\in\Add T$, and let $M$ be a module
satisfying condition (T2), i.e. $\Ext_R^1(M,M^{(\kappa)})=0$ for any
cardinal $\kappa$. Then $Y'\oplus M$ is a tilting module equivalent
to $T$ if and only if the following hold true.
\begin{enumerate}[(a)] \item $\calb\subset M^\perp$, \item
$M\in\calb$, \item $T\in \Add (Y'\oplus M)$.\end{enumerate}
\Elem \BB For the only-if-part,  note that $\Add (Y'\oplus M)=\Add
T \subset \calb=(Y'\oplus M)^\perp$, which immediately yields (a),
(b), (c).

For the if-part, we show that $Y'\oplus M$ is tilting.
Condition (T1) is trivially verified. In order to check
(T2), let $\alpha$ be a cardinal. Then $\Ext^1_R((Y'\oplus
M), (Y'\oplus M)^{(\alpha)})\cong \Ext^1_R(Y',
Y'\,^{(\alpha)})\oplus \Ext^1_R(M, M^{(\alpha)})\oplus
\Ext^1_R( M, Y'\,^{(\alpha)})\oplus \Ext^1_R(Y',
M^{(\alpha)})$. Now the first term vanishes since $Y'\in
\Add T$, the second by assumption on $M$, the third term
vanishes by (a) since $\calb$ is closed under direct sums
and therefore $Y'\,^{(\alpha)}\in\calb$, and the last term
vanishes because $M^{(\alpha)}\in\calb=T^\perp\subset
Y'\,^\perp$ by property (b). Finally, condition (T3) is
satisfied by property (c).

So, $Y'\oplus M$ is a tilting module with $\Add T\subset \Add
(Y'\oplus M)$, thus $(Y'\oplus M)^\perp\subset T^\perp$.
Conversely, $T^\perp=\calb\subset Y'^\perp\cap M^\perp=(Y'\oplus
M)^\perp$, showing that $Y'\oplus M$ is equivalent to $T$. \EB

Now we proceed with our reduction. Given a  non-homogeneous tube
$\tube_\lambda$ where $\tube_\lambda\cap\Add T\not=\emptyset$, we
want to replace $T$ by an equivalent tilting module of the form
$\tube_\lambda(Y)\oplus M$ where $M$ is a  tilting module over a
suitable universal localization $R_\calu$ of $R$. To this end, we
replace the resolving subcategory $\cals$ by its localization
$$\cals_\calu=\{A\otimes_R R_\calu \mid A\in\cals\}$$ and choose $M$
to be a tilting $R_\calu$-module with tilting class
$\calb_\calu=\cals_\calu\,^\perp$. The existence of $M$ is guaranteed by
\cite[5.2.2]{GT}. We formulate criteria that will allow to perform
the replacement.

 \Bprop{critilting2}
Let $\tube_\lambda$ be a tube of rank $r>1$, and let
$\mathcal U$ be a set of $m<r$ quasi-simples in $\tube_\lambda$ with
extension closure $\mathcal W$. Assume that $M$ is an
$R_\calu$-tilting module with tilting class
$$\calb_\calu=\{X\in\M_\calu\,\mid\,\Ext^1_{R_\calu}(A\otimes_R R_\calu, X)=0 \;{\rm for\, all}\; A\in\cals\}$$ such that
 \begin{enumerate}[(i)]
 \item $\mathcal W\cup\mathcal W^o$ contains the subset $\tube'\subset \tube$ with $\cals=\add(\p\cup\tube')$,
 \item $\Add(\tube_\lambda\cap\calb)\subset M^\perp=\{X\in\M\mid \Ext_R^1(M,X)=0\}$,
 \item every adic module in $\calb$ belonging to $\tube_\lambda$ is contained in $\mathcal U^\perp$.
   \end{enumerate} Then $\tube_\lambda(Y)\oplus M$ is a tilting $R$-module equivalent to $T$.
 \Eprop
 \BB
 As $\mathbf{t}_\lambda(Y)\in \Add T$ and  $M$
 satisfies (T2), we only have to verify the conditions in Lemma~\ref{critilting}.

 (a) We prove $\calb\subset M^\perp$ in two steps.

\underline{Step 1:} We show $\calb\cap\mathcal U^\perp\subset M
^\perp$. Take  $X\in\calb\cap\mathcal U^\perp$. We claim
$M\in{}^\perp X$.  Since the $R_\calu$-tilting module $M$ is
filtered by the modules in $\cals_\calu=\{ A\otimes_R
R_\calu\,\mid\,A\in\cals\}$ by \cite[Lemma~4.5]{T}, it suffices to
show that $\cals_\calu\subset{}^\perp X$. So, let $A\in\cals$, w.l.o.g.
$A$ indecomposable. Then $A\in\p\subset\mathcal W^o$ or
$A\in\tube'\subset\mathcal W\cup\mathcal W^o$ by (i). If
$A\in\mathcal W$, then $A\otimes_R R_\calu=0$ by
Propostion~\ref{univloc}(5), so we can assume w.l.o.g. $A\in\mathcal
W^o$.  Then we know from Proposition~\ref{univloc}(7) that there is
a short exact sequence $0\to A\to A\otimes_R R_\calu\to A\otimes_R
R_\calu/R\to 0$ where  the two outer terms $A\in\cals$ and
$A\otimes_R R_\calu/R\in {}^\perp(\mathcal U^\perp)$ belong to
${}^\perp X$, so we infer that the middle-term $A\otimes_R R_\calu$
belongs to ${}^\perp X$ as well.

\underline{Step 2:}  We  now consider an arbitrary $X\in\calb$ and  apply the structure result in  \ref{Sperp}. Since the divisible module $X_{\mathcal D}$ belongs to  $\calb\cap\mathcal U^\perp\subset M^\perp$, and $M^\perp$ is closed under extensions, it is enough to show that $X'$ and $\overline{X}$ are in $M^\perp$.
Observe first  that $X'$ and $\overline{X}$ are in $\calb$ since $\calb$ is closed under pure submodules and epimorphic images.
Furthermore, we know from Lemma \ref{pi} and condition (iii) that the torsion-free module
$\overline{X}\in\calb$ is contained in $\calu^\perp$. So, we conclude from Step 1 that $\overline{X}\in
M^\perp$.

Now let us turn to $X'=\bigoplus_{\mu\in\mathfrak{T}}X_\mu$.
Since $\calb\cap\mathcal U^\perp$ is closed under direct sums, we
have $\bigoplus_{\mu\not=\lambda}X_\mu\in\calb\cap\mathcal
U^\perp\subset M^\perp$, so we only have to consider
$X_\lambda$. Recall that there is a pure-exact sequence
$0\to A_\lambda\to X_\lambda\to Z_\lambda\to 0$ where
$A_\lambda$ is a direct sum of modules in
$\tube_\lambda\cap\calb$, and $Z_\lambda\in\mathcal{B}$ is a direct sum of
Pr\"ufer modules belonging to the tube $\tube_\lambda$. Then
$Z_\lambda$ is divisible and therefore in $ M^\perp$ by Step 1,
and $A_\lambda\in M^\perp$ by (ii), thus also $X_\lambda\in
M^\perp$.

  \smallskip

(b) We now prove $M\in\calb$. Let $A\in\cals$, and assume
w.l.o.g.~that $A$ is indecomposable. As in (a) we infer from (i)
that $A\in\mathcal W\cup\mathcal W^o$. If $A\in\mathcal W$, then
$\Ext^1_R(A,M)=0$ because $M$ is an $R_\calu$-module and thus
belongs to $\mathcal W^\wedge$ by Proposition~\ref{univloc}(4). If
$A\in\mathcal W^o$, then we know from Proposition~\ref{univloc}(7)
that $A$ embeds in $A\otimes_R R_\calu\in{}^\perp M$, hence
$A\in{}^\perp M$, and the claim is verified.

  \smallskip

(c) Finally, we check that $T\in \Add (\tube_\lambda(Y)\oplus M)$.
By Theorem~\ref{structure} there is a decomposition
$T=\bigoplus_{\mu\in\mathfrak{T}} \tube_\mu(T)\oplus \overline{T}$
where $\overline{T}$ is torsion-free, and each $\tube_\mu(T)$ is
a direct sum of copies of $\tube_\mu(Y)$ and Pr\"ufer modules
belonging to $\tube_\mu$. Moreover, a Pr\"ufer module $S[\infty]$
occurs as a direct summand in $\tube_\mu(T)$ if and only if
$\tube_\mu\cap\cals$ contains the complete ray $\{S[n]\mid
n\geq1\}$, again by Theorem~\ref{structure}. Observe that complete rays in $\mathcal{S}$
are contained in $\mathcal{W}^o$ by (i). So, we deduce
that the Pr\"ufer modules occurring    as direct summands in
$\tube_\lambda(T)$ do not belong to quasi-simples in $\mathcal U$
and are therefore contained in $\mathcal U^o$, and even in
$\mathcal U^\wedge$ as they are divisible modules. Thus
$\tube_\lambda(T)$ is the direct sum of a module in
$\Add\tube_\lambda(Y)$ with a module in $\mathcal U^\wedge$. Of
course, also the $\tube_\mu(T)$ with $\mu\not=\lambda$ belong to
$\mathcal U^\wedge$. Finally, the torsion-free module
$\overline{T}$ is contained $\mathcal U^o$, and even in $\mathcal U^\wedge$ by Lemma \ref{pi} and condition (iii). So, our claim will be
proven once we show that  $\Add T\cap\mathcal U^\wedge\subset\Add
M$.

Let us thus consider $X\in\Add T\cap\mathcal U^\wedge$. First of
all, $X\in\calb\subset M^\perp$ by (a). Moreover, $X$ is an
$R_\calu$-module, hence $\Ext^1_{R_\calu}(M,X)= \Ext^1_{R}(M,X)=0$.
Therefore $X$ belongs to the $R_\calu$-tilting class $\calb_\calu$,
and there is an  exact sequence $0\to M_1\to M_0\stackrel{f}\to X\to
0$ with $M_0,M_1\in\Add M$ by \cite[5.1.8(d)]{GT}.  Note that $\Add
M\subseteq X^\perp$ because $M\in\mathcal{B}=T^\perp$ and $X\in\Add
T$. Hence  the exact sequence splits and $X\in\Add M$.

 Now the proof of the Proposition is complete.
\EB

In order to specify the set $\calu$ at which we will localize, we  have to distinguish two cases, depending on whether $\tube_\lambda\cap\cals$ contains a complete ray or not.

\Bdef{reducer}
Let $\tube_\lambda$ be a tube of rank $r>1$, and let $S_1[m_1],\ldots,S_l[m_l]$ be the vertices of $T$ in $\tube_\lambda$. We define a  set $\calu$ of quasi-simple modules
 as follows. If $\tube_\lambda\cap\cals$ does not contain a complete ray, then $\calu$ consists of the quasi-simple modules in the union of the wings $\bigcup_{j=1}^l\mathcal{W}_{S_j[m_j]}$. If $\tube_\lambda\cap\cals$ contains a complete ray, then $\calu$ consists of the quasi-simples in $\tube_\lambda$ whose ray is not completely contained in $\cals$.

\medskip

We remark that the set $\calu$
consists of exactly $m_1+\dotsb+m_l<r$ quasi-simple modules. Indeed, this is clear in the first case by Corollary~\ref{disjointwings}. In the second case, the  rays that are not completely contained in $\cals$ correspond to the $m_1+\dotsb+m_l$ quasi-simples in $\bigcup_{j=1}^l\mathcal{W}_{\tau^-S_j[m_j]}$ by Lemma~\ref{lem:completeray}.\Edef

\Bprop{case1} Let $\tube_\lambda$ be a tube of rank $r>1$ such that $\tube_\lambda\cap\cals\neq\emptyset$  does not contain complete rays.
Let $\calu$ be defined as in Definition \ref{reducer}, and
let $M$ be an $R_\mathcal{U}$-tilting module with tilting class
 $\calb_\calu=\{X\in\Modu\,\mid\,\Ext^1_{R_\calu}(A\otimes_R R_\calu, X)=0 \;{\rm for\, all}\; A\in\cals\}$. Then $\tube_\lambda(Y)\oplus M$ is a tilting $R$-module equivalent to $T$ such that neither $\tube_\lambda$ nor
  the $R_\mathcal{U}$-tube $\tube_\lambda\otimes R_\mathcal{U}$ have modules from $\Add M$.
 \Eprop

 \BB
Let $S_1[m_1],\dotsc,S_l[m_l]$ be the vertices of $T$ in
$\tube_\lambda$.
For each $j\in\{1,\dotsc,l\}$, let $\mathcal{U}_j$
consist of the $m_j$ quasi-simples in $\mathcal{W}_{S_j[m_j]}$.
By definition, $\mathcal{U}=\bigcup_{j=1}^l\mathcal{U}_j$. We denote by $\mathcal W$ the extension closure of $\mathcal U$ and recall from Lemma \ref{wing} that $\mathcal W$ consists of all finite direct sums of modules in ${\mathcal W'}=\bigcup_{i=1}^l\mathcal{W}_{S_j[m_j]}$.

We verify conditions (i)-(iii) in Proposition \ref{critilting2}.

(i) We claim $\tube'\subset\mathcal W\cup\mathcal W^o$. Indeed, if
$A\in\tube'\cap \tube_\nu$  with $\nu\neq\lambda$, then clearly
$A\in\calu^o$, which coincides with  $\mathcal W^o$ by
Lemma~\ref{univloc}(3). If $A\in\tube'\cap\tube_\lambda$, then, by
the assumption, $A$ lies on a ray $\{S'[n]\mid n\in\mathbb{N}\}$
which is not completely contained in $\cals$. Then, by
Lemma~\ref{old}(1), $A$ is a submodule of $S'[m]\in\Add T$. By the
definition of vertices, $\mathcal W_{S'[m]}\subseteq \mathcal W$,
and therefore $A\in\mathcal W$.

(ii) Let us now verify  $\Add(\tube_\lambda\cap\calb)\subset
M^\perp$. Choose $A\in \Add(\tube_\lambda\cap\calb)$. By \cite{W}
there is an indecomposable decomposition of the form
$$A=\bigoplus_{p\in P}W_p^{(\alpha_p)} \oplus \bigoplus_{q\in Q}X_q^{(\alpha_q)}$$
where $\{W_p\,\mid\,p\in P\}$  is a complete irredundant set of
modules in ${\mathcal W'}$, and $\{X_q\,\mid\,q\in Q\}$ is a
complete irredundant set of  modules in
$(\tube_\lambda\cap\calb)\setminus {\mathcal W'}$. Note that
the index set $P$ is finite.

First of all, we prove that $\mathcal U\subset M^\perp$.
We fix a $j\in\{1,\dotsc,l\}$, and choose
the numbering $U_{1}=S_j,
U_{2}=\tau^-U_{1},\dotsc,U_{m_j}=\tau^{-}U_{m_j-1}$ for the
quasi-simples in ${\mathcal U}_j$.
For
$1\le i< m_j$ we have $\Ext^1_R(M, U_{i})\cong D{\rm
Hom}_R(U_{i+1},M)=0$ since $M$ is an $R_\calu$-module. Moreover,
$U_{m_j}\in\calb$ because it is a quotient of $S_j[m_j]\in\Add T$,
and $U_{m_j}\in\mathcal U^\perp$ because  $\Ext^1_R(U,
U_{m_j})\cong D{\rm Hom}_R(\tau^-U_{m_j},U)=0$ for all $U\in\calu$ as $\tau^-U_{m_j}\notin\calu$.
Now we infer as in Step 1 of the proof of Proposition
\ref{critilting2} that $\mathcal \calb\cap U^\perp\subset M^\perp$,
so $U_{m_j}\in M^\perp$ as well. Hence $\mathcal U\subset M^\perp$,
thus also ${\mathcal W'}\subset M^\perp$, yielding by
 Lemma \ref{inf}(2) and (3)  that
$\operatornamewithlimits{\bigoplus}\limits_{p\in
P}W_p^{(\alpha_p)}\in M^\perp$.

Next, we consider $X\in(\tube_\lambda\cap\calb)\setminus{\mathcal
W'}$. Note that $S_j[m_j]\in\Add T$ implies $S_j[t]\in\cals$ for all
$t\le m_j$, thus $0=\Ext^1_R(S_j[t], X)\cong D{\rm Hom}_R(X,\tau
S_j[t])$. Then $X$ cannot lie on a coray ending at $\tau U_{1},
U_{1},\ldots,U_{m_j-1}$. Hence $\Ext^1_R(U_{t}, X)\cong D{\rm
Hom}_R(X,\tau U_{t})=0$ for all $1\le t\le m_j$, which shows that
$X\in\mathcal U^\perp$. Therefore
$\operatornamewithlimits{\bigoplus}\limits_{q\in
Q}X_q^{(\alpha_q)}\in \mathcal{B}\cap\calu^\perp$ as $\mathcal{B}$
and $\mathcal{U}^\perp$ are closed under direct sums. Now we infer
as in Step~1 of the proof of Proposition \ref{critilting2} that
$\calb\cap\mathcal U^\perp\subset M^\perp$, and we conclude that
$A\in M^\perp$ as desired.

(iii) Finally, we check that every adic module in $\calb$ belonging
to $\tube_\lambda$ is contained in $\mathcal U^\perp$. So suppose
that  $I=U[-\infty]\in\calb$ for some quasi-simple
$U\in\tube_{\lambda}$. As in (ii), we fix a $j\in\{1,\dotsc,l\}$,
and we see that $0=\Ext^1_R(S_j[t], I)\cong D{\rm Hom}_R(I,\tau
S_j[t])$ for all $1\leq t\le m_j$, hence $U\notin\{\tau U_{1},
U_{1}, \dotsc,U_{m_j-1}\}$. This implies $\Ext^1_R(U_{t}, I)\cong
D{\rm Hom}_R(I,\tau U_{t})=0$ for all $1\le t\le m_j$,
$j=1,\dotsc,l$, that is, $I\in\mathcal U^\perp$.

Therefore $\tube_\lambda(Y)\oplus M$ is a tilting $R$-module equivalent to $T$. Now we prove the remaining assertions.

By Proposition~\ref{useful}(1) and (2), the $R$-tube
$\tube_\lambda$ contains the quasi-simple modules and therefore all modules in the $R_\calu$-tube
$\tube_\lambda\otimes R_\calu$. Moreover, since $\M_\calu$ is a full subcategory of $\M$ closed under direct sums and direct summands,  $\Add _{R_\calu}M=\Add _RM$. So, it is enough to show that $\tube_\lambda$ has
no submodules from $\Add M$.

Assume that $Z\in\tube_\lambda\cap\Add M$. Then $Z$ is an $R_\calu$-module,
because $\M_\calu$ is a full subcategory of $\M$ closed under direct sums and direct summands.
On the other hand, as  $\Add(\tube_\lambda(Y)\oplus M)=\Add T$, we deduce that
$Z$ belongs to $\tube_\lambda\cap\Add T$, thus to
$\mathcal{W}_{S_j[m_j]}$ for some $j\in\{1,\dotsc,l\}$. But then  it
follows from Proposition \ref{univloc}(5) that $Z\otimes_R
R_\calu=0$, a contradiction.
 \EB

\Bprop{case2} Let $\tube_\lambda$ be a tube of rank $r>1$ such that
$\tube_\lambda\cap\cals\not=\tube_\lambda$ contains a complete ray.
Let $\calu$ be as in Definition \ref{reducer}, and let $M$ be an
$R_\calu$-tilting module with tilting class
$\calb_\calu=\{X\in\Modu\,\mid\,\Ext^1_{R_\calu}(A\otimes_R R_\calu,
X)=0 \;{\rm for\, all}\; A\in\cals\}$. Then $\tube_\lambda(Y)\oplus
M$ is a tilting $R$-module equivalent to $T$   such that neither
$\tube_\lambda$ nor
  the $R_\mathcal{U}$-tube $\tube_\lambda\otimes R_\mathcal{U}$ have modules from $\Add M$.
 \Eprop

 \BB
Let $S_1[m_1],\dotsc,S_l[m_l]$ be the vertices of $T$ in
$\tube_\lambda$. For each $j\in\{1,\dotsc,l\}$, let $\calu_j$
consist of the $m_j$ quasi-simples lying in
$\mathcal{W}_{\tau^-S_j[m_j]}$ and choose the numbering
$U_{j1}=\tau^-S_j, U_{j2}=\tau^-U_{j1},\dotsc,
U_{jm_j}=\tau^-U_{jm_j-1}$ for these quasi-simples. By definition,
$\calu=\bigcup_{j=1}^l\calu_j$. We denote by $\mathcal W$ the
extension closure of $\mathcal U$ and recall from Lemma \ref{wing}
that $\mathcal W$ consists of all finite direct sums of modules in
${\mathcal W'}=\bigcup_{i=1}^l\mathcal{W}_{\tau^-S_j[m_j]}$.

Let us verify conditions (i)-(iii) in Proposition \ref{critilting2}.

 (i) We claim $\tube'\subset\mathcal W\cup\mathcal W^o$. If $A\in\tube'\cap\tube_\mu$ with $\mu\neq\lambda$, then clearly $A\in\calu^o$ which coincides with $\mathcal{W}^o$ by Proposition~\ref{univloc}. If $A\in\tube'\cap\tube_\lambda$ and $A$ lies on a ray which is completely contained in $\tube'$, then $A\in\calu^o$ because $\calu$ consists of the quasi-simples in $\tube_\lambda$ whose ray is not completely contained in $\cals$, cf.~Definition~\ref{reducer}. Then, as before, $A\in\mathcal{W}^o$. Assume now that $A\in\tube'\cap\tube_\lambda$ lies on a ray $\{S'[n]\mid n\in\mathbb{N}\}$ which is not completely contained in $\cals$. There exists $t\in\mathbb{N}$ such that $S'[t]\in\Add T$ and $A=S'[v]$ with $v\leq t$ by Lemma~\ref{old}(1). Note that $S'\neq S_j$ because the ray starting at $S_j$ is completely contained in $\cals$, and also that $S'\neq U_{jm_j}$ because $U_{jm_j}\notin\cals$ by Lemma~\ref{old}(2). But, by Proposition~\ref{branchmodules}, $S'[t]\in\mathcal{W}_{S_j[m_j]}$ for some $j$. Hence $S'[v]\in\mathcal{W}_{U_{j1}[m_j-1]}\subseteq\mathcal{W}_{\tau^-S_j[m_j]}$.


(ii) In order to verify $\Add(\tube_\lambda\cap\calb)\subset
M^\perp$, we first observe that no module $\tau S_j[n]$ on the ray
starting at $\tau S_j$ can belong to $\calb$. Consider now a module
$X\in\tube_\lambda\cap\calb$, and assume $X\notin\mathcal U^\perp$.
Then there are $j\in\{1,\dotsc,l\}$  and $i\in\{1,\dotsc,m_j\}$ such
that $0\not=\Ext^1_R(U_{ji}, X)\cong D{\rm Hom}_R(X, \tau U_{ji})$,
hence $X$ lies on one of the corays ending at $S_j=\tau
U_{j1},U_{j1}, \dotsc, U_{jm_j-1}$. But then $X$ must belong to the
wing $\mathcal{W}_{S_j[m_j]}$, because otherwise there is an
epimorphism from $X$ to a module  $\tau S_j[n]$ with $2\le n\le
m_j+1$, which would imply $\tau S_j[n]\in\calb$. Now recall that the
$R_\calu$-module $M$ belongs to $\mathcal U^\wedge$, hence
$\Ext^1_R(M, \tau U_{ji})\cong D{\rm Hom}_R(U_{ji},M)=0$ for all
$1\le i\le m_j$, $1\leq j\leq l$, which proves that
$\{S_j,U_{j1},\dotsc,U_{jm_j-1}\}\subset M^\perp$ and therefore
$X\in \mathcal{W}_{S_j[m_j]}\subset M^\perp$.

As in the proof of Proposition~\ref{critilting2}, step 2, we have
$\calb\cap\calu^\perp\subset M^\perp$, so  the claim follows.

(iii) is trivially satisfied, because the assumption that
$\tube_\lambda\cap\cals$ contains a complete ray implies by Lemma
\ref{pi}  that $\calb$  does not contain any adic module  belonging
to $\tube_\lambda$.

Hence $\tube_\lambda(Y)\oplus M$ is a tilting $R$-module equivalent
to $T$. As in Proposition~\ref{case1}, we observe that there are no
modules in $\tube_\lambda\cap\Add M$. Finally, since $\calu$
consists of the quasi-simples whose ray is not completely contained
in $\cals$, we infer from Proposition~\ref{useful} that the
$R_\calu$-tube $\tube_\lambda\otimes R_\calu$ is completely
contained in $\cals_\calu=\{A\otimes_R R_\calu\,\mid\,A\in\cals\}\subset{}^\perp\calb_\calu\cap\m_\calu$. But then
the tilting $R_\calu$-module $M$ cannot contain direct summands from
this tube by Proposition~\ref{branchmodules}(2). \EB

Now we are in a position to prove our main result.

  \Bth{main}
Let $R$ be a tame hereditary algebra with
$\tube=\bigcup_{\lambda\in\mathfrak{T}}\tube_\lambda$. Every tilting
$R$-module is either equivalent to a finitely generated tilting
module, or to
 precisely one module
$T_{(Y,\Lambda)}$ in the following list:

\begin{enumerate}[(1)]
\item $T_{(Y,\emptyset)}=Y\oplus (\mathbf L\otimes_R R_{\mathcal U})$
where $Y$ is a branch module, and $\calu$ is the set of quasi-simple
composition factors of $Y$.
\item  $T_{(Y,\Lambda)}=Y\oplus R_{\mathcal V}\oplus
R_{\mathcal V}/R_{\mathcal U}$ where  $Y$ is a branch module,
$\emptyset\neq\Lambda\subseteq\mathfrak{T}$, and $\mathcal U$,
$\mathcal V$
are defined as follows:
\begin{enumerate}[(i)]
\item If $\lambda\in\Lambda$, then $\tube_\lambda\cap\calv$
is the complete clique in $\tube_\lambda$, and $\tube_\lambda\cap\calu$
is the set of all the quasi-simples in $\tube_\lambda$ that appear in a regular composition series of $\tau^-Y$.
\item If $\lambda\notin\Lambda$, then
$\tube_\lambda\cap\calv=\tube_\lambda\cap\calu$
consists of all the quasi-simples in $\tube_\lambda$
that appear in a regular composition series of $Y$.
\end{enumerate}
\end{enumerate}
 Moreover,   the  large tilting modules  are  parametrized, up to equivalence, by the elements of $\mathcal Y\times\mathcal P(\mathfrak{T})$, where  $\mathcal P(\mathfrak{T})$ denotes the power set of $\mathfrak T$, and  $\mathcal Y=\{Y_1,\ldots,Y_t\}$ is a complete irredundant set of branch modules over $R$.

\Eth

\BB Let $T$ be a tilting $R$-module with tilting class
$\calb=T^\perp$ and $\cals={}^\perp\calb\cap\m$. Assume that $T$ is
not equivalent to a finitely generated tilting module. Thus there
exists $\tube'\subset \tube$ such that $\cals=\add(\p\cup\tube')$ by
Theorem~\ref{class}. By Lemma~\ref{fginAddT},  there are at most
finitely many non-isomorphic finitely generated indecomposable
modules in $\Add T$ and all of them are regular modules from some
non-homogeneous tube. Let us denote by $Y$ the direct sum of a
complete irredundant set of such modules. By Lemma~\ref{old}(3), $Y$
is a branch module. We want to show that $T$ is equivalent to
$T_{(Y,\Lambda)}$ where $\Lambda=\{\lambda\in\mathfrak{T}\mid
\tube_\lambda\cap\cals \textrm{ contains a complete ray}\}$.

 \underline{Step 1:} Assume that $\Add T$ does not contain finitely
generated  modules. Then $\tube'$ is empty or a union of tubes by
Proposition~\ref{branchmodules}. In the first case,
$\Lambda=\emptyset$,
$Y=0$ and $\calu=\emptyset$, hence $T_{(Y,\Lambda)}=\mathbf L$, which is equivalent to $T$ by   Theorem~\ref{class}.
If $\tube'$ is a union of tubes, then $\Lambda=\{\lambda\mid\tube_\lambda\subseteq \tube'\}$, $Y=0$ and $\calu=\emptyset$, and $\calv$ consists of the quasi-simples in $\tube'$. Hence $T_{(Y,\Lambda)}=R_\calv\oplus R_\calv/R$, which is equivalent to $T$ by Theorem~\ref{class}, as desired.

\underline{Step 2:} Assume now that $\Add T$ contains some finitely generated
indecomposable module.
So, let us consider   a tube $\tube_\lambda$ of rank $r>1$ and such that
$\tube_\lambda\cap\Add T\neq\emptyset$.
Let $\calu_\lambda$ be as in Definition \ref{reducer}. Set
$\calb_{\calu_\lambda}=\{X\in\M_{\calu_\lambda} \mid \Ext_R^1(A\otimes_RR_{\calu_\lambda},X)=0 \textrm{ for all } A\in\cals \}$ and $\cals_{\calu_\lambda}={}^\perp \calb_{\calu_\lambda}\cap\m_{\calu_\lambda}$.

 We have to distinguish two cases depending on whether
$\tube_\lambda\cap\cals$ contains a complete ray or not.

 Suppose first that $\tube_\lambda\cap\cals$ does not contain a complete ray, that is, $\lambda\notin\Lambda$.
Then $\calu_\lambda$ consists of the quasi-simples that appear
in a regular composition series of $\tube_\lambda(Y)$, so
$\calu_\lambda=\tube_\lambda\cap\calu=\tube_\lambda\cap\calv$.

Suppose now that $\tube_\lambda\cap\cals$ contains  a complete ray, that is, $\lambda\in\Lambda$.
Here $\calu_\lambda$ consists of the quasi-simples in $\tube_\lambda$ whose ray is not completely contained in $\cals$, which
coincide with the
 regular composition factors of $\tau^-(\tube_\lambda(Y))$, or in other words, with  the quasi-simples in $\tube_\lambda$ that appear in
the regular series of $\tau^-Y$. Thus $\calu_\lambda=\tube_\lambda\cap\calu$.

Let $M_\lambda$ be a  tilting module over the tame
hereditary algebra  $R_{\calu_\lambda}$ with tilting class $\calb_{\calu_\lambda}$.
It follows from Propositions~\ref{case1} and \ref{case2} that $T$ is equivalent to $\tube_\lambda(Y)\oplus M_\lambda$ over $R$, and $\Add M_\lambda$ has no modules from $\tube_\lambda$.
Over ${R_{\calu_\lambda}}$, we know that $\Add M_\lambda$ has no
module from the $R_{\calu_\lambda}$-tube  $\tube_\lambda\otimes
R_{\calu_\lambda}$, and
further, the modules from the other tubes that belong to $\cals_{\calu_\lambda}$ are the same as before. Indeed,
 if $\mu\neq\lambda$, then $\tube_\mu\otimes
R_{\calu_\lambda}=\tube_\mu$ because every element in $\tube_\mu$ is
already an $R_{\calu_\lambda}$-module. Hence $(\tube_\mu\otimes
R_{\calu_\lambda})\cap\cals_{\calu_\lambda}=\tube_\mu\cap\cals$.

\underline{Step 3:} Now we apply  Propositions \ref{case1} and
\ref{case2} as in Step~2 repeatedly (at most twice more) until we
obtain that $T$ is equivalent to $Y\oplus M$ where $M$ is a tilting
module over  a universal localization $R_{\mathcal U}$ at the set
$\calu$ from Definition \ref{reducer}, and
$\Add_{R_\calu} M$ does not contain finitely generated
 $R_\calu$-modules. Note that ${\mathcal U}$ is a
union of quasi-simples from different tubes, and it does not contain
a complete clique by Proposition~\ref{branchmodules}. Thus $R_\calu$
is a tame hereditary algebra, and  Step~1  yields that $M$ is
equivalent either to the Lukas tilting module over $R_\calu$, or to
a tilting module of the form $(R_\calu)_{\mathcal V'}\oplus
(R_\calu)_{\mathcal V'}/R_\calu$ for  a set $\mathcal V'$ of
quasi-simple $R_\calu$-modules which is a union of cliques over
$R_\calu$.

In the first case we know from \cite[Theorem 6]{aht3} that $M$  is
equivalent to $\mathbf L\otimes_R R_\calu$. Observe that, by
construction,  this first case holds if and only if $\cals$ does not
contain a complete ray, and that $\calu$ is the set of quasi-simples
that appear in the regular composition series of $Y$. Therefore $T$
is equivalent to $T_{(Y,\emptyset)}$.

In the second case we apply Proposition \ref{useful}. By
construction, $\mathcal V'=\{V\otimes R_\calu\,\mid\,V\in{\mathcal
R}\}$ where ${\mathcal R}$ is a set of quasi-simple $R$-modules
defined as follows: if $\lambda\in\Lambda$, then
$\tube_\lambda\cap\mathcal R$  is   the complement of
$\tube_\lambda\cap\mathcal U$, and $\tube_\lambda\cap\mathcal
R=\emptyset$ otherwise. Then $\mathcal V=\calu\cup{\mathcal R}$ and
$(R_\calu)_{\mathcal V'}\cong R_{\mathcal V}$. Thus $T$ is
equivalent to $T_{(Y,\Lambda)}$, as desired.


\underline{Step 4:} Conversely, we show that  for any branch module
$Y$ and any subset $\Lambda\subseteq \mathfrak{T}$, there exists a
tilting $R$-module of the form $T_{(Y,\Lambda)}$ as above.

First of all,  by Lemma~\ref{finitedimensionaltilting}, there exists
a finitely generated tilting $R$-module $H=H_0\oplus Y$ with
$H_0\neq0$ preprojective and
$\cals_H={}^\perp(H^\perp)\cap\m=\add(\p'\cup\tube')$ where
$\p'\subseteq \p$ and $\tube'\subseteq \tube$.

We claim that $\tube'$ does not contain any complete ray.
Indeed,  if $\tube'$ contains a ray, then
we
infer as in Example~\ref{put} that the modules in
 $\tube'^\perp$ cannot have direct summands in  $\p$. But $H^\perp=\cals_H\,^\perp\subset \tube'^\perp$ contains the preprojective module
$H_0\neq0$, a contradiction. Therefore the claim holds true.

Suppose that $\Lambda=\emptyset$. Consider
$\cals=\add(\p\cup\tube')$. Then $\cals$ is a resolving subcategory
of $\m$ because so is $\cals_H$. Hence there exists a tilting
$R$-module $T$ with $\cals={}^\perp(T^\perp)\cap\m$ by
Theorem~\ref{oneone}(1). By Remark~\ref{trivial}, $T$ has neither
preinjective nor preprojective direct summands. Since
there are no finite dimensional regular tilting $R$-modules (cf.~the
proof of Lemma~\ref{finitedimensionaltilting}), we infer that $T$ is a
large tilting $R$-module. By  Steps~1-3 above, $T$ is then
equivalent to a tilting module of the form
$T_{(Y',\emptyset)}=Y'\oplus(\mathbf L\otimes_RR_\calu)$ where
$\calu$ is the set of quasi-simple modules that appear in a regular
composition series of $Y'$. But we know from
Proposition~\ref{branchmodules} that $\Add T\cap\tube_\lambda$ is
determined by $\cals\cap\tube_\lambda$, which coincides with
$\cals_H\cap\tube_\lambda$ for all $\lambda\in\mathfrak{T}$. Hence
$Y\cong Y'$, and  $T_{(Y,\emptyset)}$ is a tilting module equivalent
to $T$.

Suppose now that $\Lambda\neq\emptyset$. By Lemma~\ref{old}(1), the
set $\tube'$ is contained in the union $\mathcal W'$ of the wings
determined by the vertices $S_1[m_1],\ldots,S_l[m_l]$  of $H$. We
now want to enlarge $\tube'$ by inserting some rays from the tubes
$\tube_\lambda$ with $\lambda\in\Lambda$, namely, the rays
corresponding to the set $\mathcal R$  of all quasi-simples  in
$\bigcup_{\lambda\in\Lambda}\tube_\lambda$  that do not appear in
the regular composition series of $\tau^-Y$. So, let
$\tube''\subseteq\tube$ be obtained from $\tube'$ by adding these
rays, that is, $\tube''=\tube'\cup\{S[n]\mid S\in\mathcal{R},\
n\in\N\}$.


We claim that $\add(\p\cup\tube'')$ is a resolving
subcategory of $\m$. To this end, we start by observing that  $\add(\tube'')$ is closed under regular submodules by construction, since so is
 $\add(\tube')$.

Next, we prove that $\add(\tube'')$ is closed under extensions.
Consider an extension $$\exs {R_1'\oplus
R_1''}{\iota}{X}{\pi}{R_2'\oplus R_2''}$$ of two modules in $\add(\tube'')$, where we suppose that
$R_i''$ is a direct sum of indecomposables from the inserted rays  $\{S[n]\mid
S\in\mathcal{R},\ n\in\N\}$,
while $R_i'$  is a direct sum of
indecomposables lying on the remaining rays of $\tube'$.
Let $Z$ be an indecomposable regular direct summand of $X$. The
module $\ker \pi_{\mid_Z}$ is a regular submodule of $Z$, hence
indecomposable. Further, $\ker\pi_{\mid Z}$ is a submodule of
$R_1'\oplus R_1''$, hence a submodule of an indecomposable summand
of $R_1'\oplus R_1''$. If it is a submodule of an indecomposable
summand of $R_1''$,  we are done, because $Z$ then belongs to $\{S[n]\mid
S\in\mathcal{R},\ n\in\N\}\subset\tube''$.
So suppose that $\ker \pi_{\mid Z}$ is a submodule of $R_1'$. Thus
$Z$ is a module lying on a ray starting at a quasi-simple
$S\in\tube'\setminus\mathcal R$. By construction, the kernel of a non-zero map from $Z$ to a module in  $\{S[n]\mid
S\in\mathcal{R},\ n\in\N\}$ cannot belong to the union of the wings $\mathcal W'$. On the other hand, recall that $\ker \pi_{\mid Z}\in\tube'\subset\mathcal W'$. This shows that $\pi{(Z)}$ is a submodule of $R_2'$ and  belongs to $\tube'$. Now $Z$ is an extension
of two elements in $\tube'$, and therefore lies in $\tube'\subseteq
\add \tube''$ as desired.

Finally, we deduce that $\add(\p\cup\tube'')$ is closed under extensions (and is therefore resolving).
Consider an extension $$\exs {P_1\oplus R_1}{\iota}{X}{\pi}
{P_2\oplus R_2}$$ where $P_i\in\add\p$ and $R_i\in \add \tube''$ for
each $i$.
Firstly, $X$ has no preinjective direct summand.
Hence $X=P\oplus R$ with $P\in\add \p$ and $R\in\add \tube$. We
have to prove that $R$ belongs to $\add\tube''$. Observe that
$\pi(R)=R_2'$ is a regular submodule of $R_2$. Thus
$R_2'\in\add\tube''$. Now $\ker \pi_{\mid R}$ is a regular module
 because $\add\tube$ is closed under kernels. Hence it is a submodule of
 $R_1$ and thus $\ker \pi_{\mid R}\in\add
 \tube''$. Therefore $R\in\add\tube''$ because it is the extension of
 two modules in $\add \tube''$.

So,  $\cals=\add(\p\cup\tube'')$ is a resolving subcategory of $\m$.
By Theorem~\ref{oneone}(1), there exists a tilting $R$-module $T$
with ${}^\perp(T^\perp)\cap\m=\cals$, and by the discussion above,
$T$ is equivalent to a tilting $R$-module $T_{(Y',\Lambda)}$ as in
(2). By construction, the vertices of $T$ are exactly the vertices
$S_1[m_1],\ldots,S_l[m_l]$  of $H$ (the only difference being that
the vertices  in the tubes $\tube_\lambda, \lambda\in\Lambda,$ now
lie on rays that are completely contained in $\cals$, while they are
not completely contained in $\cals_H$). Moreover, we know from
Proposition~\ref{branchmodules} that  $\tube_\lambda\cap\Add T$ is
determined by  the  intersections with the corresponding wings
$\cals\cap\mathcal{W}_{S_j[m_j]}$, which coincide with
$\cals_H\cap\mathcal{W}_{S_j[m_j]}$ for all $1\le j\le l$.
 Hence $Y\cong Y'$, and
$T_{(Y,\Lambda)}$ is a tilting module equivalent to $T$, as desired.

\underline{Step 5:} Finally, we establish the parametrization. Observe first that $\mathcal Y$ is indeed a finite set by Lemma~\ref{fginAddT}. Furthermore, we have just seen that the assignment $(Y,\Lambda)\to T_{(Y,\Lambda)}$ is a well-defined  surjective map from $\mathcal Y\times\mathcal P(\mathfrak{T})$ to the set of equivalence classes of tilting modules. It remains to verify the injectivity.
Suppose that $T_{(Y,\Lambda)}$ and  $T_{(Y',\Lambda')}$ are equivalent tilting modules, where
 $Y,Y'\in\mathcal Y$ and $\Lambda,\Lambda'$ are  subsets of
$\mathfrak{T}$.
 Proposition~\ref{Add} then implies that the torsion parts of $T_{(Y,\Lambda)}$ and  $T_{(Y',\Lambda')}$, which are direct sums of  modules with local endomorphism ring, must   coincide up to multiplicity of the summands. We give a  precise description of these  summands in  Remark~\ref{Lukalization} below: for
 $T_{(Y,\Lambda)}$, they are the indecomposable summands of $Y$ and  the Pr\"ufer modules
$S[\infty], {S\in\mathcal R}$, where
 $\mathcal R$ is a set of  quasi-simples  with $\mathcal R\cap\tube_\lambda\not=\emptyset$ if and only if $ \lambda\in\Lambda$, and correspondingly for $T_{(Y',\Lambda')}$.  We conclude that the torsion parts determine $Y,\Lambda$ and $Y',\Lambda'$, respectively, and  we infer that $Y=Y'$ and $\Lambda=\Lambda'$. This concludes the proof of the Theorem.
 \EB

\Brem{Lukalization} (1) Let $Y$ be a branch $R$-module and $\Lambda$
a subset of $\mathfrak{T}$. The tilting $R$-module $T_{(Y,\Lambda)}$
from Theorem~\ref{main} is equivalent to a tilting $R$-module of the
form
$$L_{(Y,\Lambda)}=\bigoplus_{S\in\mathcal R}S[\infty]\oplus Y\oplus (\mathbf{L}\otimes_RR_\calv),$$
where $\calv$ consists of  all quasi-simples in $\bigcup_{\lambda\in\Lambda}\tube_\lambda$  and all the regular composition factors of $Y$, and $\mathcal R$ is the set of  quasi-simples  in $\bigcup_{\lambda\in\Lambda}\tube_\lambda$  that are not  regular composition factors of $\tau^-Y$.  In particular, $\mathcal R\cap\tube_\lambda\not=\emptyset$ if and only if $ \lambda\in\Lambda$.

(2) Let $Z$ be a finitely generated
multiplicity-free regular exceptional  module, and let
$\Delta$ be a set of quasi-simple modules. Set
$$E=\bigoplus_{S\in\Delta}S[\infty]\oplus Z.$$ Then $E$ is a direct summand of a
 large tilting $R$-module $T$
  if
and only if no element of $\Delta$ is a regular composition factor
of $\tau^{-}Z$. In this event,  $T$ is equivalent to
$T_{(Y,\Lambda)}$ where $Y$ is a branch module having $Z$ as a
direct summand, and  $\{\lambda\in\mathfrak{T}\mid
\tube_\lambda\cap\Delta\neq\emptyset\}\subseteq \Lambda$
 \Erem

\BB (1) If $\Lambda=\emptyset$, then $\mathcal{R}=\emptyset$, and
the result holds by Theorem~\ref{main}(1).

Suppose that $\Lambda\neq\emptyset$. Let $\calu,\calv$ be defined as
in Theorem~\ref{main}(2). Then $\calv$ is as stated above, and
$\calv\setminus\calu=\mathcal R$. Moreover, we know from
Propositions~\ref{useful}(2) and \ref{prop:sumofprufer}(1) that
$R_\calv/R_\calu$ is the direct sum of all Pr\"ufer
$R_\calu$-modules corresponding to the  tubes $\tube_\lambda\otimes
R_\calu, {\lambda\in\Lambda}$, which are precisely the  Pr\"ufer
$R$-modules corresponding to the quasi-simples  in
$\bigcup_{\lambda\in\Lambda}\tube_\lambda\setminus\calu$, that is,
to the quasi-simples from $\mathcal R$.
 Hence $\Add
(R_\calv/R_\calu)=\Add
(\bigoplus_{S\in\mathcal R} S[\infty])$. Furthermore, as remarked in Definition \ref{reducer},  the cardinality of $\calu\cap\tube_\lambda$ is always strictly smaller than the rank of $\tube_\lambda$, so $\mathcal R\cap\tube_\lambda\not=\emptyset$ if and only if $ \lambda\in\Lambda$.

Let $\mathbf L$ be the Lukas tilting $R$-module. By (T3), there
exists a short exact sequence
\begin{equation}\label{eq:Lukasshortexact}
\exs R{}{L_0}{}{L_1}
\end{equation}
with $L_0,L_1\in\Add\mathbf L$. By \cite[Lemma~4(iii)]{aht3},
$\Tor_1^R(\mathbf L,R_\calv)=0$. So, applying
${}_{-}\otimes_RR_\calv$ to \eqref{eq:Lukasshortexact}, we obtain
the short exact sequence
\begin{equation}\label{eq:Lukaslocalized}
\exs {R_\calv}{}{L_0\otimes_RR_\calv}{}{L_1\otimes_RR_\calv}.
\end{equation}
By \cite[Theorem~5]{aht3}, $\mathbf{L}\otimes_RR_\calv$ is a
projective $R_\calv$-module, and therefore $L_0\otimes_RR_\calv$ and
$L_1\otimes_RR_\calv$ are projective $R_\calv$-modules. Thus
\eqref{eq:Lukaslocalized} splits, and $\Add_R R_\calv=\Add_R
(\mathbf{L}\otimes_RR_\calv)$. Hence $\Add_R T_{(Y,\Lambda)}=\Add_R
L_{(Y,\Lambda)}$, and $L_{(Y,\Lambda)}$ is a tilting
$R$-module equivalent to $T_{(Y,\Lambda)}$.

(2) Suppose that $E$ is a direct summand of a large tilting
$R$-module $T$. Let $Y$ be a branch module and $\Lambda\subseteq
\mathfrak{T}$ be such that $L_{(Y,\Lambda)}=\bigoplus_{S\in\mathcal
R}S[\infty]\oplus Y\oplus (\mathbf{L}\otimes_RR_\calv)$ is
equivalent to $T$, where $\mathcal R$ is defined as in (1). Set
$\Lambda'=\{\lambda\in\mathfrak{T}\mid
\tube_\lambda\cap\Delta\neq\emptyset\}$. By Proposition~\ref{Add},
$Z$ is a direct summand of $Y$, $\Lambda'\subseteq \Lambda$ and
$\Delta\subseteq \mathcal R$, so no element of $\Delta$ is a regular
composition factor of $\tau^{-}Z$.

Conversely, suppose that no element of $\Delta$ is a regular
composition factor of $\tau^{-}Z$. By
Lemma~\ref{finitedimensionaltilting}, there exists at least one
branch module $Y$ such that $Z$ is a direct summand of $Y$ and no
element of $\Delta$ is a regular composition factor of $\tau^{-}Y$.
For all such $Y$, and for $\Lambda$ containing $\Lambda'$, we get
that $L_{(Y,\Lambda)}$ is a tilting module with $\Delta\subseteq
\mathcal R$. Therefore $E$ is a direct summand of $L_{(Y,\Lambda)}$.
\EB

\Bcor{structure2} Let $R$ be a tame hereditary algebra with
$\tube=\bigcup_{\lambda\in\mathfrak{T}}\tube_\lambda$. Let $T$ be a
large tilting $R$-module, and let
$$T=\bigoplus_{\lambda\in\mathfrak{T}} \tube_\lambda(T)\oplus
\overline{T}$$
be a decomposition
as  in Theorem \ref{structure}.
 Set $\mathcal{B}=T^\perp$ and
$\cals={}^\perp\calb\cap \m$. There are two cases.
\begin{enumerate}[(1)]
\item If $\cals$ contains no complete rays, then
$T$ is equivalent to $T_{(Y,\emptyset)}=Y\oplus ({\mathbf L}\otimes
R_\mathcal{U})$ where $Y$ is a branch module, and $\mathcal U$ is
the set of quasi-simple composition factors of $Y$. Thus $\calu$ is
a set of quasi-simple modules that contains no complete cliques.
Moreover,
\begin{enumerate}[(a)]
\item $\mathcal{F}\cap\mathcal{B}$ consists of the
torsion-free $R_\calu$-modules with no direct summand from
$\p\otimes_RR_\calu$.
\item $\overline{T}$ is equivalent to the Lukas tilting module over $R_{\mathcal U}$.
\end{enumerate}
\item  If $\cals$ contains some rays, then $T$ is equivalent to $T_{(Y,\Lambda)}=Y\oplus R_\calv\oplus
R_\calv/R_\calu$ where $Y$ is a branch module, and $\mathcal
U\subset\mathcal V$ are sets of quasi-simple modules as in
Theorem~\ref{main}(2). Thus $\mathcal{U}$ contains no complete
clique and $\mathcal{V}$ contains complete cliques. Moreover,
\begin{enumerate}[(a)]
\item $\mathcal{F}\cap\calb$ consists of the torsion-free
$R_\mathcal{V}$-modules.
\item $\overline{T}$ is a projective generator for $R_{\mathcal V}$.
\end{enumerate}
\end{enumerate}
 \Ecor

\BB
According  to Theorem~\ref{main}, we see that  $\cals$ contains  no complete rays (respectively, does contain some ray)
if and only if $T$ is equivalent to a tilting module as in (1) (respectively, (2)).

Observe further that, given a subset
$\mathcal{Y}\subset\mathbb{U}$, an
$R_\mathcal{Y}$-module $X$, and a quasi-simple $S$, we have
$$(\ast)\qquad\Hom_R(S,X)\cong\Hom_{R_\mathcal{Y}}(S\otimes_R R_\mathcal{Y},X).$$
In case (1), $R_\calu$ is a tame hereditary algebra  with
preprojective component $\p\otimes_RR_\calu$, and $\{S\otimes_R
R_\mathcal{U}\,\mid\,S\in\mathbb U\setminus{\mathcal U}\}$ is a
complete irredundant set of quasi-simple $R_{\mathcal U}$-modules,
cf.~Proposition~\ref{useful}(2) and \cite[10.1]{GL}. So, $(\ast)$
shows  that an $R_{\mathcal U}$-module is torsion-free over
$R_{\mathcal U}$ if and only if it is torsion-free over $R$.

Now assume that  $X\in\mathcal{F}\cap\mathcal{B}$. Then $X$ is generated by
${\mathbf L}\otimes_RR_\calu$, thus $X\in\calu^\perp$
because the same holds true for the
$R_\calu$-module ${\mathbf L}\otimes R_\calu$. Hence $X$ is an $R_\calu$-module which
is generated by ${\mathbf L}\otimes_RR_\calu$. Since $\mathbf{L}\otimes_RR_\calu$ is equivalent to the Lukas tilting
module over $R_\calu$ \cite[Theorem~6]{aht3}, it follows that $X$ has no direct summand in
$\p\otimes_RR_\calu$.
Conversely, if $X$ is a torsion-free $R_\calu$-module which has no direct summand from
$\p\otimes_RR_\calu$, then it is generated by the Lukas tilting module over $R_\calu$,
whence $X\in\Gen ({\mathbf{L}}\otimes_RR_\calu)\subset{\mathcal B}$.

For  assertion (b), first note that $\overline{T}$ is an
$R_\calu$-module by (a), and $\Ext^1_{R_\calu}(\overline{T},\overline{T}\,^{(I)})=\Ext^1_{R}(\overline{T},\overline{T}\,^{(I)})=0$ for any set $I$. Next observe that $\Add \overline{T}=\Add
(\mathbf{L}\otimes_RR_\calu) $ by
Proposition~\ref{Add}. Since $\M_\calu$ is a full subcategory of
$\M$ closed under direct sums and direct summands, it follows that
$\Add_{R_\calu}
\overline{T}=\Add_{R_\calu}(\mathbf{L}\otimes_RR_\calu)$, and
therefore $\overline{T}$ is a tilting $R_\calu$-module equivalent to
$\mathbf{L}\otimes_RR_\calu$.


We now turn to case (2). Here $R_\mathcal{V}$ is a hereditary order
in $R_\mathbb{U}$ by \cite[4.2]{CB1}, and  $\{S\otimes_R
R_\mathcal{V}\,\mid\,S\in\mathbb U\setminus{\mathcal V}\}$  is a
complete irredundant set of simple $R_{\mathcal V}$-modules,
cf.~Proposition~\ref{useful}(3). Moreover, by definition an
$R_\mathcal{V}$-module $X$ is torsion-free if its torsion submodule
$\{x\in X\,\mid\, xs=0 \,\text{for some regular element}\, s\in
R_{\mathcal V}\}$  is zero, or equivalently, if the canonical map
$X\to X\otimes_RR_{\mathbb{U}}$ is an embedding.

 If
$X\in\mathcal{F}\cap\calb$, then $X\in\Gen R_\mathcal{V}$.
Thus $X\in\mathcal{V}^\perp$ because the same holds true  for $R_\mathcal{V}$. Hence $X$ is an
$R_\mathcal{V}$-module, which is torsion-free as $X\hookrightarrow
X\otimes_RR_{\mathbb{U}}$ by Proposition~\ref{univloc}. For the converse, having $\Gen R_{\mathcal V}\subset \calb$, it is enough to show  that any torsion-free
 $R_{\mathcal V}$-module is also torsion-free over $R$. This is clear in case $\mathcal V=\mathbb U$, so we can assume w.l.o.g.~that $\mathcal V$ is properly contained in $\mathbb U$. Then, as is well known, all simple $R_{\mathcal V}$-modules are torsion, so $(\ast)$ yields the claim, and (a) is verified.

For assertion (b), we show as in case (1)  that $\overline{T}$ is an $R_{\mathcal V}$-module such that $\Add
\overline{T}=\Add
R_\mathcal{V}$.  \EB

We remark that the projective $R_\mathcal{V}$-modules are  well understood, see
for example~\cite{Levysurvey} and \cite[\S 4]{S2}.

\Bex{ray2} Let $\tube_\lambda$ be a tube of rank $r>1$, and let $S$
be a quasi-simple module in $\tube_\lambda$. If $T$ is a tilting
module with $\cals=\add(\p\cup\{\,S[n]\,\mid\,n\in\N\,\})$, then
$$T\sim S\oplus S[2]\oplus\dotsb\oplus S[r-1]\oplus S[\infty]\oplus
R_{\tube_\lambda}.$$

Indeed, we  have already computed $Y=\tube_\lambda(T)= S\oplus
S[2]\oplus\dotsb\oplus S[r-1]$ in Example~\ref{ray}. Choose the
numbering $S=U_1, U_2=\tau^-U_1,\dotsc,U_r=\tau^-U_{r-1}$ for the
quasi-simples in $\mathbf{t}_\lambda$, and  set
$\calu=\{U_2,\ldots,U_r\}$. Consider the universal localization
$R_\calu$. Following the proof of Theorem~\ref{main}, $T\sim S\oplus
S[2]\oplus\dotsb\oplus S[r-1]\oplus M$ where $M$ is a tilting
$R_\calu$-module whose tilting class $\calb_\calu$ is given by
$\calb_\calu=\{X\in\M_\calu \mid \Ext_{R_\calu}^1(A\otimes
R_\calu,X)=0 \textrm{ for all } A\in\cals\}.$ Note that $\{P\otimes
R_\calu\,\mid\,P\in\p\}$ is the preprojective component of $R_\calu$
by \cite[10.1]{GL}, and $\{S[n]\otimes R_\calu\,\mid\,n\in\N\}$ is a
homogeneous $R_\calu$-tube with mouth $\mathcal V'=\{S\otimes
R_\calu\}$. Hence $M\cong (R_\calu)_{\mathcal V'}\oplus
(R_\calu)_{\mathcal V'}/R_\calu\cong R_{\mathcal V}\oplus
R_{\mathcal V}/R_{\mathcal U}$ where $\mathcal
V=\calu\cup\{S\}=\{U_1,\ldots,U_r\}$. We conclude that $R_{\mathcal
V}=R_{\tube_\lambda}$. Moreover, we deduce
as in Remark~\ref{Lukalization} that
$R_{\mathcal V}/R_{\mathcal U}$ is a direct sum of copies of
$S[\infty]$. This proves the claim. \Eex

We now turn to the tilting modules arising from  ring epimorphisms studied in \cite{AS}.

\Bcor{cor:epimoruniversal}
Let $T$ be a tilting $R$-module which is not equivalent to a finite dimensional tilting module. Set $\calb=T^\perp$ and $\cals={}^\perp\calb\cap\m$. The following statements are equivalent.
\begin{enumerate}
\item There exists an injective ring epimorphism $\lambda\colon R\rightarrow R'$ such that $\Tor_1^R(R',R')=0$ and $R'\oplus R'/R$ is a tilting $R$-module equivalent to $T$.
\item $T$
 is equivalent to a tilting module $T_\calu=R_\calu\oplus R_\calu/R$ with $\calu\subseteq\mathbb{U}$.
 \end{enumerate}
  Moreover, under these conditions,  $\cals$ must contain some rays.
\Ecor \BB In \cite[Theorem~6.1]{KS} it is proved that   $\lambda$ as
in (1) can be chosen as a universal localization of $R$. We will
give a different proof for that and   also show that $T$ is
equivalent to $T_\calu$ as stated.

By Proposition~\ref{Add}, both modules   $R',R'/R\in\Add T$  are direct sums of their torsion part and their torsion-free part. We denote by $\overline{R'}$ the torsion-free part of $R'$ and observe that $$\calb=\Gen T=\Gen R'_R=(R'/R)^\perp, \;\text{and in particular,}\; \mathcal{F}\cap\calb\subseteq\Gen {\overline{R'}}.$$

Suppose that $\cals$ contains no complete ray. Then $T$ is equivalent to a tilting $R$-module of the form $Y\oplus(\mathbf{L}\otimes_RR_\calu)$ as in Corollary~\ref{structure2}(1). Since $\mathbf{L}\otimes_RR_\calu\in\mathcal{F}\cap\calb$, we have  $\overline{R'}\neq 0$. Moreover, it follows from Proposition~\ref{Add} that any torsion-free module in $\Add T$ belongs to $\Add (\mathbf{L}\otimes_RR_\calu)$, and any torsion module in $\Add T$ belongs to $\Add Y$.
Now $\mathbf{L}\otimes_RR_\calu$ is equivalent to the Lukas tilting module over the tame hereditary algebra $R_\calu$ by \cite[Theorem~6]{aht3}, and we know  from
\cite[Lemma~3.3(a)]{L2} that $\Hom_{R_\calu}(A,B)\neq0$ for any two nonzero $A,B\in\Add (\mathbf{L}\otimes_RR_\calu)$.
This shows that any torsion-free module  $0\not=A\in\Add T$ satisfies $\Hom_R(A,\overline{R'})\not=0$.
Note that $\Hom_R(R'/R,R')=0$, see for example \cite[2.6]{AS}. So, we infer that $R'/R$ is a torsion module, hence $R'/R\in\Add Y$.
In particular, it follows  that  $(R'/R)^\perp=Y^\perp$. By Lemma~\ref{finitedimensionaltilting}, the branch module $Y$ can be completed to a finite dimensional tilting module $H$ with tilting class $\Gen H=Y^\perp$. But then $\Gen T=(R'/R)^\perp=\Gen H$,   contradicting the assumption that $T$ is not equivalent to a finite-dimensional tilting $R$-module.

\smallskip

So, $\cals$ contains some rays, and $T$ must be equivalent to a tilting module of the form  $Y\oplus R_\mathcal{V}\oplus R_\mathcal{V}/R_\mathcal{U}$ as in  Corollary~\ref{structure2}(2). Since $R_\mathcal{V}\in\mathcal{F}\cap\calb$, we have $\overline{R'}\neq 0$. Moreover, it follows from Proposition~\ref{Add} that any torsion-free module  $A\in\Add T$ belongs to $\Add R_\mathcal{V}$ and is therefore a projective $R_\mathcal{V}$-module. So $0\not=A\in \mathcal{F}\cap\calb \subseteq\Gen{\overline{R'}}$ implies $A\in\Add\overline{R'}$ and
 $\Hom_R(A, \overline{R'})\not=0$. Again, from $\Hom_R(R'/R,R')=0$ we infer that $R'/R$ is a torsion module, hence a direct sum of Pr\"ufer modules and finite-dimensional torsion modules. Observe that if $S[\infty]\in\Add T$ belongs to a tube of rank $r>1$, then it is filtered by $S[r]$, which belongs to $\{S[n]\mid n\geq 1\}\subseteq \cals$ by (the proof of) Theorem~\ref{structure}. Thus $R'/R$ is filtered by  non-projectives in $\cals$, and we can assume that $\lambda$ is a universal localization by \cite[Corollary~3.5]{AA}.

\smallskip

We then know from \cite[2.3]{S2} that $T$ is equivalent to $R_\mathbf{E}\oplus R_\mathbf{E}/R$ for some
 full exact abelian subcategory  $\mathbf{E}$  of $\m$ which is closed under extensions.
By \cite[2.6]{S2} and \cite[4.12 and 4.13]{AS} we have $T^\perp=\mathbf{E}^\perp$, hence
$\cals={}^\perp(\mathbf{E}^\perp)\cap\m$ contains $\mathbf{E}$, and from  the bijection between resolving subcategories of $\m$ and tilting classes given in Theorem~\ref{oneone} we infer that $\cals$ is the resolving closure of $\mathbf E$.
In particular, it follows that $\cals=\add(\p\cup\mathbf{E})$. So, the set $\tube'=\cals\cap\tube$ coincides with $\mathbf E\cap\tube$ and  is therefore closed under cokernels.

\smallskip

We claim that $T$ is equivalent to $T_{\calu'}$ where $\calu'$ is the set of quasi-simple modules in $\tube'$.
 Indeed, $\Gen T=  \tube'^\perp$ as $\tube'$ contains a complete ray (cf.~Example \ref{put}), and by Example
 \ref{Var} it remains to show $\tube'^\perp=(\calu') ^\perp$. Take   $S[m]\in\tube'$. Since $\cals$ is closed under submodules, all $S[n]$ with $n\le m$ are in $\tube'$ as well, and so are the cokernels of the inclusions $S[n]\hookrightarrow S[n']$  for  $n<n'\le m$. Thus $\tube'$  contains the wing $\mathcal{W}_{S[m]}$, and $\calu'$ contains the quasi-simples from that wing. But then  $(\calu')^\perp\subset S[m]^\perp$, and the proof is complete.\EB

 We know from  \cite[Corollary~9]{aht3} that the Lukas tilting module $\mathbf{L}$ is noetherian over its endomorphism ring. The following result generalizes this.

\Bcor{coro:noetherianoverendo}
Let $T$ be a tilting $R$-module which is not equivalent to a finite dimensional tilting module. Set $\calb=T^\perp$ and $\cals={}^\perp\calb\cap\m$. The following statements are equivalent.
\begin{enumerate}
\item
 $T$ is noetherian over its endomorphism ring.
 \item $T$  is equivalent to a tilting module $T_{(Y,\emptyset)}=Y\oplus(\mathbf{L}\otimes_RR_\calu)$ as in Corollary \ref{structure2}(1).
 \item $\cals$ contains no complete rays.
 \end{enumerate}
\Ecor

\BB
We know from \cite[9.9]{AH} that $T$ is noetherian over its endomorphism ring  if and only if $D(T)$ is $\Sigma$-pure-injective. For (2)$\Rightarrow$(1), we  proceed  as in the proof of \cite[Corollary~9]{aht3}. Suppose that $T$ is equivalent to $Y\oplus(\mathbf{L}\otimes_RR_\calu)$. Then $\Add T=\Add (Y\oplus (\mathbf{L}\otimes_RR_\calu))$, and $D(T)\in\Prod (D(Y)\oplus D(\mathbf{L}\otimes_RR_\calu))$. Since $Y$ is finite dimensional, $Y$ and  $D(Y)$ are (right and left, respectively) endofinite modules.
Moreover, by Lemma \ref{firstprop}, the dual of the Lukas tilting $R_\calu$-module $D(\mathbf{L}\otimes_R R_\calu)$ is a cotilting ${R_\calu}$-module whose cotilting class is the class
of $R_\calu$-modules without preinjective summands.
Then $D(\mathbf{L}\otimes_RR_\calu)$ is equivalent to the Reiten-Ringel cotilting module over $R_\calu$
and therefore it is a $\Sigma$-pure-injective $R_\calu$-module. By \cite[1.36]{Facchinibook}, $D(\mathbf{L}\otimes_RR_\calu)$ is also $\Sigma$-pure-injective over $R$.
Hence $D(T)$ is $\Sigma$-pure-injective, and the claim is proven.

For the remaining implications, we show that $T$ is not noetherian over its endomorphism ring
whenever it is equivalent to a tilting module $Y\oplus R_\mathcal{V}\oplus R_\mathcal{V}/R_\calu$ as in Corollary \ref{structure2}(2).
Indeed, in the latter case,  the indecomposable direct summands of the torsion part of $Y\oplus R_\mathcal{V}\oplus R_\mathcal{V}/R_\calu$ are direct summands of $T$, and we see as in Remark~\ref{Lukalization} that $R_\mathcal{V}/R_\calu$ is a non-trivial direct sum of Pr\"ufer modules. Hence $D(T)$ has an adic module as a direct summand. Since adic modules are not $\Sigma$-pure-injective modules, $D(T)$ is not.
\EB


Let us now describe the dual property.

\Bcor{coro:notiltingcotilting}
Let $T$ be a tilting $R$-module which is not equivalent to a finite dimensional tilting module. Set $\calb=T^\perp$ and $\cals={}^\perp\calb\cap\m$. The following statements are equivalent.
\begin{enumerate}[(1)]
\item $T$ is ($\Sigma$-)pure-injective.
\item $G\in\Add T$.
 \item $\cals$ contains a complete ray from each tube.
\item $T$ is equivalent to a tilting module $T_{(Y,\mathfrak{T})}=Y\oplus R_{\mathbb{U}}\oplus
R_{\mathbb{U}}/R_\calu$ where $Y$ is a branch module and $\calu$
consists of the quasi-simples that appear in a regular composition
series of $\tau^-Y$.
 \item $T$ is a cotilting right $R$-module with ${}^\perp T={}^\perp \calb$
 \item There exists a cotilting left $R$-module $C$ such that
 $D(C)$ is a tilting module equivalent to $T$.
 \end{enumerate}
\Ecor

\BB
Recall that $T=Y\oplus Z\oplus\overline{T}$
where $Y$ is a direct sum of copies of finitely many
finite-dimensional modules, $Z$ is a direct sum of Pr\"ufer modules,
and $\overline{T}$ is a non-zero torsion-free module. Now $Y$ is
endofinite, hence $\Sigma$-pure-injective, and $Z\in\Add \mathbf{W}$
is ($\Sigma$-)pure-injective by \cite[10.1]{RR}. So, we have that
$T$ is ($\Sigma$-)pure-injective if and only if so is
$\overline{T}$.

 (1)$\Rightarrow$(2) and (3): By Corollary \ref{structure2}, either $\overline{T}$ is equivalent  to the Lukas
tilting module over the tame hereditary algebra $R_\mathcal{U}$, where $\mathcal U$ is a set of quasi-simple modules that contains no complete cliques, or $\overline{T}$ is a projective generator for $R_{\mathcal V}$, where $\mathcal V$ is a set of
quasi-simple modules that contains  complete cliques.

In the first case,  we know from \cite[Proposition~7 and
Example~8]{aht3} that $\Add \overline{T}$  does not contain
indecomposable pure-injective $R_{\mathcal U}$-modules, and
therefore $\overline{T}$  is not a pure-injective
$R_{\mathcal U}$-module. By \cite[8.62]{JL}, an
$R_{\mathcal U}$-module is  pure-injective over
$R_{\mathcal U}$ if and only if it is  pure-injective over
$R$. So, we conclude that ${T}$ is not pure-injective.

Let us consider   the second case. If $\mathcal V$ is properly contained in $\mathbb U$, then $R_{\mathcal V}$ is a hereditary order in $R_{\mathbb U}$ which is not simple artinian, and from the classification of the indecomposable pure-injective $R_{\mathcal V}$-modules in \cite[3.3]{P} we know that no projective $R_{\mathcal V}$-module can be pure-injective. So,  $\overline{T}$ can only be pure-injective if $\mathcal V=\mathbb U$ and $\overline{T}\in\Add G$.

In particular, we see that $\Add T$ contains $G$, but does not contain  adic modules. On the other hand,  the class $\mathcal F$ of all torsion-free modules  coincides with ${}^\perp(G^\perp)$ by \cite[Proposition~7]{aht3} and is therefore contained in ${}^\perp(T^\perp)$. We infer that all adic modules are in ${}^\perp(T^\perp)\setminus \calb$. By Lemma \ref{pi} it follows that $\cals$ contains a complete ray from each tube.


(3)$\Rightarrow$(4) by Theorem~\ref{main}.

(4)$\Rightarrow$(2): $R_{\mathbb{U}}$ is a direct sum of copies of
$G$ by Proposition~\ref{RU}(2).

(2)$\Rightarrow$(1): If $G\in\Add T$, then $\overline{T}\in G^\perp$ is $\Sigma$-pure-injective by \cite[Proposition~7]{aht3}.

(1)$\Rightarrow$(5): Since $T$ is $\Sigma$-pure-injective, any
module in $\Add T$ is also $\Sigma$-pure-injective, and thus every
pure embedding into a module in $\Add T=\calb\cap{}^\perp\calb$
splits. Hence \cite[Corollary~2.3]{AngeleriSarochTrlifaj} implies
that ${}^\perp\calb$ is closed under direct limits. Now (5) follows
from \cite[Corollary~3.3]{AngeleriSarochTrlifaj}.

(5)$\Rightarrow$(1): Every  cotilting module over an arbitrary ring
is pure-injective by \cite{BK2}.

(3)$\Rightarrow$(6):  By  all the foregoing, we can suppose that $T$
is equivalent to $X=Y \oplus G\oplus \bigoplus_{S\in \mathcal{Y}}
S[\infty]$ where $Y$ is a finite-dimensional module and
$\mathcal{Y}\subseteq \mathbb{U}$. Note that,  for each
$\lambda\in\mathfrak{T}$,
 $\tube_\lambda(X)$ has precisely $r_\lambda$, the rank of
$\tube_\lambda$, pairwise non-isomorphic direct summands.
Observe that $DX$ is a cotilting module. By \cite[3.9]{BK},
(the $r_\lambda$) nonisomorphic direct summands of $DX$
that belong to ${}_R\tube_\lambda$ are precisely the duals
of the nonisomorphic direct summands of
$\tube_\lambda(X)$. Hence, again by \cite[3.9]{BK}, $DX$
is equivalent to $C=DY\oplus G\oplus
\bigoplus_{S\in\mathcal{Y}} (DS)[-\infty]$.

Condition (5) implies that $\Add T=\Prod T$ by
\cite[Corollary~3.3]{AngeleriSarochTrlifaj}. Therefore $DC$ is a
tilting module equivalent to $T$.

(6)$\Rightarrow$(1): The dual of any left $R$-module is a
pure-injective right $R$-module.
 \EB

\medskip



\section{Appendix: The classification of cotilting modules}

 Combining work of Buan and Krause \cite[3.9]{BK} with some combinatorial arguments form \cite{BK2} and  with Bazzoni's result \cite{Bazzonicotilting} stating that every  cotilting module over an arbitrary ring is pure-injective, one obtains a classification of cotilting modules over tame hereditary algebras which we recall below.
We now recover this classification by an elementary proof that only
uses the results from Sections~2--4.

\Bth{BKelementary} Let $R$ be a tame hereditary algebra with
$\tube=\bigcup_{\lambda\in\mathfrak{T}}\tube_\lambda$. Let $C$ be a
cotilting left $R$-module with an indecomposable direct summand
which is not finitely generated. The following  hold true:
\begin{enumerate}[(1)]
\item[(I)] Each indecomposable direct summand of $C$ is either
generic or of the form $S[n]$ for some quasi-simple left
$R$-module $S$ and some $n\in\mathbb{N}\cup\{\infty,-\infty\}$.
\item[(II)]  For each tube ${}_R\tube_\lambda$, $\lambda\in\mathfrak{T}$,
let $\mathcal{I}_\lambda$ be the set of non-isomorphic
indecomposable direct summands of $C$ which are of the form
$S[n]$ for some $n\in\mathbb{N}\cup\{\infty,-\infty\}$ and
quasi-simple $S\in{}_R\tube_\lambda$. Then the number of
elements in $\mathcal{I}_\lambda$ equals the rank of
${}_R\tube_\lambda$.
\end{enumerate}
\Eth

\BB
Let us fix   a cotilting left $R$-module $C$  having an indecomposable direct summand which is not
finitely generated. We
know from Theorem \ref{oneone} that the cotilting class ${}^\perp
C$ is of the form ${}^\perp(\cals^\ast)$ where $$\cals={}^\intercal (^\perp C)\cap\m$$ is a
resolving subcategory of $\m$. Furthermore, if $T$ is a
tilting module with tilting class $\cals^\perp$, then we
know from  Lemma~\ref{firstprop} that $D(T)$ is a cotilting module
equivalent to $C$. More precisely, denoting as before $\calb=T^\perp$, we have $${}^\perp
C={}^\perp(\cals^\ast)=\{{}_RX\,\mid\,D(X)\in\calb\}.$$
Moreover, $$(^\perp C)^\perp\cap\ml=\cals^\ast$$
because every finitely generated left $R$-module $X$ is of the form $X=D(W)$ for some $W\in\m$, and the condition $D(W)\in{}(^\perp C)^\perp$ means by the Ext-Tor-relations that  $W\in{}^\intercal(^\perp C)\cap\m=\cals$.

\medskip

Recall that the modules in $\Prod C=\Prod D(T)$ are pure-injective.
In particular, this implies that ${}^\perp C$
is closed under  direct limits. Since ${}^\perp C$ is also closed under submodules, it follows
 that ${}^\perp
C=\varinjlim ({}^\perp C\cap\ml)$, see \cite[1.1]{BK}. If $I$ is a
pure-injective left $R$-module, we thus have
$$(\sharp)\qquad I\in({}^\perp C)^\perp \text{ if and only if }
\Ext_R^1(A,I)=0 \text{ for all } A\in{}^\perp C\cap\ml.$$

\medskip

\underline{Step 1:} We  compute the indecomposable modules  in $\Prod C$.
First of all, we have
\begin{enumerate}\item[(0)] $\Prod C\cap\ml=\{D(W)\,\mid\,W\in\Add T\cap\m\}$.\end{enumerate}
In fact, if $X$ is a finitely generated left $R$-module of the form $X=D(W)$ with $W\in\m$,
then by the observations above, the condition
$X\in\Prod C={}^\perp C\cap(^\perp C)^\perp$ means that $W\in\cals\cap\calb=\Add T\cap\m$, so (0) is verified.

Recall that there are at most finitely many
non-isomorphic finitely generated indecomposable right $R$-modules in $\Add
T$. As before, we denote by $Y$ the direct sum of  a complete irredundant
set of such modules. Then $D(Y)$ is the direct sum of  a complete irredundant
set of finitely generated indecomposable left $R$-modules in
 $\Prod C$.

\medskip

Next, we compute the adics and the Pr\"ufer modules in $\Prod
C$. Observe that  adic and Pr\"ufer modules are dual to each other. So,  Lemma \ref{pi}, Theorem \ref{structure} and Remark \ref{missing} yield
$$D(T)\cong\prod_{\lambda\in\mathfrak{T}} D(\tube_\lambda(T))\oplus
D(\overline{T})$$ where $D(\overline{T})$ is divisible without finite dimensional direct summands, hence a direct sum of Pr\"ufer modules and copies of $G$, and
$D(\tube_\lambda(T))$ is a direct product of copies of the indecomposable
direct summands of
$D(\tube_\lambda(Y))$ and of adic modules belonging to the corresponding tube
$_R\tube_\lambda$ in $\ml$. More precisely,
 the following statements hold true for a tube $_R\tube_\lambda$ of rank $r$:
\begin{enumerate}
\item[(1)] if
$\cals^\ast$  contains some modules from $_R\tube_\lambda$, but no complete coray, then $D(\tube_\lambda(T))$
is a direct sum of copies of  $s$ pairwise non-isomorphic modules
from $_R\tube_\lambda$, and ${}^\perp C$ contains precisely $r-s$
pairwise non-isomorphic Pr\"ufer modules belonging to
$_R\tube_\lambda$;
\item[(2)] if $\cals^\ast$ contains some   corays
from $_R\tube_\lambda$, then ${}^\perp C$ does not contain any Pr\"ufer module belonging to
$_R\tube_\lambda$, and  $D(\tube_\lambda(T))$ has precisely $r$ pairwise non-isomorphic indecomposable summands: these are  the $s$
adic modules corresponding to the $s\le r$ corays from $_R\tube_\lambda$ contained in $\cals^\ast$,  and  $r-s$ modules from $_R\tube_\lambda$;
\item[(3)] $D(\tube_\lambda(T))=0$ whenever $_R\tube_\lambda\cap\cals^\ast=\emptyset$.
 \end{enumerate}

\medskip

It remains to show that
an indecomposable module belongs to $\Prod C$ if and only if it is isomorphic to a module in the following list:
\begin{enumerate}
\item[-] the indecomposable summands of $D(\tube_\lambda(T))$, $\lambda\in\mathfrak{T}$,
\item[-] the Pr\"ufer modules in ${}^\perp C$,
\item[-] the generic left $R$-module  $_RG$.
\end{enumerate}

\smallskip

For the if-part, we verify that all these modules belong to $\Prod C$. This is clear for the indecomposable summands of $D(\tube_\lambda(T))$, $\lambda\in\mathfrak{T}$.
For the other modules, recall first from Theorem \ref{class} that $\calb\subset\p^\perp$, and $\cals=\add(\p\cup\tube')$ for some subset $\tube'\subset\tube$. Then
 ${}^\perp C\subset {}^\perp{_R\q}$, and $\cals^\ast=\add({_R\q}\cup\tube'')$ for some subset $\tube''\subset {_R\tube}$.
Thus  every $A\in{}^\perp C\cap\ml$ belongs to $_R\p\cup{_R\tube}$, and $\Ext_R^1(A,I)=0$ for any divisible module $I$ without indecomposable preprojective summands. In particular, we deduce from  $(\sharp)$ that all Pr\"ufer modules and the generic module $_RG$ belong to ${}(^\perp C)^\perp$.
Furthermore, since $_RG$ is a  torsion-free module without indecomposable preinjective summands, we also have $_RG\in {}^\perp(\cals^\ast)={}^\perp C$.
 This shows that all modules in our list
belong to $\Prod C$.

For the only-if-part, let $X$ be an indecomposable module in $\Prod C$. Then $X$ is pure-injective, and we can assume w.l.o.g.~that $X$ is neither a Pr\"ufer module nor generic. If $X$ is finite dimensional, then  by (0) it is isomorphic to a finite dimensional  indecomposable summand of $D(\tube_\lambda(T))$ for some $\lambda\in\mathfrak{T}$. If $X=S[-\infty]$ is an adic module, then the class ${}(^\perp C)^\perp$, being closed under epimorphic images, must contain the whole coray ending at $S$. So $X$ is the adic module corresponding to a coray  in $\cals^\ast$, and  by (2) it is isomorphic to an  indecomposable summand of $D(\tube_\lambda(T))$ for some $\lambda\in\mathfrak{T}$. This completes the proof of the claim.

\medskip

\underline{Step 2:} Now statement (I) in the Theorem follows immediately from  Step 1.
For statement (II), we fix a tube ${}_R\tube_\lambda$ of rank $r$ and
let $\mathcal{I}_\lambda$ be the set of non-isomorphic
indecomposable direct summands of $C$ which are adic, Pr\"ufer or finite dimensional modules belonging to $_R\tube_\lambda$. We have to verify that $\mathcal I_\lambda$ has precisely $r$ elements.  By (2) and  (3) in Step 1, we need only  to show that all Pr\"ufer modules in $\Prod C$ belonging to $_R\tube_\lambda$ and all indecomposable summands  of $D(\tube_\lambda(T))$ occur as direct summands of $C$ and therefore are elements of
$\mathcal I_\lambda$. Note that this is clear for the finite dimensional summands, as it is  well known that a finite dimensional indecomposable module arises as a direct summand of a
product of modules $\prod C_j$ if and only if it is a
direct summand of one of the factors $C_j$.
For the other modules, we distinguish two cases.

\medskip

(i) Suppose first that $\cals^\ast$ contains $s>0$   corays
from $_R\tube_\lambda$.
By (2) we  have only to consider  the $s$
adic modules corresponding to these corays. Let $X$ be one of these adic modules.
Then there is a quasi-simple right $R$-module $S\in\tube_\lambda$ such that $X=D(S[\infty])$ and $\cals$ contains the complete ray starting at $S$.
Choose the numbering $S=U_1,\ U_2=\tau^-U_1,\,\dotsc,\
U_r=\tau^-[U_{r-1}]$ for the quasi-simple modules in
$\tube_\lambda$.
Then we have the  numbering $D(U_r)=\tau [D(U_{r-1})],\,\ldots,\, D(U_2)=\tau[D(U_1)],\, D(U_1)=D(S)$  for the quasi-simple modules in
$_R\tube_\lambda$.

Let $m$ be the greatest positive integer such that
$S[m]\in\Add T$, or $m=0$ if $S[m]\notin\Add T$ for all
$m\geq1$. Consider $A= D(U_r[m+1])\in{_R\tube_\lambda}$. Since $S=\tau^- U_r$, we have
$\Ext_R^1(A,X)=\Ext_R^1(S[\infty],U_r[m+1])=D\Hom_R(S[m+1],S[\infty])\neq 0$. But this means that
$A\notin{}^\perp C$ as  $X\in\Prod C$. Now ${}^\perp C=\bigcap_{N} {}^\perp N$ where $N$ runs through all indecomposable direct summands of $C$ by \cite[2.2]{BK}. Thus there must be an indecomposable direct summand $N$ of $C$ with
$\Ext_R^1(A,N)\neq 0$, and of course, $N$ cannot be divisible, nor can it belong to a tube $_R\tube\mu$ with $\mu\neq \lambda$, so  $N$ is a finite dimensional or an adic module belonging to ${_R\tube_\lambda}$.

Note that $D\Hom_R(N,\tau A)\cong \Ext_R^1(A,N)\neq 0$, and
$\tau A\cong D(S[m+1])$ lies on the coray ending at $D(S)$. Moreover,
$U_2[m],U_3[m-1],\dotsc,U_{m+1}\notin\cals$ by Lemma~\ref{old}(2), hence
$D(U_2[m]), D(U_3[m-1]),\dotsc, D(U_{m+1})\notin\cals^\ast$. Since the finite dimensional quotients of $N$ lie in $(^\perp C)^\perp\cap\ml=\cals^\ast$, we deduce that $N$ does neither lie on one of the corays ending at
$D(U_2), D(U_3),\dotsc, D(U_{m+1})$ nor it is an adic module determined by one of these corays. Further, $N$ does not lie on the coray ending at $D(S)$ by choice of $m$. It follows that   $X=D(S)[-\infty]=N$ is the desired direct summand of $C$.

\medskip

(ii)
Suppose now that $\cals^\ast$  contains no complete coray from $_R\tube_\lambda$.
By  (3) we  have only to consider  the $r-s$  Pr\"ufer modules in $\Prod C$ belonging to
$_R\tube_\lambda$.
Let $X=S[\infty]$ be one of these Pr\"ufer left $R$-modules.
Take the greatest positive integer $m$ such that
$S[m]\in\Prod C$, or $m=0$ if $S[m]\notin\Prod T$ for all
$m\geq1$. Then $A= S[m+1]\in{_R\tube_\lambda}$ is cogenerated by $C$, so  there must be an indecomposable direct summand $N$ of $C$ with
$\Hom_R(A,N)\neq 0$. Of course, $N$ cannot be torsion-free, nor can it belong to a tube $_R\tube\mu$ with $\mu\neq \lambda$, so  $N$ is a finite dimensional or a Pr\"ufer module belonging to ${_R\tube_\lambda}$.

Choose the numbering $S=U_1,\ U_2=\tau^-U_1,\,\dotsc,\
U_r=\tau^-[U_{r-1}]$ for the quasi-simple modules in
$_R\tube_\lambda$.
As in Lemma \ref{old}(2), we show that $U_2[m],U_3[m-1],\dotsc,U_{m+1}\notin {}^\perp C$. Since  the finite dimensional submodules  of $N$ lie in $ {}^\perp C$, we deduce that  $N$ does neither lie on one of the rays starting at
$U_2, U_3,\dotsc, U_{m+1}$ nor it is a Pr\"ufer module determined by one these rays.
Further, $N$ does not lie on the ray starting at $S$ by choice of $m$. It follows that   $X=S[\infty]=N$ is the desired direct summand of $C$.
\EB

\Brem{alternative}
Assume that $_R\tube_\lambda$ is a tube of rank $r$ having no   complete coray  in $\cals^\ast$ and having precisely $s\geq0$  non-isomorphic indecomposable modules in $\Prod C$.   As we have seen above, the set $\mathcal I_\lambda$ contains $r-s$
Pr\"ufer modules. They arise as duals of the $r-s$ adic modules in $\calb$ established by Lemma \ref{pi}(2), see also Remark \ref{missing}.
 We will now give an alternative explanation for the occurrence of these Pr\"ufer modules by using Proposition \ref{useful}.

 \medskip

 Let the notation be as above.
According to Theorem~\ref{main} and Corollary~\ref{structure2},  we  distinguish two cases.

\smallskip

(1)  $\cals^\ast$ contains no complete corays. Then,  up to equivalence, $T=Y\oplus (\mathbf{L}\otimes_RR_\calu)$  as in Corollary \ref{structure2}(1).
By \cite[Theorem~6]{aht3}, $\mathbf{L}\otimes_RR_\calu$ is equivalent to the Lukas tilting module over $R_\calu$,
so $D(\mathbf{L}\otimes_RR_\calu)$ is a cotilting $R_\calu$-module equivalent to the Reiten-Ringel tilting $R_\calu$-module $\mathbf{W}_\calu$. Hence $\Prod_{R_\calu}{D(\mathbf{L}\otimes_RR_\calu)}=\Prod_{R_\calu} {\mathbf{W}_\calu}=\Add_{R_\calu} \mathbf{W}_\calu$, and  $\Prod D(T)=\Add (D(Y)\oplus \mathbf{W}_\calu)$. Therefore any module in
$\Prod C$ is a direct sum of indcomposable direct summands of $D(Y)$ and of Pr\"ufer $R_\calu$-modules.

By assumption, there are precisely $s\geq0$  non-isomorphic indecomposable modules in $\tube_\lambda\cap \Add T$    (in fact, in $\Add \tube_\lambda(Y)$),  whose duals give the indecomposables in  $\ltube_\lambda\cap\Prod C$. By construction, $\calu\cap\tube_\lambda$ has $s$ elements. Hence the $R_\calu$-tube  $\tube_\lambda\otimes R_\calu$ has $r-s$ quasi-simples, and $\Prod C$ has precisely $r-s$ Pr\"ufer left $R_\calu$-modules belonging to this tube.

\medskip

(2) $\cals^\ast$ contains some corays. Then,  up to equivalence,
 $T=Y\oplus R_\calv/R_\calu\oplus R_\calv$  as in Corollary \ref{structure2}(2).
Thus $\Prod C=\Prod (D(Y)\oplus D(R_\calv/R_\calu)\oplus D(R_\calv))$, and the Pr\"ufer modules in $\Prod C$ are all in $\Prod D(R_\calv)$ because there are no nonzero morphism from a Pr\"ufer module neither to a torsion-free module nor to a regular module.

By assumption, $\tube_\lambda\cap\cals$ does not contain a complete ray, and according  to the construction, $\calv$ cannot
contain all  quasi-simple $R$-modules in $\tube_\lambda$.
More precisely, $\tube_\lambda(T)$ has $s$ pairwise non-isomorphic indecomposable direct summands, whose duals give the indecomposables in  $\ltube_\lambda\cap\Prod C$. They are arranged in disjoint wings from $\tube_\lambda$, and
the quasi-simple modules in $\calv\cap\tube_\lambda$ are
precisely the quasi-simples in these wings. So, there are
exactly $s$ quasi-simple modules in
$\calv\cap\tube_\lambda$. Each of the remaining $r-s$ quasi-simple modules
$S\in\tube_\lambda\setminus\calv$ gives rise to a
simple $R_\calv$-module $S\otimes_RR_\calv$ with projective
presentation  $0\to\mathfrak{m}\to
R_\calv\to S\otimes_RR_\calv\to 0$ for some   maximal right
ideal   $\mathfrak{m}$. Applying $D$, we
obtain the exact sequence $0\to D(S\otimes_RR_\calv)\to
D(R_\calv)\to D(\mathfrak{m})\to 0$.

Observe that $D(R_\calv)$ is an injective left
$R_\calv$-module \cite[Corollary~3.6C]{Lam2} that contains
the simple left $R_\calv$-module $D(S\otimes_RR_\calv)$. Thus the injective
envelope $_RS[\infty]$ of $D(S\otimes_RR_\calv)$ is a direct summand of
$D(R_\calv)$. We conclude that $\Prod C$ has precisely $r-s$ Pr\"ufer left $R_\calv$-modules belonging to this tube.\Erem


\def\cprime{$'$} \def\cprime{$'$} \def\cprime{$'$} \def\cprime{$'$}
  \def\cprime{$'$} \def\cprime{$'$} \def\acento{\'a}

  \renewcommand{\baselinestretch}{1}



\end{document}